\documentclass[leqno]{article}

\usepackage[english]{babel}
\usepackage[utf8]{inputenc}
\usepackage{amsmath}
\usepackage{amssymb}
\usepackage{amsfonts}
\usepackage{enumerate}
\usepackage{vmargin}
\usepackage[all]{xy}
\usepackage{mathrsfs}
\usepackage{mathtools}
\usepackage{lmodern}
\usepackage{slashed}
\usepackage{xurl}
\usepackage[colorlinks=true,linkcolor=blue,pagebackref=true,hyperfootnotes=false]{hyperref}%
\usepackage{array}
\setmarginsrb{3cm}{3cm}{3.5cm}{3cm}{0cm}{0cm}{1.5cm}{3cm}
\usepackage{comment}

\IfFileExists{accents.sty}{\usepackage{accents}}{}
\newlength{\dhatheight}

\newcommand{\Cbb}{\ensuremath{\mathbb{C}}}

\newcommand{\N}{\ensuremath{\mathbb{N}}}
\let\L\relax %\L ne fait plus rien !
\newcommand{\L}{\mathrm{L}}
\newcommand{\Q}{\mathrm{Q}}
\newcommand{\M}{\mathrm{M}}
\let\cal\relax
\newcommand{\cal}{\mathcal}

\newcommand{\Z}{\ensuremath{\mathbb{Z}}}

\newcommand{\R}{\ensuremath{\mathbb{R}}}
\newcommand{\T}{\ensuremath{\mathbb{T}}}

\renewcommand{\leq}{\ensuremath{\leqslant}}
\renewcommand{\geq}{\ensuremath{\geqslant}}
\newcommand{\qed}{\hfill \vrule height6pt  width6pt depth0pt}

\newcommand{\co}{\colon}

\newcommand{\ot}{\otimes}
\newcommand{\ovl}{\overline}

\let\i\relax %\i ne fait plus rien !
\newcommand{\i}{\mathrm{i}}
\newcommand{\ov}{\overset}

\newcommand{\epsi}{\varepsilon}
\renewcommand{\d}{\mathop{}\mathopen{}\mathrm{d}} %op\'erateur diff\'erentiel
\newcommand{\e}{\mathrm{e}} %constante e
\renewcommand{\d}{\mathop{}\mathopen{}\mathrm{d}}

\DeclareMathOperator{\Irr}{Irr} %Irreductible
\DeclareMathOperator{\Rep}{Rep} %Irreductible
\newtheorem{thm}{Theorem}[section]

\newtheorem{prop}[thm]{Proposition}

\newtheorem{cor}[thm]{Corollary}
\newtheorem{lemma}[thm]{Lemma}

\newtheorem{remark}[thm]{Remark}

\newenvironment{proof}[1][]{\noindent {\it Proof #1} : }{\hbox{~}\qed
\smallskip
}

\usepackage{tocloft}
\numberwithin{equation}{section}
\usepackage[nottoc,notlot,notlof]{tocbibind}%biblio dans sommaire
\let\OLDthebibliography\thebibliography
\renewcommand\thebibliography[1]{
  \OLDthebibliography{#1}
  \setlength{\parskip}{0pt}
  \setlength{\itemsep}{0pt plus 0.3ex}
}

\newcommand\reallywidehat[1]{\arraycolsep=0pt\relax%
\begin{array}{c}
\stretchto{
  \scaleto{
    \scalerel*[\widthof{\ensuremath{#1}}]{\kern-.5pt\bigwedge\kern-.5pt}
    {\rule[-\textheight/2]{1ex}{\textheight}} %WIDTH-LIMITED BIG WEDGE
  }{\textheight} % 
}{0.5ex}\\           % THIS SQUEEZES THE WEDGE TO 0.5ex HEIGHT
#1\\                 % THIS STACKS THE WEDGE ATOP THE ARGUMENT
\rule{-1ex}{0ex}
\end{array}
}

\begin{document}
\selectlanguage{english}
\title{\bfseries{Haagerup symmetries of the $E_8$ theory and conformal nets at central charge eight}}
\date{}
\author{\bfseries{C\'edric Arhancet}}
\maketitle

%%%%%%%%%%%%%%%%%%%%%%%%%%%%%%%%%%%%%%%%%%%%%%%%%%%%%%%%%%%%%%%%
%%%%%%%%%%%%%%%%%%%%%%%%%%%%%%%%%%%%%%%%%%%%%%%%%%%%%%%%%%%%%%%%

\begin{abstract}
A Haagerup algebra is a positive strongly rational vertex operator algebra $U$ of central charge eight with $\Rep(U)$ braided equivalent to the Drinfeld center $Z(\cal H_3)$ of the Haagerup fusion category. We prove that $U$ has exactly three holomorphic extensions, each isomorphic to $V_{E_8}$, and that the corresponding Lagrangian algebras have module categories $\cal H_1$, $\cal H_2$ and $\cal H_3$. The extensions attached to $\cal H_2$ and $\cal H_3$ are the two lattice extensions of $V_{E_6 \oplus A_2}$. We prove a recognition theorem for $\cal H_2$ and $\cal H_3$ among unitary fusion categories of rank six with three invertible objects, a noncommutative Grothendieck ring and trivial pointed associator. For simple regular unitary vertex operator algebras with bounded localized vertex operators, we construct the transport module of two ordinary unitary modules when the second integrates in Tener's annular sense. For a unitary conformal Haagerup subalgebra of $V_{E_6 \oplus A_2}$, these results give a completely rational net with representation category $Z(\cal H_3)$, a realization of the Haagerup subfactor and the generalized orbifold of Bischoff. In particular, they apply to the strongly local explicit algebra $S$ of levels $(3,13)$, whose Haagerup property is established in \cite[Theorem 6.1]{ArhIII}. In \cite{ArhIII}, the net of $S$ is recognized from strong rationality and lattice branching, and a diagonal extension of the net is then used to identify $\Rep(S)$. Here an exact construction of the compatible unitary structure on the diagonal lattice algebra $V_{E_8 \oplus E_8}$, followed by passage to the twin algebra, proves complete unitarity and Condition I+ for every Haagerup algebra. This also gives the VOA--net comparison through smeared vertex operators for unitary conformal embeddings in $V_{E_6 \oplus A_2}$.
\end{abstract}

%%%%%%%%%%%%%%%%%%%%%%%%%%%%%%%%%%%%%%%%%%%%%%%%%%%%%%%%%%%%%%%%
%%%%%%%%%%%%%%%%%%%%%%%%%%%%%%%%%%%%%%%%%%%%%%%%%%%%%%%%%%%%%%%%

\makeatletter
 \renewcommand{\@makefntext}[1]{#1}
 \makeatother

{
  \hypersetup{linkcolor=blue}
 \tableofcontents
}

%%%%%%%%%%%%%%%%%%%%%%%%%%%%%%%%%%%%%%%%%%%%%%%%%%%%%%%%%%%%%%%%%%%%%%%%%%%%%
\section{Introduction}
\label{sec:Introduction}

From a completely rational conformal net one obtains finite index subfactors, and the representation category of the net is a unitary modular tensor category \cite{KLM01}. The connection with finite depth subfactors goes back to Jones \cite[Section 1.3 p.~2]{Bro20}. In the formulation of \cite[Definition 4.1 p.~20]{Bis17a}, a subfactor $N \subset M$ arises from a completely rational net $\cal A$ if it is of the form $\beta(\cal B_a(I)) \subset \cal B_b(I)$ for two irreducible relatively local extensions $\cal B_a, \cal B_b \supset \cal A$ and a sector $\beta$ contained in $\iota_b\rho\ovl{\iota}_a$ for some $\rho \in \Rep(\cal A)$, where $\iota_a$ and $\iota_b$ are the inclusions $\cal A(I) \subset \cal B_a(I)$ and $\cal A(I) \subset \cal B_b(I)$. He proved that $N \subset M$ arises from $\cal A$ as soon as $\Rep(\cal A)$ is braided equivalent to the quantum double $D(N \subset M)$ \cite[Proposition 4.7 p.~22]{Bis17a}, and he formulated, following \cite{EvG11}, the problem of realizing such quantum doubles by conformal nets \cite[Problem 1.3 p.~16]{Bis17a}. For all subfactors of index less than four this problem was solved in \cite{Bis15}. The Haagerup subfactor of index $\frac{5+\sqrt{13}}{2}$ \cite{AsH99}, which has the smallest index above four among finite depth subfactors \cite[Section 4.2.1 p.~22]{JMS14}, is the first case of index above four. Its quantum double is the Drinfeld center $Z(\cal H_3)$ of any of the three Morita equivalent Haagerup fusion categories $\cal H_1$, $\cal H_2$, $\cal H_3$ of Grossman and Snyder \cite{GrS12}.

Evans and Gannon predicted that $Z(\cal H_3)$ is realized by a chiral theory of central charge eight \cite{EvG11}, and numerical work suggests a Haagerup-symmetric theory with $c \approx 2$ \cite{HL+22,V+22}. At central charge eight, explicit vertex operator algebras inside $V_{E_6 \oplus A_2}$ have recently been proposed by Jiao \cite{Jia26}, Riesen \cite{Rie26} and Huang \cite{Hua26}. Jiao gives a computer-assisted proof that his algebra realizes $Z(\cal H_3)$, in a preprint dated 4 October 2026, and Riesen and Huang announce the same property for their algebras, with outlines of proof. In \cite{Arh26} we placed every Haagerup algebra in the chain $V_{E_6 \oplus A_2} \subset V_{E_8}$, determined its currents and constructed its twin. There are exactly two isomorphism classes as soon as one of the explicit candidates of \cite{Jia26,Rie26} is a Haagerup algebra \cite[Theorems 12.3 and 12.4]{Arh26}. In \cite[Theorem 6.1]{ArhIII} we prove that the case (II) algebra $S$, of levels $(3,13)$, is a Haagerup algebra, and in \cite[Theorem 7.13(1)]{ArhIII} we identify its case (I) twin, of levels $(1,39)$, with the algebras of \cite{Rie26,Hua26}. Albert, Honda, Kaidi and Zheng suggested that the chiral $(\mathfrak e_8)_1$ theory carries the symmetries $\cal H_1$, $\cal H_2$ and $\cal H_3$, and listed the three Lagrangian algebras of $Z(\cal H_3)$ through the tube algebra \cite[Appendix A pp.~6--7]{AHKZ26}, see also \cite{BoS25}. Huang announced, with an outline of proof, that his algebra is strongly rational and completely unitary with $\Rep(U) \cong Z(\cal H_3)$ as unitary ribbon categories, that its conformal net is completely rational, and that $\cal H_3$ acts by topological defects on a holomorphic extension isomorphic to $V_{E_8}$ \cite[Theorem 1 p.~1 and Sections 6 and 7 pp.~8--10]{Hua26}.

We determine the three Lagrangian algebras of $Z(\cal H_3)$ and identify their module categories intrinsically. We realize them by the holomorphic extensions of an arbitrary Haagerup algebra. For unitary conformal subalgebras of $V_{E_6 \oplus A_2}$ we construct the associated nets, and for those which are Haagerup algebras we prove complete rationality and identify the net category. The finite-index step uses a general construction of unitary transport modules. The identification uses a recognition theorem for unitary fusion categories. We obtain the realization of the Haagerup subfactor and the generalized orbifold of Bischoff. Finally, the exact lattice computation in Section \ref{subsec:complete-unitarity} gives a compatible unitary structure on a holomorphic diagonal extension for any Haagerup algebra of case (II). The coset theorems of Gui and the twin construction then give the complete unitarity theorem for arbitrary Haagerup algebras (Theorem \ref{rem:all-haagerup}).

Throughout, $d = \frac{3+\sqrt{13}}{2}$, so that $d^2 = 3d+1$. We use the labels $0, b, a, c_0, c_1, c_2, d_1, \ldots, d_6$ of \cite{EvG11} for the simple objects of $Z(\cal H_3)$, whose quantum dimensions are $1$, $3d+1$, $3d+2$ (four times) and $3d$ (six times). For a Haagerup algebra $U$ we denote by $M_z$ the simple $U$-module attached to $z$, and we put $V' = V_{E_6 \oplus A_2}$. By \cite[Proposition 4.1]{Arh26}, $U$ is contained in $V'$ and $V' = U \oplus M_b$ as a $U$-module.

\begin{thm}
\label{thm:intro-A}
Let $U$ be a Haagerup algebra.
\begin{enumerate}
\item Up to isomorphism of extensions, $U$ has exactly three holomorphic extensions $B_1$, $B_2$, $B_3$, that is, simple vertex operator algebras $B$ containing $U$ with the same conformal vector and with $\dim B_0 = 1$, every $B$-module of which is a direct sum of copies of $B$. Each $B_i$ is isomorphic to $V_{E_8}$. The corresponding Lagrangian algebras $L_i$ of $\Rep(U)$ satisfy $\Rep(U)_{L_i} \cong \cal H_i$, and
\[
L_1 \cong 0 \oplus b \oplus a \oplus c_0,\quad L_2 \cong 0 \oplus b \oplus 2a,\quad L_3 \cong 0 \oplus b \oplus 2c_0
\]
as objects, in the labels of \cite{EvG11}.
\item The extensions $B_2$ and $B_3$ are the two lattice vertex operator algebras $V' \oplus V'_{w} \oplus V'_{-w}$ attached to the two isotropic lines $\{0,\pm w\}$ of the discriminant group $\Z_3 \times \Z_3$ of $E_6 \oplus A_2$. Here $V'_w$ denotes the simple $V'$-module attached to $w$, and $B_2$ is the extension for which the $U$-modules $V'_{\pm w}$ are isomorphic to $M_a$.
\item The subspace $U \oplus M_b$ of $B_1$ is not a vertex subalgebra, and $U$ and $B_1$ are the only vertex subalgebras of $B_1$ containing $U$. For $i \in \{2,3\}$, the simple vertex subalgebras of $B_i$ containing $U$ are $U$, $V'$ and $B_i$.
\end{enumerate}
\end{thm}

This is Theorem \ref{thm:lagrangian-extensions} with Propositions \ref{prop:three-lagrangian} and \ref{prop:which-line}. The Lagrangian algebras attached to the two even parts $\cal H_2$ and $\cal H_1$ of the Haagerup subfactor were computed by Evans and Gannon in the language of Longo--Rehren inclusions \cite[pp.~18--19]{EvG11}, the second one by a Frobenius reciprocity argument of which Lemma \ref{lem:rank} is a categorical version. The three Lagrangian algebras, with the categories $\cal H_1$, $\cal H_2$ and $\cal H_3$ attached to them, are listed in \cite[(A8) p.~7 and pp.~3--4]{AHKZ26}, in labels $1$, $\pi_1$, $\pi_2$, $\sigma_1$ coming from the tube algebra of $\cal H_3$. We prove this identification in the labels of \cite{EvG11}, which are defined through the half-braidings on the even part $\cal H_2$ of the Haagerup subfactor, by means of the classification of module categories \cite{GrS12} and of a formula of Ostrik for the dimensions of the spaces of morphisms between Lagrangian algebras \cite[Remark 2.19 pp.~7--8]{Ost15}. This proof does not use vertex algebras (Corollary \ref{cor:categorical}), and it gives $\pi_2 = c_0$ and $\sigma_1 = a$ (Remark \ref{rem:which-line}). The realization of these Lagrangian algebras by vertex operator algebras uses the restriction of the lattice modules computed in \cite[Lemma 4.3]{Arh26}. For the algebra of \cite{Hua26}, which is a Haagerup algebra by \cite[Theorem 7.13(1)]{ArhIII}, part 2 shows that the two types of lattice extension discussed in \cite[Section 7 p.~10]{Hua26} are realized by the two isotropic lines, and that the extension attached in \cite[Section 7 pp.~9--10]{Hua26} to the algebra $I_{\cal H_3}(1)$ is the lattice extension of type $c_0$. This agrees with \cite[Remark 2.1]{Arh26}. Part 3 shows that there is no vertex subalgebra strictly between $U$ and the extension $B_1$ attached to $\cal H_1$.

We now turn to conformal nets. Let $S \subset V'$ be the vertex algebra generated by the weight-one and weight-two spaces of \cite[Theorem 1.3]{Arh26}. It is unitary for an explicit unitary structure of $V'$ \cite[Theorem 10.10]{Arh26}, and it is isomorphic to the algebra of Jiao \cite[Theorem 10.24]{Arh26}. We write $\mathrm{(H_J)}$ for the assertion that Jiao's algebra is a Haagerup algebra \cite[Section 1]{Arh26}, as stated in \cite{Jia26}, and $\mathrm{(H_U)}$ for the assertion that every simple $S$-module admits a scalar product which is invariant with respect to the PCT operator of $S$ induced by the scalar product of $V'$. Under $\mathrm{(H_J)}$, all the simple $S$-modules except the six modules $M_{d_l}$ are submodules of unitary $V'$-modules \cite[Proposition 4.1, Lemma 4.3 and Lemma 11.8]{Arh26}, so that $\mathrm{(H_U)}$ concerns these six modules. We prove $\mathrm{(H_U)}$ under $\mathrm{(H_J)}$ (Theorem \ref{thm:intro-U}). We write $\mathrm{(H_R)}$ for the analogue of $\mathrm{(H_J)}$ for the algebra of \cite{Rie26}. Let $\cal A_{E_8}$ be the conformal net of $V_{E_8}$, that is, the loop group net of $E_8$ at level one. For an $S$-submodule $M$ of a lattice extension of $V'$, we denote by $\pi_M$ the representation of the conformal net $\cal A_S$ of $S$ on the closure of $M$ (Lemma \ref{lem:submodules}), and we put $\pi_0 = \pi_S$, $\pi_b = \pi_{M_b}$, $\pi_a = \pi_{M_a}$, $\pi_{c_0} = \pi_{M_{c_0}}$, where $M_a$ and $M_{c_0}$ are taken in the two lattice extensions.

\begin{thm}
\label{thm:intro-B}
\begin{enumerate}
\item The unitary vertex operator algebra $S$ is strongly local. Its conformal net $\cal A_S$ is a covariant subnet of the conformal net $\cal A_{V'}$ of $V'$, which is itself a covariant subnet of index three of the net $\cal A_{H_\Lambda} \cong \cal A_{E_8}$ of each of the two lattice extensions $H_\Lambda = V' \oplus V'_w \oplus V'_{-w}$ attached to the isotropic lines $\Lambda = \{0,\pm w\}$ of the discriminant group of $E_6 \oplus A_2$. These inclusions are irreducible, that is, $\cal A_S(I)' \cap \cal A_{E_8}(I) = \Cbb$ and $\cal A_{V'}(I)' \cap \cal A_{E_8}(I) = \Cbb$ for any interval $I$.
\item Assume $\mathrm{(H_J)}$. Then $\cal A_S$ is completely rational. The vacuum representation of $\cal A_{E_8}$ restricts to $\cal A_S$ as $\pi_0 \oplus \pi_b \oplus 2\pi_x$, where $x = a$ for one lattice extension and $x = c_0$ for the other, and $\pi_0$, $\pi_b$, $\pi_a$, $\pi_{c_0}$ are pairwise inequivalent irreducible representations with finite statistical dimensions $1$, $\delta$, $1+\delta$, $1+\delta$ for some $\delta > 1$. The indices are $[\cal A_{V'} : \cal A_S] = 1+\delta$ and $[\cal A_{E_8} : \cal A_S] = 3(1+\delta)$. Finally, $\Rep(\cal A_S)$ is braided equivalent to the Drinfeld center of a unitary fusion category of rank six with a pointed fusion subcategory of dimension three.
\end{enumerate}
\end{thm}

This is Theorems \ref{thm:strong-locality} and \ref{thm:complete-rationality}, Proposition \ref{prop:delta-big} and Corollary \ref{cor:S-rational}. Theorem \ref{thm:strong-locality} applies to any unitary conformal subalgebra of $V'$, and Theorem \ref{thm:complete-rationality} to any such subalgebra which is a Haagerup algebra. The unitary structure on $S$ used in part 1 is that of \cite[Theorem 10.10]{Arh26}. For part 2, the simple summands $M$ of the lattice extensions and their conjugates are unitary and have annular representations. The next theorem constructs a transport module on $\ovl M \boxtimes M$, and the creation operators are localized by restriction from the lattice extensions. Tener's finite-index theorem \cite[Theorem 4.19]{Ten24} then applies to these sectors. The proof of Theorem \ref{thm:complete-rationality} passes from the finite-index sectors to the subnet inclusion. It obtains complete rationality by \cite[Theorem 24 p.~20]{Lon03} and identifies Tener's net with $\cal A_S$ by \cite[Lemma 7.6 p.~72]{Ten19a}.

\begin{thm}
\label{thm:intro-U}
Let $V$ be a simple regular unitary vertex operator algebra with bounded localized vertex operators, and let $M$ and $N$ be ordinary unitary $V$-modules such that $N$ integrates to a representation $\pi^N$ of the conformal net of $V$ in the sense of \cite[Definition 2.18 p.~20]{Ten24}. Then the fusion product $M \boxtimes N$ admits a unitary structure for which it is a transport module for $M$ and $N$. Consequently, under $\mathrm{(H_J)}$, every simple $S$-module is unitarizable, that is, $\mathrm{(H_U)}$ holds.
\end{thm}

This is Proposition \ref{prop:fusion-unitary} and Corollary \ref{cor:modules-unitary}. Tener proved that, when $\pi^N$ exists, the transport form is positive semidefinite near the unit circle \cite[Proposition 4.2 p.~43]{Ten24}, and that, if $M \boxtimes N$ admits a unitary structure, it is a transport module as soon as the transport form is positive semidefinite at one point $(z,z)$ of $\cal R$ \cite[Proposition 3.16 p.~36]{Ten24}, which gives \cite[Theorem 4.4 and Corollary 4.5 p.~44]{Ten24} under this hypothesis. He noted that the unitarizability of $M \boxtimes N$ is open in general \cite[Remark 3.18 p.~37]{Ten24}. Our proof replaces the density argument in a Hilbert space completion of $M \boxtimes N$ by an argument which does not use any Hilbert space completion of $M \boxtimes N$. We use a nondegenerate invariant Hermitian form on $M \boxtimes N$, the invariance of the relevant spaces under polynomials in $L_0$, and a moment argument (Lemma \ref{lem:moment-positivity}). In particular, the unitarity of $\ovl{M} \boxtimes M$ used in the proofs of the finite index theorems \cite[Theorem 4.20 p.~55 and Corollaries 4.21 and 4.24 pp.~55--56]{Ten24} holds in their setting (Remark \ref{rem:krein}), and Remark \ref{rem:delta} explains how this bears on the outline of \cite[Section 6 p.~8]{Hua26}.

Theorem \ref{thm:intro-E} proves the VOA--net correspondence for $S$ under $\mathrm{(H_J)}$, in the setting of \cite[Conjecture 7.6 p.~14]{Kaw18} and \cite{Gui26,Gui25,Ten24}. Denote by $\mathrm{(N_0)}$ the statement that $\Rep(\cal A_S)$ is braided equivalent to $Z(\cal H_3)$, and by $\mathrm{(N)}$ the assertion that there is a braided equivalence $\Rep(\cal A_S) \to \Rep(S)$ which maps $\pi_M$ to $M$ for $M \in \{M_b,M_a,M_{c_0}\}$. Under $\mathrm{(H_J)}$, $\mathrm{(N)}$ implies $\mathrm{(N_0)}$. We first prove $\mathrm{(N_0)}$ under $\mathrm{(H_J)}$, through the following recognition theorem, which does not use vertex algebras.

\begin{thm}
\label{thm:intro-D}
Let $\cal F$ be a unitary fusion category with six simple objects up to isomorphism, three of which are invertible, and with noncommutative Grothendieck ring. Then $\cal F$ has the fusion rules of $\cal H_2$ and $\cal H_3$. If moreover the fusion subcategory generated by the invertible objects is equivalent to $\mathrm{Vec}_{\Z_3}$ with trivial associator, then $\cal F$ is equivalent to $\cal H_2$ or to $\cal H_3$.
\end{thm}

This is Theorem \ref{thm:recognition}. We follow the strategy outlined by Huang \cite[Sections 5 and 6 pp.~7--9]{Hua26} and give complete proofs (Remark \ref{rem:huang-recognition}). The hypotheses force $\rho^2 \cong 1 \oplus m(\rho \oplus \alpha\rho \oplus \alpha^2\rho)$ for a non-invertible simple object $\rho$, an invertible object $\alpha$ of order three and some $m \geq 1$ (Lemma \ref{lem:ring}). The Frobenius--Schur indicators $\nu_1$, $\nu_2$, $\nu_3$ \cite{NS07,NS07b} give three identities for the twists of $Z(\cal F)$ (Lemma \ref{lem:moments}), and lower bounds for sums of roots of unity \cite{Lar14}, completed by an average over the Galois group of $\mathbb{Q}(\sqrt{10})$ when the squarefree part of $9m^2+4$ is $10$ (Lemma \ref{lem:galois-average}), give $m = 1$ (Proposition \ref{prop:m-equals-one}). The classification of generalized Haagerup categories by Izumi \cite{Izu18} then gives $\cal H_2$ or $\cal H_3$. The assumption on the associator is detected by the Drinfeld center. The third Frobenius--Schur indicator of $\alpha$ is equal to $1$ if and only if this associator is trivial (Lemma \ref{lem:nu3}), and $\cal F$ is equivalent to $\cal H_2$ or $\cal H_3$ if and only if no simple object of $Z(\cal F)$ has a twist of order nine (Proposition \ref{prop:twisted}). For unitary fusion categories this gives a recognition result for the Grothendieck ring studied in \cite[p.~5]{GrS12}. For the rank six categories of part 2 of Theorem \ref{thm:intro-B}, denoted $\cal F_x$ in part 4 of Theorem \ref{thm:complete-rationality}, one for each lattice extension, the triviality of the associator follows from the existence of two lattice extensions of $V'$ (Lemma \ref{lem:trivial-associator}), an input also used in \cite[p.~7]{Hua26}, and Theorem \ref{thm:intro-D} gives $\mathrm{(N_0)}$ (Theorem \ref{thm:Fx-haagerup}).

\begin{thm}
\label{thm:intro-C}
Assume $\mathrm{(H_J)}$.
\begin{enumerate}
\item $\Rep(\cal A_S)$ is equivalent to $Z(\cal H_3)$ as a unitary modular tensor category, in particular $\mathrm{(N_0)}$ holds. Every braided equivalence maps $\pi_b$ to $b$, and $\delta = d^2$.
\item $\cal A_S$ has exactly three irreducible local extensions with trivial representation category, up to isomorphism of extensions. Two of them are the nets $\cal A_{H_\Lambda}$ of the two lattice extensions, both isomorphic to $\cal A_{E_8}$, and the categories of modules over their dual canonical endomorphisms are $\cal H_2$ and $\cal H_3$. For the third one, this category is $\cal H_1$, and the only intermediate nets are $\cal A_S$ and the extension itself. For each lattice extension, the only intermediate nets are $\cal A_S$, $\cal A_{V'}$ and the extension itself. The net of the lattice extension of type $a$ carries $\cal H_2$, as the vertex algebra does.
\item The Haagerup subfactor arises from $\cal A_S$ in the sense of \cite[Definition 4.1 p.~20]{Bis17a}.
\item $\cal A_S$ is the generalized orbifold of $\cal A_{V'}$ by the hypergroup $K = \{k_0,k_1\}$ with $k_1k_1 = d^{-2}k_0+(1-d^{-2})k_1$. Hence, under these hypotheses, $\cal A_S$ solves \cite[Problem 4.27 p.~27]{Bis17b}.
\end{enumerate}
\end{thm}

This is Theorem \ref{thm:consequences} with Theorem \ref{thm:Fx-haagerup} and Corollaries \ref{cor:bischoff}, \ref{cor:orbifold} and \ref{cor:N-consequences}. The last assertion of part 2 uses Theorem \ref{thm:intro-E} below. Part 1 follows from Theorem \ref{thm:intro-D} and the dimensions computed in Theorem \ref{thm:intro-B} (Theorem \ref{thm:Fx-haagerup} and Lemma \ref{lem:N0}). Part 4 rests on the fact, which does not use Theorem \ref{thm:intro-D}, that $\cal A_S$ is a generalized orbifold of $\cal A_{V'}$ by a two-element hypergroup with parameter $\delta$ (Proposition \ref{prop:hypergroup}). By Theorem \ref{thm:intro-E}, the three extensions of part 2 are all isomorphic to $\cal A_{E_8}$, the third one being the net of the holomorphic extension of $S$ attached to $\cal H_1$ (Proposition \ref{prop:third-net} and Corollary \ref{cor:N-consequences}). For the algebra of Riesen in case (I), \cite[Theorem 7.13(1)--(2)]{ArhIII} gives the Haagerup property and complete unitarity. The assertions of Theorem \ref{thm:intro-C} for its specified inclusion in $V'$ use the compatible unitary structure asserted in \cite[Theorem 1 p.~1 and p.~4]{Rie26} (Section \ref{sec:case-I}). For the last assertion of part 2 and for the isomorphism of the third extension of its net with $\cal A_{E_8}$, the proof goes through the computer-assisted construction of Sections \ref{subsec:diagonal} and \ref{subsec:complete-unitarity} applied to the twin algebra (Theorem \ref{rem:all-haagerup}).

\begin{thm}
\label{thm:intro-E}
Assume $\mathrm{(H_J)}$.
\begin{enumerate}
\item The algebra $S \ot S$ has a holomorphic extension isomorphic to $V_{E_8 \oplus E_8}$ in which $S \ot 1$ and $1 \ot S$ are mutual commutants, and $V_{E_8 \oplus E_8}$ has a structure of unitary vertex operator algebra for which $S \ot 1$ and $1 \ot S$ are unitary subalgebras.
\item The algebra $S$ is completely unitary, every unitary $S$-module is strongly integrable, and the functor of \cite[Section 3.4 p.~45]{Gui26} is an equivalence of braided $C^*$-tensor categories from the category $\Rep^u(S)$ of unitary $S$-modules to the category $\Rep^f(\cal A_S)$ of representations of $\cal A_S$ with finite statistical dimension, which maps $M$ to $\pi_M$ for every $S$-submodule $M$ of a lattice extension. In particular, $\mathrm{(N)}$ holds.
\item The twin algebra of $S$ in case (I) and the common Heisenberg commutant of the two algebras are completely unitary.
\end{enumerate}
\end{thm}

This is Propositions \ref{prop:diagonal} and \ref{prop:diagonal-unitary}, Theorems \ref{thm:S-cu} and \ref{thm:twin-cu} and Remark \ref{rem:twin-cu}. Part 1 rests on exact computer calculations in the weight-one and weight-two spaces of $V_{E_8 \oplus E_8}$, which reduce, in the first implementation, to spaces of dimensions $496$ and $644$, carried out with two separate implementations (Remark \ref{rem:S-cu-computations}). We determine the embedding of the weight-one space of $S \ot S$ up to automorphisms, the conformal vector of the coset of the affine subalgebra of $S \ot 1$ up to four choices, and a positive rescaling of the standard antilinear involution of the lattice algebra which fixes these four vectors. Since $S \ot 1$ and $1 \ot S$ are commutants of explicit sets of vectors, they are then stable under this involution. Part 2 follows from the theorems of Gui on the strongly rational cosets of lattice vertex operator algebras \cite[Theorems 3.34 and 3.46 pp.~51 and 57]{Gui26}, \cite[Theorem 6.18 p.~70]{Gui25}. It gives a second proof of $\mathrm{(H_U)}$, independent of Theorem \ref{thm:intro-U} but computer-assisted. Part 3 follows from part 2 and the twin construction (Theorem \ref{thm:twin-cu}). In addition, Corollary \ref{cor:twin-comparison} transfers Condition I+ of \cite[Definition 6.17 p.~70]{Gui25} from $S$ to the twin algebra $S^\natural$ and to the common Heisenberg commutant $\Omega_S$ of $S$ and $S^\natural$, and it proves the analogue of $\mathrm{(N)}$ for every Haagerup algebra realized as a unitary conformal subalgebra of $V'$. The proofs of Propositions \ref{prop:diagonal} and \ref{prop:diagonal-unitary}, of Lemmas \ref{lem:sl2-four} and \ref{lem:omega-C} and of parts 1 and 2 of Theorem \ref{thm:S-cu} use no property of $S$ beyond those of an arbitrary Haagerup algebra in case (II). Consequently, every Haagerup algebra is completely unitary, and the analogue of $\mathrm{(N)}$ holds for the algebra $U^{\mathrm R}$ of \cite{Rie26} under $\mathrm{(H_R)}$ and the unitarity asserted in \cite{Rie26}, without $\mathrm{(H_J)}$ (Theorem \ref{rem:all-haagerup}). For the explicit twin of $S$, the identification with the algebra of \cite{Hua26} and its complete unitarity are given by \cite[Theorem 7.13(1)--(2)]{ArhIII}.

\paragraph{Application to the explicit algebra.} The inputs for the net of $S$ are its unitary realization and weight-one data, denoted \textup{(I1)} in \cite[Theorem 2.1]{ArhIII}, and its strong rationality and lattice branching, \textup{(F1)} and \textup{(F3)} of \cite[Theorem 2.2 and Appendix A]{ArhIII}. From these, \cite[Theorem 3.1(1)--(3)]{ArhIII} proves finite index, complete rationality and $\Rep(\cal A_S)\cong Z(\cal H_3)$. The proof uses the lattice net and submodule results of Section \ref{sec:net}, Proposition \ref{prop:fusion-unitary} for the transport modules of the lattice summands, and Theorem \ref{thm:recognition} for the unitary fusion categories of the net extensions.

For the vertex operator algebra, \cite[Theorems 4.25, 5.1 and 6.1]{ArhIII} constructs a compatible unitary diagonal extension from the net and uses the coset comparison. With \textup{(I1)} and \textup{(F1)}--\textup{(F3)} this proves $\mathrm{(H_J)}$, complete unitarity and Condition I+ for $S$. The hypothesis $\mathrm{(H_J)}$ of Theorems \ref{thm:intro-B}, \ref{thm:intro-U}, \ref{thm:intro-C} and \ref{thm:intro-E} holds by \cite[Theorem 6.1]{ArhIII}, whose proof uses no statement of this article which assumes $\mathrm{(H_J)}$ (Section \ref{sec:N}).

Table \ref{tab:dependencies} records the hypotheses of the proofs and the sources for the application to $S$.

\begin{table}[!htbp]
\centering
\small
\begin{tabular}{>{\raggedright\arraybackslash}p{7.9cm}>{\raggedright\arraybackslash}p{5.6cm}}
\hline
Statement & Status \\
\hline
Three Lagrangian algebras of $Z(\cal H_3)$, their objects in the labels of \cite{EvG11} and their categories of modules (Propositions \ref{prop:three-lagrangian} and \ref{prop:which-line}, Corollary \ref{cor:categorical}) & categorical proof \\
Holomorphic extensions of a Haagerup algebra, realization of $L_2$ and $L_3$ by the lattice extensions (Theorem \ref{thm:lagrangian-extensions}, Propositions \ref{prop:lagrangian-objects} and \ref{prop:which-line}) & for any Haagerup algebra. For $S$, \cite[Theorem 6.1]{ArhIII} supplies $\mathrm{(H_J)}$ \\
Strong locality of $S$, inclusions of nets (Theorem \ref{thm:strong-locality}) & proved, using the computer-assisted \cite[Theorem 10.10]{Arh26} \\
Recognition theorem for $\cal H_2$ and $\cal H_3$, twists of order nine (Theorem \ref{thm:recognition}, Proposition \ref{prop:twisted}) & categorical proof \\
Transport modules for two ordinary unitary modules (Proposition \ref{prop:fusion-unitary}) & regular simple unitary VOA with bounded localized vertex operators and an annular representation of the second module \\
Unitarizability of all simple modules (Corollary \ref{cor:modules-unitary}) & unitary conformal Haagerup subalgebra of $V'$. For $S$, also \cite[Theorem 6.1]{ArhIII} \\
Finite index, complete rationality, branching of sectors and $\delta>1$ (Theorem \ref{thm:complete-rationality} and Proposition \ref{prop:delta-big}) & unitary conformal Haagerup subalgebra of $V'$. For $S$, \cite[Theorem 3.1(1)--(2)]{ArhIII} uses \textup{(I1)}, \textup{(F1)}, \textup{(F3)} \\
$\cal F_x \cong \cal H_2$ or $\cal H_3$, $\mathrm{(N_0)}$, $\delta=d^2$, local extensions and the consequences of Bischoff (Theorem \ref{thm:Fx-haagerup} and Section \ref{sec:N}) & same hypotheses for Theorem \ref{thm:Fx-haagerup}. For $S$, \cite[Theorem 3.1(2)--(4)]{ArhIII} gives these net results \\
Agreement of the VOA and net labels, and the third holomorphic extension of $\cal A_S$ is $\cal A_{E_8}$ (Corollary \ref{cor:N-consequences}) & uses $\mathrm{(N)}$ and complete unitarity from Theorem \ref{thm:S-cu}. For $S$, see \cite[Corollary 6.2(1)--(2)]{ArhIII} \\
Diagonal extension $V_{E_8 \oplus E_8} \supset S \ot S$ with compatible unitary structure (Propositions \ref{prop:diagonal} and \ref{prop:diagonal-unitary}) & algebraic extension under $\mathrm{(H_J)}$. Section \ref{subsec:complete-unitarity} computes its PCT operator \\
Complete unitarity of $S$ and the comparison through $M\mapsto\pi_M$ (Theorem \ref{thm:S-cu}) & $\mathrm{(H_J)}$ and the computed unitary diagonal. The analytic proof for $S$ is in \cite[Theorems 5.1 and 6.1 and Corollary 6.2(1)]{ArhIII} \\
Equivalence of complete unitarity for a Haagerup algebra, its twin and their Heisenberg commutant (Theorem \ref{thm:twin-cu}(2)) & uses the gluing data of \cite[Section 12]{Arh26} \\
Condition I+ for $S^\natural$ and $\Omega_S$ (Corollary \ref{cor:twin-comparison}(1)) & transfer from $S$. For the explicit pair, see \cite[Theorem 7.13(2)]{ArhIII} \\
Complete unitarity and Condition I+ for every Haagerup algebra, and the analogue of $\mathrm{(N)}$ for unitary conformal embeddings (Theorem \ref{rem:all-haagerup}) & uses Section \ref{subsec:complete-unitarity} in case (II) and transfers through the twin in case (I) \\
Complete unitarity of the explicit case (I) algebra and comparison for its lattice realization (Section \ref{sec:case-I}) & complete unitarity by \cite[Theorem 7.13(2)]{ArhIII}. The comparison uses unitarity for the specified embedding \\
\hline
\end{tabular}
\caption{Dependence of the main statements on the hypotheses.}
\label{tab:dependencies}
\end{table}

\paragraph{Structure of the paper.} Section \ref{sec:preliminaries} collects the facts on Haagerup categories, Lagrangian algebras, vertex operator algebras and conformal nets that we use. Section \ref{sec:lagrangian} proves Theorem \ref{thm:intro-A}. Section \ref{sec:net} constructs the net of $S$. Section \ref{sec:rationality} proves the unitarity of fusion products, the unitarizability of the simple modules and the complete rationality. Section \ref{sec:recognition} proves Theorem \ref{thm:intro-D} and deduces $\mathrm{(N_0)}$. Section \ref{sec:N} derives the consequences of $\mathrm{(N_0)}$, constructs the unitary structure on $V_{E_8 \oplus E_8}$ and proves the complete unitarity of $S$ and $\mathrm{(N)}$. Section \ref{sec:case-I} discusses case (I) and the twin algebra.

%%%%%%%%%%%%%%%%%%%%%%%%%%%%%%%%%%%%%%%%%%%%%%%%%%%%%%%%%%%%%%%%%%%%%%%%%%%%%
\section{Preliminaries}
\label{sec:preliminaries}

%%%%%%%%%%%%%%%%%%%%%%%%%%%%%%%%%%%%%%%%%%%%%%%%%%%%%%%%%%%%%%%%%%%%%%%%%%%%%
\subsection{The Haagerup fusion categories}
\label{subsec:haagerup-categories}

We follow the conventions of Grossman and Snyder \cite{GrS12}. The two even parts of the Haagerup subfactor are unitary fusion categories $\cal H_1$ and $\cal H_2$ with four and six simple objects respectively \cite[p.~14]{GrS12}.

The simple objects of $\cal H_1$ are $1$, $\nu$, $\eta$, $\mu$, with Frobenius--Perron dimensions $1$, $d+1$, $d$, $d-1$. They are self-dual and the fusion rules are given in \cite[Table 2 p.~15]{GrS12}. In particular
\begin{equation}
\label{eq:H1-rules}
\mu \ot \mu \cong 1 \oplus \nu,\quad \eta \ot \eta \cong 1 \oplus \nu \oplus \eta \oplus \mu,\quad \nu \ot \nu \cong 1 \oplus 2\nu \oplus 2\eta \oplus \mu.
\end{equation}

The simple objects of $\cal H_2$ are $1$, $\alpha$, $\alpha^2$, $\xi$, $\alpha\xi$, $\alpha^2\xi$, with Frobenius--Perron dimensions $1$, $1$, $1$, $d$, $d$, $d$ and fusion rules
\begin{equation}
\label{eq:H2-rules}
\alpha^3 \cong 1,\quad \alpha\xi \cong \xi\alpha^2,\quad \xi^2 \cong 1 \oplus \xi \oplus \alpha\xi \oplus \alpha^2\xi.
\end{equation}
The category $\cal H_3$ is defined as the category of $(1+\alpha+\alpha^2)$-$(1+\alpha+\alpha^2)$-bimodules in $\cal H_2$ \cite[Definition 3.19 p.~23]{GrS12}. It has six simple objects, three of dimension $1$ and three of dimension $d$, and the same fusion ring as $\cal H_2$ \cite[Theorem 1.4 p.~4]{GrS12}. The three categories have the same global dimension
\begin{equation}
\label{eq:global-dimension}
\dim \cal H_i = 3+3d^2 = 3(3d+2), \quad i \in \{1,2,3\}.
\end{equation}

We will use the following facts.

\begin{prop}
\label{prop:GS}
\begin{enumerate}
\item The Morita equivalence class of $\cal H_1$ consists of exactly the three fusion categories $\cal H_1$, $\cal H_2$, $\cal H_3$, which are pairwise inequivalent.
\item For any $i \in \{1,2,3\}$, the category $\cal H_i$ has exactly three indecomposable module categories up to equivalence, and their dual categories are $\cal H_1$, $\cal H_2$ and $\cal H_3$.
\item The category $\cal H_1$ has rank four, and $\cal H_2$, $\cal H_3$ have rank six.
\item The only fusion subcategories of $\cal H_1$ are $\mathrm{Vec}$ and $\cal H_1$. For $i \in \{2,3\}$, the only fusion subcategories of $\cal H_i$ are $\mathrm{Vec}$, the pointed subcategory $\cal P_i$ generated by the invertible objects, and $\cal H_i$. Moreover $\cal P_i$ has dimension three.
\end{enumerate}
\end{prop}

\begin{proof}
1. This is \cite[Theorem 1.4 p.~4 and Corollary 3.22 p.~24]{GrS12}.

2. For $\cal H_2$, the three module categories are given by \cite[Corollary 3.16 p.~22]{GrS12}, and their duals are described in \cite[Lemmas 3.17 and 3.18 p.~22, Definition 3.19 p.~23 and Corollary 3.22 p.~24]{GrS12}. For $\cal H_1$ and $\cal H_3$ this is \cite[Theorem 3.25 p.~24]{GrS12}.

3. This was recalled above.

4. A fusion subcategory is determined by the set of its simple objects, which is closed under fusion and duality. By \eqref{eq:H1-rules}, any fusion subcategory of $\cal H_1$ containing $\mu$ contains $\nu$, any fusion subcategory containing $\nu$ or $\eta$ contains $\eta$, $\mu$ and $\nu$. Hence a nontrivial fusion subcategory of $\cal H_1$ is $\cal H_1$. By \eqref{eq:H2-rules}, each object $\alpha^k\xi$ is self-dual and $\alpha^k\xi \ot \alpha^k\xi \cong \alpha^{3k}\xi^2 \cong \xi^2$. Hence a fusion subcategory of $\cal H_2$ containing one of the objects $\alpha^k\xi$ contains the three objects of dimension $d$. It then contains $\xi \ot \alpha\xi$, which contains $\alpha^2$ by \cite[Table 1 p.~15]{GrS12}, so that it is $\cal H_2$. The same argument applies to $\cal H_3$, which has the same fusion ring. The remaining subcategories are $\mathrm{Vec}$ and $\cal P_i$, and $\dim \cal P_i = 3$.
\end{proof}

Since the categories $\cal H_i$ are Morita equivalent, their Drinfeld centers $Z(\cal H_1)$, $Z(\cal H_2)$ and $Z(\cal H_3)$ are braided equivalent \cite[Theorem 3.1 p.~9]{ENO11}, and we denote any of them by $Z(\cal H_3)$. By \cite[Theorem 1.5 p.~5]{GrS12}, there is exactly one Morita equivalence between any two of the categories $\cal H_i$. Hence the Brauer--Picard group of $\cal H_3$ is trivial, so that every braided autoequivalence of $Z(\cal H_3)$ is isomorphic to the identity functor \cite[Theorem 1.1 p.~3]{ENO10}. Consequently, any two braided equivalences between a braided fusion category and $Z(\cal H_3)$ are monoidally naturally isomorphic, and the labels of the simple objects of $Z(\cal H_3)$ used below are intrinsic.

%%%%%%%%%%%%%%%%%%%%%%%%%%%%%%%%%%%%%%%%%%%%%%%%%%%%%%%%%%%%%%%%%%%%%%%%%%%%%
\subsection{The Drinfeld center}
\label{subsec:center}

We use the modular data of $Z(\cal H_3)$ in the form of \cite[equations (2.5a) and (2.5b) p.~20]{EvG11}, as recalled in \cite[Section 2.1]{Arh26}. The simple objects are labelled
\[
0,\ b,\ a,\ c_0,\ c_1,\ c_2,\ d_1, \ldots, d_6,
\]
with $0$ the unit. Their quantum dimensions are
\begin{equation}
\label{eq:center-dimensions}
1,\quad 3d+1,\quad 3d+2,\quad 3d+2,\quad 3d+2,\quad 3d+2,\quad 3d, \ldots, 3d,
\end{equation}
and their twists are
\begin{equation}
\label{eq:center-twists}
1,\ 1,\ 1,\ 1,\ \e^{\frac{2\pi\i}{3}},\ \e^{-\frac{2\pi\i}{3}},\ \e^{\frac{2\pi\i 6l^2}{13}}\quad (1 \leq l \leq 6).
\end{equation}
The charge conjugation $S^2$ is trivial, since the $S$-matrix of \cite[(2.5a)--(2.5b) p.~20]{EvG11} is real, symmetric and unitary, so that every simple object is self-dual, see also \cite[Section 2.1]{Arh26}. Since $Z(\cal H_3)$ is unitary, the quantum dimensions coincide with the Frobenius--Perron dimensions, and
\begin{equation}
\label{eq:dim-center}
\operatorname{FPdim} Z(\cal H_3) = (\dim \cal H_3)^2 = 9(3d+2)^2.
\end{equation}
The twists $\e^{\frac{2\pi\i 6l^2}{13}}$ are different from $1$, since $13$ does not divide $6l^2$ for $1 \leq l \leq 6$. Hence the simple objects with trivial twist are exactly $0$, $b$, $a$ and $c_0$.

The modular data do not distinguish $a$ from $c_0$. In \cite{EvG11}, the labels are defined through the half-braidings of \cite[Section 8]{Izu01} on Izumi's model $\{1,\alpha,\alpha^2,\xi,\alpha\xi,\alpha^2\xi\}$ of the six-object even part of the Haagerup subfactor \cite[p.~17]{EvG11}, \cite[pp.~56--57]{Izu18}. This even part is $\cal H_2$, which is defined in \cite[p.~14]{GrS12} as the even part with six simple objects of the Haagerup subfactor. With $\mu = \xi \oplus \alpha\xi \oplus \alpha^2\xi$, the forgetful functor $F \co Z(\cal H_2) \to \cal H_2$ maps $b$ to $1 \oplus \mu$, $a$ to $2 \cdot 1 \oplus \mu$, the objects $c_j$ to $\alpha \oplus \alpha^2 \oplus \mu$ and the objects $d_l$ to $\mu$ \cite[p.~18]{EvG11}. By Section \ref{subsec:haagerup-categories}, these labels are intrinsic. We use them throughout.

%%%%%%%%%%%%%%%%%%%%%%%%%%%%%%%%%%%%%%%%%%%%%%%%%%%%%%%%%%%%%%%%%%%%%%%%%%%%%
\subsection{\'Etale and Lagrangian algebras}
\label{subsec:etale}

Let $\cal C$ be a non-degenerate braided fusion category. An \'etale algebra in $\cal C$ is a commutative separable algebra, and it is connected if $\dim \operatorname{Hom}_{\cal C}(1,A) = 1$. For a connected \'etale algebra $A$, the category $\cal C_A$ of right $A$-modules is a fusion category and the full subcategory $\cal C_A^0$ of local (dyslectic) modules is a non-degenerate braided fusion category \cite[Corollary 3.30 p.~15]{DMNO13}, with
\begin{equation}
\label{eq:dims-local}
\operatorname{FPdim}(\cal C_A) = \frac{\operatorname{FPdim}(\cal C)}{\operatorname{FPdim}(A)},\quad \operatorname{FPdim}(\cal C_A^0) = \frac{\operatorname{FPdim}(\cal C)}{\operatorname{FPdim}(A)^2},
\end{equation}
see \cite[Lemma 3.11 p.~10 and Corollary 3.32 p.~15]{DMNO13}. A connected \'etale algebra $L$ is called Lagrangian if $\operatorname{FPdim}(L)^2 = \operatorname{FPdim}(\cal C)$ \cite[Definition 4.6 p.~17]{DMNO13}, which is equivalent to $\cal C_L^0 \cong \mathrm{Vec}$ \cite[Remark 4.3 p.~16]{DMNO13}.

\begin{prop}
\label{prop:DMNO}
Let $\cal C$ be a non-degenerate braided fusion category.
\begin{enumerate}
\item Let $L$ be a Lagrangian algebra in $\cal C$. The free module functor $F_L \co \cal C \to \cal C_L$, $X \mapsto X \ot L$, lifts to a braided equivalence $\cal C \to Z(\cal C_L)$, under which $F_L$ corresponds to the forgetful functor $Z(\cal C_L) \to \cal C_L$.
\item Assume that $\cal C = Z(\cal D)$ for a fusion category $\cal D$. Then the isomorphism classes of Lagrangian algebras of $\cal C$ are in bijection with the equivalence classes of indecomposable $\cal D$-module categories, and for the Lagrangian algebra $L$ corresponding to $\cal M$ the fusion category $\cal C_L$ is equivalent to the dual category $\cal D^*_{\cal M}$.
\item Let $\cal D$ be a fusion category and $I \co \cal D \to Z(\cal D)$ the right adjoint of the forgetful functor. The \'etale subalgebras of the Lagrangian algebra $I(1)$ are in order-reversing bijection with the fusion subcategories of $\cal D$, and if $B \subset I(1)$ corresponds to $\cal B \subset \cal D$ then $\operatorname{FPdim}(B)\operatorname{FPdim}(\cal B) = \operatorname{FPdim}(\cal D)$.
\item Let $A$ be a connected \'etale algebra in $\cal C$. The commutative algebras $B$ in $\cal C$ equipped with a unital algebra morphism $A \to B$, which is automatically injective \cite[p.~12]{DMNO13}, correspond to the commutative algebras in $\cal C_A^0$. Under this correspondence, connected \'etale algebras correspond to connected \'etale algebras.
\item Let $\cal D$ be a fusion category, let $\cal B \subset \cal D$ be a fusion subcategory and let $A_{\cal B} \subset I(1)$ be the \'etale subalgebra attached to $\cal B$ in part 3. Then the category $Z(\cal D)^0_{A_{\cal B}}$ of local $A_{\cal B}$-modules is braided equivalent to $Z(\cal B)$.
\end{enumerate}
\end{prop}

\begin{proof}
1. This is \cite[Corollary 4.1 p.~16]{DMNO13}, with the central structure \cite[(17) p.~10]{DMNO13}.

2. This is \cite[Proposition 4.8 and Remark 4.9 p.~17]{DMNO13}, with the $2$-equivalence of \cite[Proposition 4.4 p.~16]{DMNO13}.

3. This is \cite[Theorem 4.10 p.~18]{DMNO13}.

4. This is \cite[Section 3.6 and Proposition 3.16 p.~12]{DMNO13}.

5. Put $\cal C = Z(\cal D)$ and $A = A_{\cal B}$. In the proof of \cite[Theorem 4.10 p.~18]{DMNO13}, the category $\cal C_A$ is identified with the relative center $Z_{\cal B}(\cal D)$, whose objects are the objects of $\cal D$ equipped with a half-braiding with respect to the objects of $\cal B$, in such a way that the free module functor $\cal C \to \cal C_A$, with its central structure \cite[(17) p.~10]{DMNO13}, corresponds to the forgetful functor $Z(\cal D) \to Z_{\cal B}(\cal D)$ with the central structure given by the half-braidings. Since $\cal C$ is non-degenerate, there is a braided equivalence $Z(\cal C_A) \cong \cal C \boxtimes (\cal C_A^0)^{\mathrm{rev}}$ under which $\cal C \boxtimes 1$ corresponds to the image of the central structure of the free module functor \cite[Corollary 3.30 and Remark 3.31 p.~15]{DMNO13}. On the other hand, $Z_{\cal B}(\cal D)$ is the dual of $\cal D \boxtimes \cal B^{\mathrm{op}}$ with respect to the module category $\cal D$ \cite[p.~18]{DMNO13}, so that Schauenburg's equivalence \cite[(12) p.~7]{DMNO13} gives a braided equivalence $Z(Z_{\cal B}(\cal D)) \cong Z(\cal D) \boxtimes Z(\cal B)^{\mathrm{rev}}$, under which the central structure of the forgetful functor $Z(\cal D) \to Z_{\cal B}(\cal D)$ corresponds to $Z(\cal D) \boxtimes 1$ \cite[pp.~18--19]{DMNO13}. The centralizer of $\cal C \boxtimes 1$ in $\cal C \boxtimes \cal E$ is $1 \boxtimes \cal E$ for any braided fusion category $\cal E$, since $\cal C$ is non-degenerate. Comparing the centralizers of $\cal C \boxtimes 1$ in the two descriptions of $Z(\cal C_A)$, we obtain a braided equivalence $(\cal C_A^0)^{\mathrm{rev}} \cong Z(\cal B)^{\mathrm{rev}}$, hence $\cal C_A^0 \cong Z(\cal B)$. This is the claim proved in \cite[proof of Theorem 4.10 p.~18]{DMNO13}.
\end{proof}

%%%%%%%%%%%%%%%%%%%%%%%%%%%%%%%%%%%%%%%%%%%%%%%%%%%%%%%%%%%%%%%%%%%%%%%%%%%%%
\subsection{Haagerup algebras and their extensions}
\label{subsec:haagerup-algebras}

A vertex operator algebra is strongly rational if it is simple, of CFT type, self-contragredient, rational and $C_2$-cofinite. It is positive if every irreducible module not isomorphic to the vertex operator algebra itself has positive lowest conformal weight. Following \cite{Arh26}, a Haagerup algebra is a positive strongly rational vertex operator algebra $U$ of central charge eight together with a braided equivalence $\Rep(U) \cong Z(\cal H_3)$. The modular tensor category $\Rep(U)$ is given by \cite{Hua08}. Since $Z(\cal H_3)$ has no nontrivial invertible object, such an equivalence preserves the ribbon structures \cite[Section 2.2 and Remark 12.6]{Arh26}. By Section \ref{subsec:haagerup-categories}, it is unique up to natural isomorphism, so that the simple $U$-modules $M_z$ below are well defined.

Let $U$ be a strongly rational vertex operator algebra. An extension of $U$ is a vertex operator algebra $B$ containing $U$ as a vertex subalgebra with the same conformal vector and with $\dim B_0 = 1$. Two extensions are isomorphic if there is an isomorphism of vertex operator algebras between them which restricts to the identity of $U$. An extension $B$ is holomorphic if it is simple and if every $B$-module is a direct sum of copies of $B$, that is, $\Rep(B) \cong \mathrm{Vec}$. By \cite[Theorem 3.2 p.~10]{HKL15}, extensions of $U$ correspond to haploid commutative associative algebras with trivial twist in $\Rep(U)$, and by \cite[Theorem 3.4 p.~16]{HKL15} the category of modules of an extension $B$ is isomorphic to the category of local modules of the corresponding algebra. If the algebra is moreover rigid in the sense of \cite[Definition 1.11 p.~10]{KiO02}, then its category of modules is semisimple \cite[Theorem 3.3 p.~17]{KiO02}, so that it is separable \cite[Proposition 2.7 p.~5]{DMNO13}, hence \'etale.

\begin{lemma}
\label{lem:intermediate-etale}
Let $U$ be a positive strongly rational vertex operator algebra and $B$ a simple extension of $U$ with $B \cong B'$. Let $W$ be a vertex subalgebra of $B$ containing $U$. Then $W$ is simple if and only if the invariant bilinear form of $B$ restricts to a nondegenerate form on $W$. In this case, the algebra in $\Rep(U)$ corresponding to $W$ is a connected \'etale subalgebra of the algebra corresponding to $B$.
\end{lemma}

\begin{proof}
Since $B$ is simple and isomorphic to $B'$, it carries a nondegenerate invariant symmetric bilinear form $(\cdot,\cdot)$, unique up to a scalar \cite{FHL93}, \cite{Li94}. Let $R = \{w \in W : (w,W) = 0\}$ be the radical of its restriction to $W$. The vertex subalgebra $W$ contains the conformal vector, hence it is stable under $L_1$. The invariance of the form expresses the adjoint of a mode $u_{(n)}$, for $u \in W$, as a finite combination of modes of the vectors $L_1^ju \in W$. Hence $R$ is an ideal of $W$. Since $(\mathbf 1,\mathbf 1) \neq 0$, we have $R \neq W$. Consequently, if $W$ is simple then $R = 0$.

Conversely, assume that $R = 0$, and let $J \subset W$ be a proper ideal. The space $J$ does not contain $\mathbf 1$, and since $\dim W_0 = 1$ and $J$ is graded, we have $J \cap W_0 = 0$. For $j \in J$ and $w \in W$, the invariance of the form gives $(j,w) = (j_{(-1)}\mathbf 1,w) = (\mathbf 1,j^\dagger w)$, where $j^\dagger$ is a finite linear combination of modes of the vectors $L_1^kj$. These vectors lie in $J$, since $J$ is an ideal containing $j$ and the conformal vector belongs to $W$. Hence $j^\dagger w \in J$. Since $J \cap W_0 = 0$, $J$ is graded and $(\mathbf 1,W_n) = 0$ for $n \neq 0$, we obtain $(j,w) = 0$. Thus $J \subset R = 0$, and $W$ is simple.

Finally, assume that $W$ is simple. By \cite[Theorem 3.2 p.~10]{HKL15}, $W$ corresponds to a haploid commutative algebra $A_W$ with trivial twist, which is a subalgebra of the algebra $A_B$ attached to $B$. The multiplication of $A_W$ corresponds to the vertex operator of $W$ \cite[proof of Theorem 3.2 p.~10]{HKL15}, so that the subobjects of $A_W$ which are $A_W$-submodules correspond to the ideals of $W$. Since $W$ is simple, $A_W$ is simple as an $A_W$-module. The categorical dimension of $A_W$ is a sum of quantum dimensions of simple $U$-modules, which are positive for positive strongly rational $U$ \cite[Lemma 4.2 p.~13]{DJX13}, \cite{Hua08}, so that it is nonzero. Since $\theta_{A_W} = \mathrm{id}$, \cite[Lemma 1.20 p.~15]{KiO02} shows that $A_W$ is rigid in the sense of \cite[Definition 1.11 p.~10]{KiO02}, its category of modules is semisimple \cite[Theorem 3.3 p.~17]{KiO02}, and $A_W$ is separable \cite[Proposition 2.7 p.~5]{DMNO13}, hence \'etale.
\end{proof}

We recall the facts from \cite{Arh26} that we use. Let $U$ be a Haagerup algebra.
\begin{enumerate}
\item By \cite[Proposition 4.1]{Arh26}, the subspace $V' = U \oplus M_b$ of the holomorphic extension of $U$ attached to an algebra isomorphic to $0 \oplus b \oplus 2a$ is a vertex operator algebra isomorphic to $V_{E_6 \oplus A_2}$, and $U$ is a conformal subalgebra of $V'$. We fix such an isomorphism. The discriminant group of $E_6 \oplus A_2$ is $\Z_3 \times \Z_3$, with quadratic form $q(i,j) = \frac{2i^2+j^2}{3}$ modulo one. We denote by $V'_{w}$, $w \in \Z_3 \times \Z_3$, the simple $V'$-modules.
\item Let $A_0$ be the algebra in $\Rep(U)$ attached to $V'$. By \cite[Lemma 4.3]{Arh26} and its proof, for $w \neq 0$ the restriction of $V'_w$ to $U$ is simple. This restriction is isomorphic to $M_{c_1}$ if $w \in \{(0,1),(0,2)\}$, to $M_{c_2}$ if $w \in \{(1,0),(2,0)\}$, and the four remaining classes, which are the nonzero elements of the two isotropic lines, restrict to $M_a$ for the two classes of one line and to $M_{c_0}$ for the two classes of the other line. Indeed, the proof of \cite[Lemma 4.3]{Arh26} shows that exactly two of these four classes restrict to $M_{c_0}$, and the set of such classes is stable under $w \mapsto -w$ since $V'_{-w}$ is the contragredient of $V'_w$ and $c_0$ is self-dual. In the coordinates chosen in \cite[Section 4.2]{Arh26}, for which the holomorphic extension of \cite[Proposition 4.1]{Arh26}, attached to an algebra isomorphic to $0 \oplus b \oplus 2a$, corresponds to the line generated by $(1,1)$, the classes restricting to $M_a$ are $(1,1)$ and $(2,2)$.
\end{enumerate}
We denote by $\Lambda_a$ and $\Lambda_{c_0}$ the isotropic lines whose nonzero classes restrict to $M_a$ and $M_{c_0}$ respectively, and we put
\begin{equation}
\label{eq:Hx}
H_a \ov{\mathrm{def}}{=} \bigoplus_{w \in \Lambda_a} V'_w,\quad H_{c_0} \ov{\mathrm{def}}{=} \bigoplus_{w \in \Lambda_{c_0}} V'_w.
\end{equation}
Each isotropic line $\Lambda$ defines an even unimodular overlattice of $E_6 \oplus A_2$, which is isometric to $E_8$, and $\bigoplus_{w \in \Lambda} V'_w$ is the corresponding lattice vertex operator algebra, see Lemma \ref{lem:cocycle} and the paragraph following it. Hence $H_a \cong H_{c_0} \cong V_{E_8}$, and as $U$-modules
\begin{equation}
\label{eq:Hx-decomposition}
H_a = U \oplus M_b \oplus M_a \oplus M_a,\quad H_{c_0} = U \oplus M_b \oplus M_{c_0} \oplus M_{c_0}.
\end{equation}

For the explicit algebras we use the following external hypotheses, introduced in \cite[Section 1]{Arh26}. We write $\mathrm{(H_J)}$ for the assertion that the explicit conformal algebra of \cite{Jia26} is a Haagerup algebra, and $\mathrm{(H_R)}$ for the same assertion for the algebra $U^{\mathrm R}$ of \cite{Rie26}. Positivity is part of these hypotheses. By \cite[Theorem 10.24]{Arh26}, Jiao's algebra is isomorphic to the algebra $S$, so that $\mathrm{(H_J)}$ implies that $S$ is a Haagerup algebra. Both hypotheses are proved in \cite[Theorems 6.1 and 7.13(1)]{ArhIII}, using results of this article, and we keep them as hypotheses below.

%%%%%%%%%%%%%%%%%%%%%%%%%%%%%%%%%%%%%%%%%%%%%%%%%%%%%%%%%%%%%%%%%%%%%%%%%%%%%
\subsection{Unitary vertex operator algebras and strong locality}
\label{subsec:unitary}

We use the framework of \cite{CKLW18}. A unitary vertex operator algebra is a vertex operator algebra $V$ with a normalized invariant scalar product, whose PCT operator $\theta$ is an antilinear involutive automorphism of $V$ \cite[Section 5]{CKLW18}. A vertex subalgebra $W$ of a unitary vertex operator algebra $V$ is a unitary subalgebra if and only if $\theta W \subset W$ and $L_1W \subset W$ \cite[Proposition 5.23 p.~43]{CKLW18}. A unitary vertex operator algebra with $V_0 = \Cbb\mathbf 1$ is simple \cite[Proposition 5.3 p.~34]{CKLW18}.

A simple unitary vertex operator algebra is strongly local if its vertex operators satisfy energy bounds and the von Neumann algebras generated by the smeared vertex operators localized in disjoint intervals commute. A strongly local vertex operator algebra $V$ defines an irreducible conformal net $\cal A_V$ on the Hilbert space completion $\cal H_V$ of $V$, by
\[
\cal A_V(I) = W^*\big(\{Y(a,f) : a \in V,\ f \in C^\infty(S^1),\ \operatorname{supp} f \subset I\}\big),
\]
the von Neumann algebra generated by the closures $Y(a,f)$ of the smeared vertex operators $Y_0(a,f)$ with domain $V$ \cite[(106) p.~49 and (114) p.~50]{CKLW18}, see \cite[Theorem 6.8 p.~53]{CKLW18}. We will use the following results.

\begin{prop}
\label{prop:CKLW}
Let $V$ be a strongly local simple unitary vertex operator algebra.
\begin{enumerate}
\item If $\{\cdot|\cdot\}$ is another normalized invariant scalar product on $V$, then $(V,\{\cdot|\cdot\})$ is strongly local, and the two conformal nets are isomorphic \cite[Theorem 6.8 p.~53]{CKLW18}.
\item A unitary subalgebra $W$ of $V$ is strongly local, and $\cal A_W$ embeds as a covariant subnet of $\cal A_V$ \cite[Theorem 7.1 p.~55]{CKLW18}. The map $W \mapsto \cal A_W$ is a bijection between unitary subalgebras of $V$ and Möbius covariant subnets of $\cal A_V$ \cite[Theorem 7.5 p.~58]{CKLW18}.
\item If $G$ is a closed subgroup of the group of unitary automorphisms of $V$, then $\cal A_{V^G} = \cal A_V^G$ \cite[Proposition 7.6 p.~58]{CKLW18}.
\item For a unitary subalgebra $W$ of $V$, the coset subnet $\cal A_W^c$, given by $\cal A_W^c(I) = \cal A_W(S^1)' \cap \cal A_V(I)$, is the net of the coset $W^c$ \cite[Proposition 7.8 p.~58]{CKLW18}.
\item A tensor product $V^\alpha \ot V^\beta$ of strongly local vertex operator algebras is strongly local, and $\cal A_{V^\alpha \ot V^\beta} = \cal A_{V^\alpha} \ot \cal A_{V^\beta}$ \cite[Corollary 8.2 p.~62]{CKLW18}. For a complex simple Lie algebra $\mathfrak g$ and a positive integer $k$, the simple affine vertex operator algebra $V_{\mathfrak g,k}$ is strongly local and its net is the loop group net $\cal A_{G_k}$ \cite[Example 8.7 p.~63]{CKLW18}.
\end{enumerate}
\end{prop}

%%%%%%%%%%%%%%%%%%%%%%%%%%%%%%%%%%%%%%%%%%%%%%%%%%%%%%%%%%%%%%%%%%%%%%%%%%%%%
\subsection{Conformal nets}
\label{subsec:nets}

Throughout, $\Rep(\cal A)$ denotes the category of finite direct sums of irreducible representations of a completely rational net. Equivalently, it is the category of representations of finite statistical dimension, also denoted by $\Rep^f(\cal A)$ below. This convention excludes infinite direct sums when we refer to a fusion category.

We refer to \cite{KLM01,BKLR15} for conformal nets, their representations and $\Q$-systems. A conformal net $\cal A$ is completely rational if it is split, strongly additive and has finite $\mu$-index. In this case every representation on a separable Hilbert space is a direct sum of irreducible representations with finite statistical dimension \cite[Corollary 39 p.~27]{KLM01}, there are finitely many irreducible sectors, $\Rep(\cal A)$ is a unitary modular tensor category, and
\begin{equation}
\label{eq:mu-index}
\mu_{\cal A} = \sum_{\pi \in \Irr(\Rep(\cal A))} d(\pi)^2
\end{equation}
by \cite[Theorem 33 p.~24]{KLM01}. By \cite{LoX04,MTW18}, complete rationality is equivalent to rationality in the sense of \cite[Section 4.4 p.~54]{Ten24}, namely finitely many irreducible representations, all with finite index.

For an inclusion $\cal A \subset \cal B$ of conformal nets, the index $[\cal B : \cal A]$ is the Jones index of $\cal A(I) \subset \cal B(I)$, which does not depend on $I$. We will use the following facts.

\begin{prop}
\label{prop:nets}
Let $\cal A \subset \cal B$ be a conformal subnet.
\begin{enumerate}
\item If $[\cal B : \cal A] < \infty$, then $\cal A$ is completely rational if and only if $\cal B$ is \cite[Theorem 24 p.~20]{Lon03}. In this case $\mu_{\cal A} = [\cal B : \cal A]^2\mu_{\cal B}$ \cite[Proposition 24 p.~19]{KLM01}.
\item If $\cal A$ is completely rational and $\cal B$ is an irreducible local extension of $\cal A$, in the sense that $\cal A(S^1)' \cap \cal B(I) = \Cbb$ \cite[p.~9]{KL04}, then $[\cal B : \cal A] < \infty$ \cite[Proposition 2.3 p.~12]{KL04}. The hypotheses of \cite[Proposition 2.3]{KL04} on the sectors of $\cal A$ hold by \cite[Theorem 33 p.~24 and Corollary 39 p.~27]{KLM01}.
\item If $\cal A$ is completely rational and $[\cal B : \cal A] < \infty$, the restriction of the vacuum representation of $\cal B$ to $\cal A$ is the dual canonical endomorphism $\theta$ of the inclusion, which carries the structure of a commutative $\Q$-system in $\Rep(\cal A)$, and $[\cal B : \cal A] = d(\theta)$ \cite[Theorem 4.9 p.~25]{LoR95}. Conversely, a commutative $\Q$-system $\Theta$ in $\Rep(\cal A)$ with $\dim\operatorname{Hom}(\mathrm{id},\Theta) = 1$ defines an irreducible local extension of $\cal A$, unique up to isomorphism of extensions \cite[Theorem 4.9 p.~25]{LoR95}, \cite[Section 5.2.1 p.~77]{BKLR15}. Intermediate extensions correspond to intermediate $\Q$-systems \cite[Proposition 4.11 p.~49]{BKLR15}.
\end{enumerate}
\end{prop}

A commutative $\Q$-system $\Theta$ with $\dim \operatorname{Hom}(\mathrm{id},\Theta) = 1$ is a connected \'etale algebra in the braided fusion category $\Rep(\cal A)$. Conversely, by \cite[Theorem 3.2 p.~11]{CGGH23}, a haploid algebra in a $C^*$-tensor category is equivalent to a $\Q$-system if and only if it is rigid. Since a connected \'etale algebra is rigid \cite[Remark 3.4 p.~7]{DMNO13}, a connected \'etale algebra in $\Rep(\cal A)$ is equivalent to a commutative $\Q$-system. A braided equivalence between unitary braided fusion categories is monoidally naturally isomorphic to a unitary braided equivalence \cite[Theorem 5 p.~2]{Reu23}.

We recall the results of Bischoff that we use. A conformal net $\cal A$ is called holomorphic if it is completely rational with $\mu_{\cal A} = 1$, that is, if $\Rep(\cal A) \cong \mathrm{Vec}$.

\begin{prop}
\label{prop:bischoff}
\begin{enumerate}
\item Let $\cal A$ be a completely rational net. Then $\Rep(\cal A)$ is braided equivalent to $Z(\cal F)$ for some unitary fusion category $\cal F$ if and only if $\cal A$ has a holomorphic finite index extension. In this case, $\cal F$ gives a holomorphic extension $\cal B_{\cal F}$ of $\cal A$, and there is a Galois correspondence between full fusion subcategories of $\cal F$ and intermediate nets $\cal A \subset \cal B \subset \cal B_{\cal F}$ \cite[Proposition 4.5 p.~21]{Bis17a}.
\item Let $N \subset M$ be a finite index, finite depth subfactor. If there is a completely rational net $\cal A$ such that $\Rep(\cal A)$ is braided equivalent to $D(N \subset M)$, then $N \subset M$ arises from $\cal A$ in the sense of \cite[Definition 4.1 p.~20]{Bis17a}, see \cite[Proposition 4.7 p.~22]{Bis17a}.
\item Let $\cal A \subset \cal B$ be a finite index subnet. There is a canonical hypergroup $K$ acting properly on $\cal B$ by vacuum preserving unital completely positive maps such that $\cal A = \cal B^K$, and $[\cal B : \cal A]$ is the weight of $K$ \cite[Theorem 1.3 p.~3 and Theorem 4.22 p.~24]{Bis17b}.
\end{enumerate}
\end{prop}

%%%%%%%%%%%%%%%%%%%%%%%%%%%%%%%%%%%%%%%%%%%%%%%%%%%%%%%%%%%%%%%%%%%%%%%%%%%%%
\section{Lagrangian algebras and holomorphic extensions}
\label{sec:lagrangian}

%%%%%%%%%%%%%%%%%%%%%%%%%%%%%%%%%%%%%%%%%%%%%%%%%%%%%%%%%%%%%%%%%%%%%%%%%%%%%
\subsection{A rank formula}

The following standard fact is the main tool of this section. In the language of Longo--Rehren inclusions, it is the Frobenius reciprocity argument used in \cite[p.~19]{EvG11} for the four-object even part of the Haagerup subfactor, and it is also used in \cite[(23) p.~7]{Hua26}. We include a proof.

\begin{lemma}
\label{lem:rank}
Let $\cal C$ be a non-degenerate braided fusion category and let $L$ be a Lagrangian algebra in $\cal C$. Then $\dim \operatorname{Hom}_{\cal C}(L,L)$ is equal to the number of isomorphism classes of simple objects of the fusion category $\cal D = \cal C_L$. Moreover, the braided equivalence $\Phi \co \cal C \to Z(\cal D)$ of part 1 of Proposition \ref{prop:DMNO} maps the algebra $L$ to an algebra isomorphic to $I(1)$, where $I$ is the right adjoint of the forgetful functor $Z(\cal D) \to \cal D$.
\end{lemma}

\begin{proof}
Put $\cal D = \cal C_L$. By part 1 of Proposition \ref{prop:DMNO}, there is a braided equivalence $\Phi \co \cal C \to Z(\cal D)$ such that $F \circ \Phi$ is isomorphic to the free module functor $F_L \co X \mapsto X \ot L$, where $F \co Z(\cal D) \to \cal D$ is the forgetful functor. The right adjoint of $F_L$ is the functor $\cal D \to \cal C$ which forgets the module structure, and it maps the unit object $L$ of $\cal D$ to $L$. Since $\Phi$ is an equivalence, the right adjoint of $F \circ \Phi$ is $\Phi^{-1} \circ I$, where $I$ is the right adjoint of $F$. By uniqueness of adjoints, we obtain an isomorphism $\Phi(L) \cong I(1)$. It is an isomorphism of algebras, since the algebra structures of $L$ and $I(1)$ are induced by the lax monoidal structures of the right adjoints of the tensor functors $F_L$ and $F$, and $F \circ \Phi \cong F_L$ as tensor functors \cite[Lemma 3.9 p.~10]{DMNO13}. Consequently, we have
\[
\operatorname{Hom}_{\cal C}(L,L)
\cong \operatorname{Hom}_{Z(\cal D)}(I(1),I(1))
\cong \operatorname{Hom}_{\cal D}(F(I(1)),1).
\]
By \cite[Proposition 5.4 p.~21]{ENO05}, \cite[Proposition 9.2.2 p.~278]{EGNO15}, we have $F(I(1)) \cong \bigoplus_{X} X \ot X^*$, where $X$ runs through a set of representatives of the isomorphism classes of simple objects of $\cal D$. Since $\dim \operatorname{Hom}_{\cal D}(X \ot X^*,1) = \dim \operatorname{Hom}_{\cal D}(X,X) = 1$ for any simple object $X$, the conclusion follows. This is \cite[Proposition 5.6 p.~21]{ENO05}, and the identification of $L$ with $I(1)$ is also a consequence of \cite[Proposition 4.4 p.~16]{DMNO13}.
\end{proof}

%%%%%%%%%%%%%%%%%%%%%%%%%%%%%%%%%%%%%%%%%%%%%%%%%%%%%%%%%%%%%%%%%%%%%%%%%%%%%
\subsection{The three Lagrangian algebras of the Haagerup double}

\begin{lemma}
\label{lem:possible-objects}
Let $L$ be a Lagrangian algebra in $Z(\cal H_3)$. Then $L$ is isomorphic, as an object, to $0 \oplus b \oplus X$ with $X \in \{a \oplus a,\ c_0 \oplus c_0,\ a \oplus c_0\}$.
\end{lemma}

\begin{proof}
The pairing $L \ot L \to L \to 1$ is nondegenerate \cite[Remark 3.4 p.~7]{DMNO13}, so that $L$ is a rigid haploid algebra in the unitary braided fusion category $Z(\cal H_3)$. By \cite[Theorem 3.2 p.~11]{CGGH23}, $L$ is equivalent to a standard commutative $\Q$-system, which has trivial twist by \cite[Proposition 2.22 p.~9]{CGGH23}. Hence $L$ is a direct sum of simple objects with trivial twist. By Section \ref{subsec:center}, we can write $L \cong 0 \oplus n_b b \oplus n_a a \oplus n_{c_0}c_0$ with nonnegative integers $n_b$, $n_a$, $n_{c_0}$, where the multiplicity of $0$ is one since $L$ is connected. By \eqref{eq:dim-center}, the Lagrangian condition reads
\[
1 + n_b(3d+1) + (n_a+n_{c_0})(3d+2) = 3(3d+2) = 9d+6.
\]
Since $d$ is irrational, we can compare the coefficients of $1$ and $d$ and obtain $1+n_b+2(n_a+n_{c_0}) = 6$ and $n_b+n_a+n_{c_0} = 3$. Hence $n_a+n_{c_0} = 2$ and $n_b = 1$.
\end{proof}

\begin{prop}
\label{prop:three-lagrangian}
The braided fusion category $Z(\cal H_3)$ has exactly three Lagrangian algebras up to isomorphism, $L_1$, $L_2$ and $L_3$, characterized by $Z(\cal H_3)_{L_i} \cong \cal H_i$. As objects, $L_1 \cong 0 \oplus b \oplus a \oplus c_0$, and $L_2$, $L_3$ belong to $\{0 \oplus b \oplus 2a,\ 0 \oplus b \oplus 2c_0\}$.
\end{prop}

\begin{proof}
By part 2 of Proposition \ref{prop:DMNO} applied to $\cal D = \cal H_3$, the Lagrangian algebras of $Z(\cal H_3)$ are in bijection with the indecomposable $\cal H_3$-module categories, and the category of modules over the Lagrangian algebra attached to $\cal M$ is the dual category $(\cal H_3)^*_{\cal M}$. By part 2 of Proposition \ref{prop:GS}, there are three such module categories and their duals are $\cal H_1$, $\cal H_2$ and $\cal H_3$. Since these three categories are pairwise inequivalent, the three Lagrangian algebras are pairwise non-isomorphic, and we can label them so that $Z(\cal H_3)_{L_i} \cong \cal H_i$.

By Lemma \ref{lem:rank} and part 3 of Proposition \ref{prop:GS}, we have
\begin{equation}
\label{eq:end-dims}
\dim\operatorname{Hom}(L_1,L_1) = 4,\quad \dim\operatorname{Hom}(L_2,L_2) = \dim\operatorname{Hom}(L_3,L_3) = 6.
\end{equation}
For the three objects of Lemma \ref{lem:possible-objects}, the dimensions of the endomorphism spaces are $1+1+4 = 6$ for $0 \oplus b \oplus 2a$ and $0 \oplus b \oplus 2c_0$, and $1+1+1+1 = 4$ for $0 \oplus b \oplus a \oplus c_0$, since $0$, $b$, $a$, $c_0$ are pairwise non-isomorphic simple objects. This proves the assertion.
\end{proof}

\begin{remark}\normalfont
The objects of $L_1$ and $L_2$ appear in \cite[pp.~18--19]{EvG11} as the dual canonical endomorphisms of the Longo--Rehren inclusions of the two even parts of the Haagerup subfactor, together with the three monomial modular invariants $|\chi_0+\chi_b+2\chi_a|^2$, $|\chi_0+\chi_b+\chi_a+\chi_{c_0}|^2$ and $|\chi_0+\chi_b+2\chi_{c_0}|^2$ \cite[Proposition 4 p.~32]{EvG11}. Proposition \ref{prop:which-line} below attaches the third invariant to the category $\cal H_3$, which was found later \cite{GrS12}.
\end{remark}

We now realize the Lagrangian algebras of objects $0 \oplus b \oplus 2a$ and $0 \oplus b \oplus 2c_0$ by means of a Haagerup algebra.

\begin{prop}
\label{prop:lagrangian-objects}
Let $U$ be a Haagerup algebra, and let $L_{\Lambda_a}$ and $L_{\Lambda_{c_0}}$ be the algebras in $\Rep(U)$ attached to the extensions $H_a$ and $H_{c_0}$ of \eqref{eq:Hx}.
\begin{enumerate}
\item $L_{\Lambda_a}$ and $L_{\Lambda_{c_0}}$ are Lagrangian algebras, isomorphic as objects to $0 \oplus b \oplus 2a$ and $0 \oplus b \oplus 2c_0$ respectively.
\item Under the braided equivalence $\Rep(U) \cong Z(\cal H_3)$, we have $\{L_{\Lambda_a},L_{\Lambda_{c_0}}\} = \{L_2,L_3\}$ and the remaining Lagrangian algebra $L_1$ is isomorphic to $0 \oplus b \oplus a \oplus c_0$ as an object.
\item The algebra $A_0$ attached to $V' = U \oplus M_b$ is a connected \'etale algebra, and $L_{\Lambda_a}$, $L_{\Lambda_{c_0}}$ are the only Lagrangian algebras of $\Rep(U)$ containing $A_0$. In particular $L_1$ does not contain $A_0$ as a subalgebra.
\end{enumerate}
\end{prop}

\begin{proof}
1. The vertex operator algebra $H_a$ is a simple extension of $U$ isomorphic to $V_{E_8}$, hence holomorphic and self-contragredient. By \cite[Theorems 3.2 and 3.4 pp.~10 and 16]{HKL15} and Lemma \ref{lem:intermediate-etale} applied to $W = B = H_a$, the algebra $L_{\Lambda_a}$ is a connected \'etale algebra in $\Rep(U)$ whose category of local modules is equivalent to $\Rep(H_a) \cong \mathrm{Vec}$. By \cite[Remark 4.3 p.~16]{DMNO13}, it is Lagrangian. Its object is read off from \eqref{eq:Hx-decomposition}. The same argument applies to $H_{c_0}$.

2. Since $a$ and $c_0$ are non-isomorphic, the objects of $L_{\Lambda_a}$ and $L_{\Lambda_{c_0}}$ are non-isomorphic, so these two Lagrangian algebras are distinct. By Proposition \ref{prop:three-lagrangian}, they have endomorphism spaces of dimension six, hence they are $L_2$ and $L_3$ in some order, and $L_1$ is the remaining one, with the stated object.

3. The algebra $A_0$ corresponds to the simple vertex subalgebra $V' \subset H_a$, so it is a connected \'etale subalgebra of $L_{\Lambda_a}$ by Lemma \ref{lem:intermediate-etale}. By \cite[Theorem 3.4 p.~16]{HKL15} and \cite[Theorem 3.65 p.~85]{CKM24}, the braided category $\Rep(U)_{A_0}^0$ is braided equivalent to $\Rep(V')$, which is the pointed modular category of the discriminant form $(\Z_3 \times \Z_3,q)$. By part 4 of Proposition \ref{prop:DMNO}, the connected \'etale algebras $L$ of $\Rep(U)$ containing $A_0$ correspond to the connected \'etale algebras $B$ of $\Rep(U)^0_{A_0} \cong \Rep(V')$. If $W$ is the corresponding extension of $U$ and of $V'$, then \cite[Theorems 3.2 and 3.4 pp.~10 and 16]{HKL15}, applied to $U \subset W$ and to $V' \subset W$, give isomorphisms of categories $\Rep(U)^0_L \cong \Rep(W) \cong \Rep(V')^0_B$. Hence $L$ is Lagrangian if and only if $B$ is, by \cite[Remark 4.3 p.~16]{DMNO13}, and the Lagrangian algebras of $\Rep(U)$ containing $A_0$ correspond to the holomorphic extensions of $V'$, by \cite[Theorems 3.2 and 3.4 pp.~10 and 16]{HKL15} and Lemma \ref{lem:intermediate-etale}. A holomorphic extension of the lattice vertex operator algebra $V'$ is a direct sum of simple currents $V'_w$ with $w$ in a subgroup $\Lambda$ of $\Z_3 \times \Z_3$ of order three on which $q$ vanishes, since its Frobenius--Perron dimension is three. Indeed, the simple summands of a connected \'etale algebra in the pointed category $\Rep(V')$ have trivial twist, appear with multiplicity at most one and form a subgroup, by the stabilizer argument of the proof of \cite[Proposition 4.1]{Arh26}. Since $q$ vanishes on $\Lambda$, the restriction of the pointed braided category $\Rep(V')$ to $\Lambda$ is equivalent to $\Rep(\Z_3)$ as a symmetric category, with trivial associator and braiding, so that an algebra structure on $\bigoplus_{w \in \Lambda}V'_w$ is given by a normalized $2$-cochain with values in $\Cbb^\times$ and trivial coboundary, the components of the multiplication being nonzero since the pairing of a connected \'etale algebra is nondegenerate \cite[Remark 3.4 p.~7]{DMNO13}. Hence this algebra structure is unique up to isomorphism since $H^2(\Z_3,\Cbb^\times) = 0$. The form $q(i,j) = \frac{2i^2+j^2}{3}$ vanishes exactly on the two lines generated by $(1,1)$ and $(1,2)$. Hence there are exactly two Lagrangian algebras of $\Rep(V')$, giving $H_a$ and $H_{c_0}$. The last statement follows from part 2.
\end{proof}

\begin{prop}
\label{prop:which-line}
In the labels of \cite{EvG11} recalled in Section \ref{subsec:center}, we have $L_2 \cong 0 \oplus b \oplus 2a$ and $L_3 \cong 0 \oplus b \oplus 2c_0$ as objects. If $U$ is a Haagerup algebra, then the extension $H_a$ of $U$ carries $\cal H_2$ while $H_{c_0}$ carries $\cal H_3$, that is, $L_{\Lambda_a} \cong L_2$ and $L_{\Lambda_{c_0}} \cong L_3$.
\end{prop}

\begin{proof}
By Lemma \ref{lem:rank} applied to $\cal D = \cal H_2$, the Lagrangian algebra $L_2$ corresponds to $I(1)$ under $Z(\cal H_3) \cong Z(\cal H_2)$. By adjunction, the multiplicity of a simple object $z$ of $Z(\cal H_2)$ in $I(1)$ is $\dim\operatorname{Hom}_{\cal H_2}(F(z),1)$, where $F$ is the forgetful functor. By Section \ref{subsec:center}, this multiplicity is $1$ for $z \in \{0,b\}$, $2$ for $z = a$ and $0$ for the other simple objects. Hence $L_2 \cong 0 \oplus b \oplus 2a$. This is the computation of \cite[p.~18]{EvG11}.

We now compute $L_3$. Let $\Phi \co Z(\cal H_3) \to Z(\cal H_2)$ be a braided equivalence with $\Phi(L_2) \cong I(1)$, given by Lemma \ref{lem:rank}. For $i \in \{1,2,3\}$, the algebra $\Phi(L_i)$ is a Lagrangian algebra of $Z(\cal H_2)$ with category of modules equivalent to $\cal H_i$. By part 2 of Proposition \ref{prop:DMNO}, it is the Lagrangian algebra attached in \cite[Section 4.2 p.~17]{DMNO13} to an indecomposable $\cal H_2$-module category $\cal N_i$, and $\cal N_1$, $\cal N_2$, $\cal N_3$ are pairwise inequivalent, since the algebras $\Phi(L_i)$ are pairwise non-isomorphic, and they are, up to equivalence, all the indecomposable $\cal H_2$-module categories (part 2 of Proposition \ref{prop:GS}). Here $\cal N_2$ is the regular module category, since $\Phi(L_2) \cong I(1)$ and $I(1)$ is the Lagrangian algebra attached to $\cal H_2$ as a module category over itself \cite[Remark 4.9 p.~17]{DMNO13}. We equip $\cal H_2$ with its spherical structure with positive dimensions \cite[Proposition 8.23 p.~34]{ENO05}. The induced dimensions of the simple objects of $Z(\cal H_2)$ are positive, so that the twists of the Lagrangian algebras $\Phi(L_i)$ are trivial \cite[Remark 2.19 p.~8]{Ost15}. Hence, by \cite[Remark 2.19 pp.~7--8]{Ost15}, the number of simple objects of the category $\operatorname{Fun}_{\cal H_2}(\cal N_2,\cal N_i)$ of module functors is $\dim\operatorname{Hom}(\Phi(L_2),\Phi(L_i)) = \dim\operatorname{Hom}(L_2,L_i)$. Since a module functor from the regular module category is determined by the image of the unit object, $\operatorname{Fun}_{\cal H_2}(\cal N_2,\cal N_i)$ is equivalent to $\cal N_i$, and $\dim\operatorname{Hom}(L_2,L_i)$ is the rank of $\cal N_i$. By Proposition \ref{prop:three-lagrangian} and the first paragraph, $\dim\operatorname{Hom}(L_2,L_1) = \dim\operatorname{Hom}(0 \oplus b \oplus 2a,0 \oplus b \oplus a \oplus c_0) = 4$, so that $\cal N_1$ has rank four.

Let $A = 1 \oplus \alpha \oplus \alpha^2$ with the structure of $\Q$-system in $\cal H_2$ used to define $\cal H_3$ \cite[Lemma 3.14 p.~21 and Definition 3.19 p.~23]{GrS12}, and let $\cal M$ be the category of right $A$-modules in $\cal H_2$, which is a semisimple left $\cal H_2$-module category since $A$ is separable \cite[Proposition 2.7 p.~5]{DMNO13}. For objects $Y$, $Y'$ of $\cal H_2$, we have $\operatorname{Hom}_A(Y \ot A,Y' \ot A) \cong \operatorname{Hom}(Y,Y' \ot A)$, and since $A$ is a $\Q$-system, every simple right $A$-module is a summand of some $Y \ot A$ with $Y$ simple. Since $\xi\alpha \cong \alpha^2\xi$ and $\xi\alpha^2 \cong \alpha\xi$ by \eqref{eq:H2-rules}, we have $\alpha^j\xi \ot A \cong \alpha^j\xi \oplus \alpha^{j+1}\xi \oplus \alpha^{j+2}\xi$ for $j \in \{0,1,2\}$. Hence $\dim\operatorname{Hom}_A(Y \ot A,Y' \ot A)$ is equal to $1$ if $Y$ and $Y'$ belong both to $\{1,\alpha,\alpha^2\}$ or both to $\{\xi,\alpha\xi,\alpha^2\xi\}$, and to $0$ otherwise. We obtain that the simple objects of $\cal M$ are $A \cong \alpha \ot A \cong \alpha^2 \ot A$ and $\xi \ot A \cong \alpha\xi \ot A \cong \alpha^2\xi \ot A$, up to isomorphism. Since every simple object of $\cal M$ is a summand of $Y \ot A$ for some $Y$, the module category $\cal M$ is generated by $A$, hence indecomposable. It has rank two, so that it is equivalent neither to $\cal N_2$, of rank six, nor to $\cal N_1$, of rank four. Hence $\cal M \cong \cal N_3$ and $\dim\operatorname{Hom}(L_2,L_3) = 2$. Since $\dim\operatorname{Hom}(0 \oplus b \oplus 2a,0 \oplus b \oplus 2a) = 6$, Proposition \ref{prop:three-lagrangian} gives $L_3 \cong 0 \oplus b \oplus 2c_0$.

Finally, let $U$ be a Haagerup algebra. By part 1 of Proposition \ref{prop:lagrangian-objects}, $L_{\Lambda_a}$ and $L_{\Lambda_{c_0}}$ are Lagrangian algebras isomorphic to $0 \oplus b \oplus 2a$ and $0 \oplus b \oplus 2c_0$ as objects. By the above and Proposition \ref{prop:three-lagrangian}, $L_2$ and $L_3$ are the only Lagrangian algebras with these objects, so that $L_{\Lambda_a} \cong L_2$ and $L_{\Lambda_{c_0}} \cong L_3$.
\end{proof}

\begin{remark}\normalfont
\label{rem:which-line}
The conventions in \cite[Remark 2.1]{Arh26} agree with Proposition \ref{prop:which-line}. In the tube-algebra labelling of \cite[(A8) p.~7]{AHKZ26}, the Lagrangian algebra of $\cal H_3$ is $1 \oplus \pi_1 \oplus 2\pi_2$. This is the category $\cal H_3$ of \cite{GrS12}, since the algebra $1 \oplus \rho \oplus \alpha\rho$ relates it to $\cal H_1$ at index $4+\sqrt{13}$ \cite[p.~4]{AHKZ26}, \cite[Theorems 1.3 and 1.4 pp.~3--4]{GrS12}. Comparing its Lagrangian algebra with $0 \oplus b \oplus 2c_0$ gives $\pi_2 = c_0$ and $\sigma_1 = a$, as in \cite[Remark 2.1]{Arh26}. The transposition of $a$ and $c_0$ preserves the modular data of Section \ref{subsec:center}.
\end{remark}

\begin{cor}
\label{cor:categorical}
The three Lagrangian algebras of $Z(\cal H_3)$ are isomorphic, as objects, to
\[
0 \oplus b \oplus a \oplus c_0,\quad 0 \oplus b \oplus 2a,\quad 0 \oplus b \oplus 2c_0,
\]
with categories of modules $\cal H_1$, $\cal H_2$ and $\cal H_3$ respectively.
\end{cor}

\begin{proof}
This is Proposition \ref{prop:three-lagrangian} with Proposition \ref{prop:which-line}.
\end{proof}

%%%%%%%%%%%%%%%%%%%%%%%%%%%%%%%%%%%%%%%%%%%%%%%%%%%%%%%%%%%%%%%%%%%%%%%%%%%%%
\subsection{Holomorphic extensions of a Haagerup algebra}

\begin{thm}
\label{thm:lagrangian-extensions}
Let $U$ be a Haagerup algebra.
\begin{enumerate}
\item Up to isomorphism of extensions, $U$ has exactly three holomorphic extensions $B_1$, $B_2$, $B_3$. The algebra $L_i$ in $\Rep(U)$ attached to $B_i$ is the Lagrangian algebra of Proposition \ref{prop:three-lagrangian}, so that $\Rep(U)_{L_i} \cong \cal H_i$. Each $B_i$ is a strongly rational holomorphic vertex operator algebra isomorphic to $V_{E_8}$.
\item $B_2 = H_a$ and $B_3 = H_{c_0}$ as extensions of $U$.
\item The subspace $U \oplus M_b$ of $B_1$ is not a vertex subalgebra. The only vertex subalgebras of $B_1$ containing $U$ are $U$ and $B_1$. For $i \in \{2,3\}$, the only simple vertex subalgebras of $B_i$ containing $U$ are $U$, $V'$ and $B_i$.
\end{enumerate}
\end{thm}

\begin{proof}
1. Let $B$ be a holomorphic extension of $U$. Since $\Rep(B)$ is trivial, the contragredient module $B'$ is isomorphic to $B$. By \cite[Theorems 3.2 and 3.4 pp.~10 and 16]{HKL15} and Lemma \ref{lem:intermediate-etale}, the algebra $L_B$ attached to $B$ is a connected \'etale algebra with $\Rep(U)_{L_B}^0 \cong \Rep(B) \cong \mathrm{Vec}$, hence a Lagrangian algebra. Conversely, a Lagrangian algebra $L$ of $\Rep(U)$ is a haploid commutative algebra with trivial twist (Lemma \ref{lem:possible-objects} and its proof), so it defines an extension $B_L$ of $U$ with $\Rep(B_L) \cong \Rep(U)^0_L \cong \mathrm{Vec}$ by \cite[Theorems 3.2 and 3.4 pp.~10 and 16]{HKL15}. Since $L$ is connected, $L$ is a simple object of $\Rep(U)_L$ \cite[p.~10]{DMNO13}, and the ideals of $B_L$ correspond to the subobjects of $L$ in $\Rep(U)_L$ \cite[proof of Theorem 3.2 p.~10]{HKL15}. Hence $B_L$ is simple. These two constructions are inverse to each other up to isomorphism, by \cite[Theorem 3.2 p.~10]{HKL15}. By Proposition \ref{prop:three-lagrangian}, there are exactly three Lagrangian algebras, which proves the first two assertions. The nonvacuum summands of $L_i$ are $b$, $a$ and $c_0$, whose lowest conformal weights are equal to one \cite[(4.1)]{Arh26}. Hence $B_i$ is of CFT type. It is self-contragredient since $\Rep(B_i) \cong \mathrm{Vec}$, and it is $C_2$-cofinite and rational by the argument of the proof of \cite[Corollary 5.13 p.~53]{CKM22}, since its category of modules is semisimple, see also \cite[proof of Proposition 4.1]{Arh26}. A strongly rational holomorphic vertex operator algebra of central charge eight is isomorphic to $V_{E_8}$ \cite[Theorem 1 p.~2]{DoM04b}.

2. This is Proposition \ref{prop:which-line}.

3. Let $W$ be a simple vertex subalgebra of $B_i$ containing $U$. By Lemma \ref{lem:intermediate-etale}, the corresponding algebra $A_W$ is a connected \'etale subalgebra of $L_i$. By Lemma \ref{lem:rank}, under the braided equivalence $\Rep(U) \cong Z(\cal H_i)$ given by part 1 of Proposition \ref{prop:DMNO}, the algebra $L_i$ corresponds to $I_{\cal H_i}(1)$. By part 3 of Proposition \ref{prop:DMNO}, the connected \'etale subalgebras of $L_i$ correspond to the fusion subcategories of $\cal H_i$, with $\operatorname{FPdim}(A_W)\operatorname{FPdim}(\cal B) = \dim \cal H_i$ if $A_W$ corresponds to $\cal B$. By part 4 of Proposition \ref{prop:GS}, the connected \'etale subalgebras of $L_1$ are the unit object $0$ and $L_1$, and for $i \in \{2,3\}$ the connected \'etale subalgebras of $L_i$ are $0$, $L_i$ and a unique subalgebra of Frobenius--Perron dimension $\frac{3(3d+2)}{3} = 3d+2$, corresponding to $\cal P_i$. For $i \in \{2,3\}$, the extension $B_i$ contains $V'$, which is simple, so this subalgebra is $A_0$. This proves the statements on simple vertex subalgebras.

Let now $W$ be any vertex subalgebra of $B_1$ containing $U$, and let $(\cdot,\cdot)$ be the invariant bilinear form of $B_1$. It is $U$-invariant, so it pairs the isotypic component of type $z$ with that of type $z^*$, and its restriction to each isotypic component of a self-dual type is nondegenerate. Since $L_1 \cong 0 \oplus b \oplus a \oplus c_0$ is multiplicity free, $W$ is a sum of isotypic components of $B_1$, all of self-dual types. Hence the restriction of $(\cdot,\cdot)$ to $W$ is nondegenerate, $W$ is simple by Lemma \ref{lem:intermediate-etale}, and the description of the connected \'etale subalgebras of $L_1$ obtained above gives $W \in \{U,B_1\}$. In particular, $U \oplus M_b$ is not a vertex subalgebra of $B_1$.
\end{proof}

\begin{remark}\normalfont
By part 3 of Theorem \ref{thm:lagrangian-extensions}, there is no vertex subalgebra strictly between $U$ and $B_1 \cong V_{E_8}$. The weight-one space of $B_1$ is $U_1 \oplus R_b \oplus R_a \oplus R_{c_0}$, where $R_z$ is the lowest weight space of $M_z$, and its dimension is $(6+13\gamma)+(80-13\gamma)+81+81 = 248 = \dim \mathfrak e_8$ by \cite[(4.1)]{Arh26}, for both values $\gamma \in \{0,1\}$ of the parameter of loc.\ cit. In fact $\gamma = 0$ by \cite[Theorem 7.1]{Arh26}.
\end{remark}

\begin{remark}\normalfont
\label{rem:symmetries}
Let $L$ be a Lagrangian algebra in $\Rep(U)$ with holomorphic extension $B$. In the language of symmetry topological field theories, the fusion category $\Rep(U)_L \cong \cal H_i$ of $L$-modules is interpreted as the category of the topological defect lines of $B$ which preserve $U$, see \cite{AHKZ26,GaR26,Hua26}. Theorem \ref{thm:lagrangian-extensions} thus gives a vertex algebraic version of the suggestion of \cite{AHKZ26} that the chiral $(\mathfrak e_8)_1$ theory carries the three symmetries $\cal H_1$, $\cal H_2$ and $\cal H_3$, under the assumption that a Haagerup algebra exists. An operator algebraic version is given in Section \ref{sec:N}, under $\mathrm{(H_J)}$.
\end{remark}

%%%%%%%%%%%%%%%%%%%%%%%%%%%%%%%%%%%%%%%%%%%%%%%%%%%%%%%%%%%%%%%%%%%%%%%%%%%%%
\section{The conformal net of \texorpdfstring{$S$}{S}}
\label{sec:net}

Throughout this section, $L = E_6 \oplus A_2$ denotes the root lattice and $V' = V_L = M(1) \ot \Cbb_\epsi[L]$ is the lattice vertex operator algebra of \cite[Section 10]{Arh26}, built as in \cite[Chapter 8]{FLM88} with a bimultiplicative cocycle $\epsi \co L \times L \to \{\pm 1\}$ satisfying $\epsi(\alpha,\beta)\epsi(\beta,\alpha)^{-1} = (-1)^{\langle \alpha,\beta\rangle}$ and $\epsi(\alpha,\alpha) = (-1)^{\langle\alpha,\alpha\rangle/2}$ for any $\alpha, \beta \in L$. The second property holds for the cocycle of \cite[Section 10.1]{Arh26}, which satisfies it on the simple roots, since $\alpha \mapsto \epsi(\alpha,\alpha)(-1)^{\langle\alpha,\alpha\rangle/2}$ is a group homomorphism $L \to \{\pm 1\}$ by the first property. We denote by $L^\circ$ the dual lattice, so that $L^\circ/L \cong \Z_3 \times \Z_3$ with the quadratic form $q$ of Section \ref{subsec:haagerup-algebras}. No hypothesis on Haagerup algebras is made in this section.

%%%%%%%%%%%%%%%%%%%%%%%%%%%%%%%%%%%%%%%%%%%%%%%%%%%%%%%%%%%%%%%%%%%%%%%%%%%%%
\subsection{The lattice extensions and their unitary structures}
\label{subsec:lattice-extensions}

Let $\Lambda \subset L^\circ/L$ be one of the two isotropic lines and let $L_\Lambda \subset L^\circ$ be its preimage. Since $q$ vanishes on $\Lambda$, the lattice $L_\Lambda$ is even, and it is unimodular since $[L_\Lambda : L] = 3$ and $|L^\circ/L| = 9$. Hence $L_\Lambda$ is isometric to $E_8$.

\begin{lemma}
\label{lem:cocycle}
\begin{enumerate}
\item The cocycle $\epsi$ extends uniquely to a bimultiplicative map $\epsi_\Lambda \co L_\Lambda \times L_\Lambda \to \{\pm 1\}$. We have $\epsi_\Lambda(\alpha,\beta)\epsi_\Lambda(\beta,\alpha)^{-1} = (-1)^{\langle \alpha,\beta\rangle}$ and $\epsi_\Lambda(\alpha,\alpha) = (-1)^{\langle\alpha,\alpha\rangle/2}$ for any $\alpha,\beta \in L_\Lambda$.
\item The lattice vertex operator algebra
\[
H_\Lambda \ov{\mathrm{def}}{=} M(1) \ot \Cbb_{\epsi_\Lambda}[L_\Lambda] = \bigoplus_{w \in \Lambda} V'_w,\quad V'_w \ov{\mathrm{def}}{=} M(1) \ot \Cbb_{\epsi_\Lambda}[L+\lambda_w],
\]
where $\lambda_w \in L^\circ$ is a representative of $w$, contains $V' = V'_0$ as a vertex subalgebra with the same conformal vector, and $H_\Lambda$ is isomorphic to $V_{E_8}$.
\end{enumerate}
\end{lemma}

\begin{proof}
1. A bimultiplicative map $M \times M \to \{\pm 1\}$ on a lattice $M$ is the same thing as a bilinear form on the $\mathbb F_2$-vector space $M/2M$. The inclusion $L \subset L_\Lambda$ induces an isomorphism $L/2L \to L_\Lambda/2L_\Lambda$. Indeed, both spaces have dimension eight, and if $x \in L$ satisfies $x = 2y$ with $y \in L_\Lambda$, then $3y \in L$ and $2y \in L$, hence $y \in L$. This gives the existence and uniqueness of $\epsi_\Lambda$. The maps $(\alpha,\beta) \mapsto \epsi_\Lambda(\alpha,\beta)\epsi_\Lambda(\beta,\alpha)^{-1}$ and $(\alpha,\beta) \mapsto (-1)^{\langle \alpha,\beta\rangle}$ are bimultiplicative on $L_\Lambda$, since $L_\Lambda$ is integral, and they coincide on $L$. By the same isomorphism, they coincide on $L_\Lambda$. Finally, the map $s(\alpha) = \epsi_\Lambda(\alpha,\alpha)(-1)^{\langle\alpha,\alpha\rangle/2}$ satisfies $s(\alpha+\beta) = s(\alpha)s(\beta)$ by the commutator identity, so it is a group homomorphism $L_\Lambda \to \{\pm 1\}$. It is trivial on $L$, and $[L_\Lambda:L] = 3$ is odd, so that $s$ is trivial.

2. By part 1, the construction of lattice vertex operator algebras \cite[Chapter 8]{FLM88} applies to $(L_\Lambda,\epsi_\Lambda)$, and the subspace $M(1) \ot \Cbb_{\epsi}[L]$ is a vertex subalgebra equal to $V'$ with the same conformal vector. Up to isomorphism, the lattice vertex operator algebra does not depend on the choice of the cocycle, and $L_\Lambda \cong E_8$, so that $H_\Lambda \cong V_{E_8}$.
\end{proof}

The subspace $V'_w$ of $H_\Lambda$ realizes the simple $V'$-module attached to $w$. An extension of $V'$ whose underlying $V'$-module is $\bigoplus_{w \in \Lambda}V'_w$ is unique up to isomorphism of extensions, by the argument of the proof of part 3 of Proposition \ref{prop:lagrangian-objects}, since $H^2(\Z_3,\Cbb^\times) = 0$.

Hence, in the situation of Section \ref{subsec:haagerup-algebras}, the vertex operator algebras $H_a$ and $H_{c_0}$ of \eqref{eq:Hx} are isomorphic, as extensions of $V'$, to $H_{\Lambda_a}$ and $H_{\Lambda_{c_0}}$. For a character $\chi \co L \to \R_{>0}$, we write $\chi(\lambda) = \e^{\langle h_\chi,\lambda\rangle}$ with $h_\chi$ in the real span of $L$, and we extend $\chi$ to $L^\circ$ by the same formula. Let $\phi_{1}$ be the antilinear map of $H_\Lambda$ defined by \cite[(10.3)]{Arh26} for the trivial character, that is
\[
\phi_1\big(c  h_1(-n_1)\cdots h_k(-n_k)\e^\lambda\big) = (-1)^k  \ovl{c}  h_1(-n_1)\cdots h_k(-n_k)\e^{-\lambda},\quad c \in \Cbb,\ \lambda \in L_\Lambda,
\]
where $h_1, \ldots, h_k$ belong to the real span of $L$.

\begin{lemma}
\label{lem:unitary-extension}
\begin{enumerate}
\item The map $\phi_1$ is an antilinear involutive automorphism of $H_\Lambda$. If $\beta$ denotes the normalized invariant bilinear form of $H_\Lambda$, then $H_1(u,v) \ov{\mathrm{def}}{=} \beta(u,\phi_1(v))$ is a normalized invariant scalar product on $H_\Lambda$.
\item Let $\{\cdot|\cdot\}$ be a normalized invariant scalar product on $V'$. There exists $x$ in $(V')_1$ such that $\tilde h = \exp(x(0))$ is an automorphism of $H_\Lambda$ preserving $V'$, $\{u|v\} = H_1(\tilde hu,\tilde hv)$ for any $u,v \in V'$ and $H_1(u,\tilde hu) > 0$ for any nonzero $u \in V'$. Consequently, $\{u|v\}_\Lambda \ov{\mathrm{def}}{=} H_1(\tilde hu,\tilde hv)$ is a normalized invariant scalar product on $H_\Lambda$ which extends $\{\cdot|\cdot\}$, and $V'$ is a unitary subalgebra of $(H_\Lambda,\{\cdot|\cdot\}_\Lambda)$. The scalar product $\{\cdot|\cdot\}_\Lambda$ does not depend on the choice of $x$.
\item For the unitary structure $\phi_\chi$ of \cite[Theorem 10.10]{Arh26}, one can take $x = \frac12 h_\chi(-1)\mathbf 1$ in part 2. The corresponding PCT operator of $H_\Lambda$ is $\phi_\chi = \phi_1 \circ \exp(h_\chi(0))$.
\end{enumerate}
\end{lemma}

\begin{proof}
1. The cocycle $\epsi_\Lambda$ is bimultiplicative, so that $\epsi_\Lambda(-\alpha,-\beta) = \epsi_\Lambda(\alpha,\beta)$, and $\epsi_\Lambda(\alpha,\alpha) = (-1)^{\langle\alpha,\alpha\rangle/2}$ by Lemma \ref{lem:cocycle}. Hence $\e^\alpha\e^{-\alpha} = \epsi_\Lambda(\alpha,-\alpha)\e^0 = (-1)^{\langle\alpha,\alpha\rangle/2}\e^0$ for any $\alpha \in L_\Lambda$, which is the relation used in \cite[Lemma 4.11 p.~23]{DoL14}, and $\phi_1$ is the antilinear involution for which lattice vertex operator algebras are unitary \cite[Theorem 4.12 p.~24]{DoL14}, \cite[Example 5.9 p.~35]{CKLW18}, as in the proof of \cite[Theorem 10.10]{Arh26}.

2. By part 1, the restriction of $H_1$ to $V'$ is a normalized invariant scalar product on $V'$. Since $V'$ is simple, there exists by \cite[Proposition 5.19 p.~40]{CKLW18} an automorphism $h$ of $V'$ with $\{u|v\} = H_1(hu,hv)$ for any $u,v \in V'$ and $H_1(u,hu) > 0$ for any nonzero $u \in V'$. In particular $h$ preserves each $(V')_n$ and is diagonalizable with positive eigenvalues. Let $D$ be the linear map acting by $\log \lambda$ on the eigenspace of $h$ for the eigenvalue $\lambda$. Since $h(u_{(n)}v) = \lambda\mu  u_{(n)}v$ if $hu = \lambda u$ and $hv = \mu v$, the map $D$ is a derivation of $V'$ preserving the grading. Its restriction to $(V')_1$ is a derivation of the semisimple Lie algebra $((V')_1,[u,v] = u_{(0)}v) \cong \mathfrak e_6 \oplus \mathfrak{sl}_3$, so that there is $x \in (V')_1$ with $D(u) = x_{(0)}u$ for $u \in (V')_1$. The derivations $D$ and $x_{(0)} = x(0)$ of $V'$ coincide on $(V')_1$ and vanish on the vacuum. The space $(V')_1$ contains the vectors $h(-1)\mathbf 1$ and $\e^\alpha$ for the roots $\alpha$ of $L$, the lattice $L$ is spanned by its roots, and $\e^{\alpha+\beta}$ is a nonzero multiple of $\e^\alpha_{(-\langle\alpha,\beta\rangle-1)}\e^\beta$ \cite[Chapter 8]{FLM88}. Hence $(V')_1$ generates $V'$, and $D = x(0)$. Now $x(0)$ is a derivation of $H_\Lambda$ which preserves the finite-dimensional homogeneous spaces, so that $\tilde h = \exp(x(0))$ is an automorphism of $H_\Lambda$ which extends $h = \exp(D)$. The transport of a normalized invariant scalar product by an automorphism is a normalized invariant scalar product, and the PCT operator of $\{\cdot|\cdot\}_\Lambda$ is $\tilde h^{-1}\phi_1\tilde h$, which preserves $V'$. Since $L_1V' \subset V'$, the last assertion follows from \cite[Proposition 5.23 p.~43]{CKLW18}. Two choices of $x$ give automorphisms of $H_\Lambda$ which coincide on $V'$, since a positive square root is unique. They differ by an automorphism of $H_\Lambda$ fixing $V'$ pointwise, which acts on each simple $V'$-module $V'_w$, $w \in \Lambda$, by a scalar, hence by a power of the automorphism $g_\Lambda$ of Proposition \ref{prop:lattice-nets}. This automorphism is unitary for $H_1$ and commutes with $\tilde h$, so that $\{\cdot|\cdot\}_\Lambda$ does not depend on the choice of $x$.

3. With $g = \exp(h_\chi(0))$, which acts by $\chi(\lambda)$ on $M(1) \ot \e^\lambda$, we have $\phi_\chi = \phi_1 \circ g = g^{-1/2} \circ \phi_1 \circ g^{1/2}$ on $V'$ by the proof of \cite[Theorem 10.10]{Arh26}, and the same computation holds on $H_\Lambda$. We still denote by $H_\chi(u,v) = H_1(g^{1/2}u,g^{1/2}v)$ the corresponding scalar product on $H_\Lambda$, which extends the scalar product $H_\chi$ of $V'$ of \cite[Theorem 10.10]{Arh26}.
\end{proof}

\begin{prop}
\label{prop:lattice-nets}
Let $\Lambda$ be an isotropic line and let $\{\cdot|\cdot\}$ be a normalized invariant scalar product on $H_\Lambda$.
\begin{enumerate}
\item The unitary vertex operator algebra $(H_\Lambda,\{\cdot|\cdot\})$ is strongly local and has bounded localized vertex operators in the sense of \cite[Definition 4.12 p.~42]{Ten19a}. The net of \cite{Ten19a} and the net of \cite{CKLW18} of $(H_\Lambda,\{\cdot|\cdot\})$ are isomorphic, and they are isomorphic to the loop group net $\cal A_{E_8}$ of $E_8$ at level one. This net is completely rational and $\mu_{\cal A_{E_8}} = 1$.
\item Assume that $\{\cdot|\cdot\}$ is one of the scalar products of part 2 of Lemma \ref{lem:unitary-extension}. Let $\beta_\Lambda \in L^\circ$ with $\langle \beta_\Lambda,\lambda_w\rangle \in \frac13+\Z$ for a generator $w$ of $\Lambda$, and put $g_\Lambda \ov{\mathrm{def}}{=} \e^{2\pi\i \beta_\Lambda(0)}$. Then $g_\Lambda$ is a unitary automorphism of order three of $H_\Lambda$ with fixed point algebra $V'$, and $\cal A_{V'} = \cal A_{H_\Lambda}^{\langle g_\Lambda\rangle}$. The net $\cal A_{V'}$ is completely rational, $\mu_{\cal A_{V'}} = 9$, $[\cal A_{H_\Lambda} : \cal A_{V'}] = 3$ and $\cal A_{V'}(I)' \cap \cal A_{H_\Lambda}(I) = \Cbb$ for any interval $I$. Moreover $\cal A_{V'}$ is isomorphic to $\cal A_{E_6,1} \ot \cal A_{SU(3),1}$.
\end{enumerate}
\end{prop}

\begin{proof}
1. The lattice vertex operator algebra of $E_8$ coincides with $V(E_8,1)$ \cite[Proof of Theorem 7.4 p.~71]{Ten19a}. With its unitary structure, it has bounded localized vertex operators \cite[Theorem 7.4 p.~71]{Ten19a}, it is strongly local with net $\cal A_{E_8}$ \cite[Example 8.7 p.~63]{CKLW18}, and the nets of \cite{Ten19a} and \cite{CKLW18} are isomorphic \cite[Corollary 7.7 p.~72]{Ten19a}. Let $\psi \co H_\Lambda \to V(E_8,1)$ be an isomorphism of vertex operator algebras. The pull back by $\psi$ of the scalar product of $V(E_8,1)$ is a normalized invariant scalar product on $H_\Lambda$. By \cite[Proposition 5.19 p.~40]{CKLW18}, it is of the form $\{hu|hv\}$ for some automorphism $h$ of $H_\Lambda$. Then $\psi \circ h^{-1}$ is a unitary isomorphism from $(H_\Lambda,\{\cdot|\cdot\})$ onto $V(E_8,1)$, and all the properties are transported. Here bounded localized vertex operators refer to the system of generalized annuli fixed in \cite[Section 7]{Ten19a}. A unitary isomorphism of vertex operator algebras preserves the conformal vectors, hence intertwines the representations of the Virasoro algebra and the operators attached to the annuli, so that it transports these properties and the nets. Finally, by \cite[Proposition 7.5 p.~72]{Ten19a}, the net of \cite{Ten19a} of $V(E_8,1)$ is isomorphic to the lattice net of $E_8$ of Dong and Xu, since $E_8$ is the code lattice of the extended Hamming code \cite[Example 7.2 p.~70]{Ten19a}. By \cite[Corollary 3.19 p.~14]{DoX06}, this net is completely rational with $\mu$-index $|E_8^\circ/E_8| = 1$.

2. The map $g_\Lambda$ acts on $M(1) \ot \e^\lambda$ by the scalar $\e^{2\pi\i\langle \beta_\Lambda,\lambda\rangle}$, which is $1$ on $V'$ and a primitive cube root of unity on $V'_w$ for $w \neq 0$. Such an element $\beta_\Lambda$ exists since the pairing $L^\circ/L_\Lambda \times L_\Lambda/L \to \mathbb{Q}/\Z$ induced by $\langle\cdot,\cdot\rangle$ is perfect. Hence $g_\Lambda$ is an automorphism of order three with fixed point algebra $V'$. The sectors $M(1) \ot \e^\lambda$ are pairwise orthogonal for $H_1$, since $\beta$ pairs the sector of $\lambda$ with the sector of $-\lambda$ and $\phi_1$ exchanges them. So $g_\Lambda$ is unitary for $H_1$. The automorphism $\tilde h = \exp(x(0))$, $x \in (V')_1$, of Lemma \ref{lem:unitary-extension} commutes with $g_\Lambda$, since $g_\Lambda x(0)g_\Lambda^{-1} = (g_\Lambda x)(0) = x(0)$. Hence $g_\Lambda$ is also unitary for $\{\cdot|\cdot\}$. By \cite[Theorem 6.9 p.~53 and Proposition 7.6 p.~58]{CKLW18}, $g_\Lambda$ defines a vacuum preserving automorphism of the net $\cal A_{H_\Lambda}$ and $\cal A_{V'} = \cal A_{H_\Lambda}^{\langle g_\Lambda\rangle}$. The group $\langle g_\Lambda\rangle$ acts faithfully and properly in the sense of \cite[Definition 2.1 p.~5]{DoX06}. Indeed, if $g \in \langle g_\Lambda\rangle$ acts trivially on $\cal A_{H_\Lambda}(I)$ for some interval $I$, then $gx\Omega = x\Omega$ for any $x \in \cal A_{H_\Lambda}(I)$, so that $g = 1$ since $\Omega$ is cyclic for $\cal A_{H_\Lambda}(I)$. By part 1 and \cite[Theorem 2.4 p.~6]{DoX06}, the net $\cal A_{V'}$ is completely rational with $\mu_{\cal A_{V'}} = 3^2\mu_{\cal A_{H_\Lambda}} = 9$.

The vertex operator algebra $V'$ contains the conformal vector of $H_\Lambda$, so that $\cal A_{V'}(I)' \cap \cal A_{H_\Lambda}(I) = \Cbb$ by Lemma \ref{lem:trivial-coset} below. In particular $\cal A_{V'}(I) \subset \cal A_{H_\Lambda}(I)$ is an irreducible subfactor. By part 2 of Proposition \ref{prop:nets}, its index is finite, and by part 1 of the same proposition, $[\cal A_{H_\Lambda}:\cal A_{V'}]^2 = \mu_{\cal A_{V'}}/\mu_{\cal A_{H_\Lambda}} = 9$. Finally $V' \cong V_{E_6} \ot V_{A_2}$, and the Frenkel--Kac construction gives $V_{E_6} \cong V_{\mathfrak e_6,1}$ and $V_{A_2} \cong V_{\mathfrak{sl}_3,1}$. Transporting the scalar product along this isomorphism as in part 1, the last assertion follows from parts 1 and 5 of Proposition \ref{prop:CKLW}.
\end{proof}

%%%%%%%%%%%%%%%%%%%%%%%%%%%%%%%%%%%%%%%%%%%%%%%%%%%%%%%%%%%%%%%%%%%%%%%%%%%%%
\subsection{Strong locality of \texorpdfstring{$S$}{S}}
\label{subsec:strong-locality}

\begin{lemma}
\label{lem:trivial-coset}
Let $V$ be a strongly local simple unitary vertex operator algebra and let $W$ be a unitary subalgebra of $V$ containing the conformal vector of $V$. Then for any interval $I$ we have
\[
\cal A_W(I)' \cap \cal A_V(I) = \cal A_W(S^1)' \cap \cal A_V(I) = \Cbb,
\]
where $\cal A_W(S^1)$ is the von Neumann algebra generated by the algebras $\cal A_W(J)$.
\end{lemma}

\begin{proof}
The first equality holds for any covariant subnet of an irreducible conformal net \cite[Example 3.7 p.~24]{CKLW18}, \cite{Kos04}. By \cite[Proposition 7.8 p.~58]{CKLW18}, the coset subnet $\cal A_W(S^1)' \cap \cal A_V(I)$ is $\cal A_{W^c}(I)$, where $W^c$ is the coset of $W$ in $V$. By \cite[Proposition 5.31 p.~45]{CKLW18}, the conformal vector $\omega$ of $V$ decomposes as $\omega = \omega^W+\omega^{W^c}$. Here $\omega^W = \omega$ since $\omega \in W$, so that $\omega^{W^c} = 0$ and $W^c = \Cbb\Omega$ by \cite[Remark 5.30 p.~45]{CKLW18}. Hence $\cal A_{W^c}(I) = \Cbb$.
\end{proof}

\begin{thm}
\label{thm:strong-locality}
Let $W \subset V'$ be a vertex subalgebra containing the conformal vector of $V'$ which is a unitary subalgebra for a normalized invariant scalar product $\{\cdot|\cdot\}$ on $V'$. Let $\Lambda$ be an isotropic line and equip $H_\Lambda$ with the scalar product $\{\cdot|\cdot\}_\Lambda$ of Lemma \ref{lem:unitary-extension}.
\begin{enumerate}
\item $W$ is strongly local, and $\cal A_W \subset \cal A_{V'} \subset \cal A_{H_\Lambda}$ are covariant subnets, where $\cal A_{H_\Lambda} \cong \cal A_{E_8}$ and $[\cal A_{H_\Lambda} : \cal A_{V'}] = 3$.
\item The coset subnets of $\cal A_W$ in $\cal A_{V'}$ and in $\cal A_{H_\Lambda}$ are trivial, and the inclusions are irreducible: $\cal A_W(I)' \cap \cal A_{H_\Lambda}(I) = \Cbb$ and $\cal A_{V'}(I)' \cap \cal A_{H_\Lambda}(I) = \Cbb$ for any interval $I$.
\end{enumerate}
In particular, all these statements hold for the algebra $S$ of \cite[Theorem 10.10]{Arh26}, with the scalar product $H_\chi$.
\end{thm}

\begin{proof}
By Lemma \ref{lem:unitary-extension}, the vertex operator algebra $V'$ is a unitary subalgebra of $(H_\Lambda,\{\cdot|\cdot\}_\Lambda)$. Since $W$ contains the conformal vector, $L_1W \subset W$, and the PCT operator of $H_\Lambda$ restricts to that of $V'$, which preserves $W$. Hence $W$ is a unitary subalgebra of $H_\Lambda$ by \cite[Proposition 5.23 p.~43]{CKLW18}. By Proposition \ref{prop:lattice-nets}, $H_\Lambda$ is strongly local, and by part 2 of Proposition \ref{prop:CKLW}, the unitary subalgebras $W$ and $V'$ are strongly local with nets embedded as covariant subnets. The remaining assertions follow from Proposition \ref{prop:lattice-nets} and Lemma \ref{lem:trivial-coset}. For $S$, the hypothesis holds by \cite[Theorem 10.10]{Arh26}, and the scalar product $\{\cdot|\cdot\}_\Lambda$ is $H_\chi$ by part 3 of Lemma \ref{lem:unitary-extension}.
\end{proof}

\begin{remark}\normalfont
By part 1 of Proposition \ref{prop:CKLW}, the net $\cal A_W$ does not depend, up to isomorphism, on the normalized invariant scalar product on $W$.
\end{remark}

%%%%%%%%%%%%%%%%%%%%%%%%%%%%%%%%%%%%%%%%%%%%%%%%%%%%%%%%%%%%%%%%%%%%%%%%%%%%%
\subsection{Submodules and representations}
\label{subsec:submodules}

Let $V$ be a strongly local simple unitary vertex operator algebra and $W$ a unitary subalgebra containing the conformal vector of $V$. Let $e_W$ be the orthogonal projection onto the closure $\cal H_W$ of $W$. By \cite[Theorem 7.1 p.~55]{CKLW18}, the conformal net of $W$ is $I \mapsto \cal B_W(I)e_W$, where $\cal B_W(I) = \cal A_V(I) \cap \{e_W\}'$. Let $\cal C_W(I)$ be the von Neumann algebra generated by the closed operators $Y(w,f)$, $w \in W$, $\operatorname{supp} f \subset I$. These operators are affiliated with $\cal A_V(I)$ and commute with $e_W$, and the restriction of $Y(w,f)$ to $\cal H_W$, with domain $D(Y(w,f)) \cap \cal H_W$, is the smeared vertex operator of $W$ \cite[proof of Theorem 7.1 p.~55]{CKLW18}. Hence $\cal C_W(I) \subset \cal B_W(I)$. Since the reduction by $e_W$ of the von Neumann algebra generated by closed operators commuting with $e_W$ is the von Neumann algebra generated by their restrictions to $\cal H_W$, the algebra $\cal C_W(I)e_W$ is generated by the smeared vertex operators of $W$, so that $\cal C_W(I)e_W = \cal B_W(I)e_W$ \cite[Theorem 7.1 p.~55]{CKLW18}. Since $\Omega$ is separating for $\cal A_V(I)$, the map $x \mapsto xe_W$ is injective on $\cal B_W(I)$, and we obtain $\cal C_W(I) = \cal B_W(I)$. We denote this algebra by $\cal A_W(I)$. The restriction map $r_I \co \cal A_W(I) \to \cal A_W(I)|_{\cal H_W}$ is an isomorphism of von Neumann algebras onto the local algebra of the net of $W$. Let $M \subset V$ be a $W$-submodule. Since $\omega \in W$, the subspace $M$ is graded by the eigenvalues of $L_0$. We denote by $\cal H_M$ its closure in $\cal H_V$.

\begin{lemma}
\label{lem:submodules}
\begin{enumerate}
\item We have $V = M \oplus (M^\perp \cap V)$ as $W$-modules. For any interval $I$, the subspace $\cal H_M$ is invariant under $\cal A_W(I)$, and, for any $w \in W$ and any $f \in C^\infty(S^1)$, the restriction $Y_M(w,f) \ov{\mathrm{def}}{=} Y(w,f)|_{\cal H_M}$, with domain $D(Y(w,f)) \cap \cal H_M$, is a closed operator on $\cal H_M$, and $M$ is a core for it. The formula $\pi_M(x) \ov{\mathrm{def}}{=} r_I^{-1}(x)|_{\cal H_M}$, for $x \in \cal A_W(I)|_{\cal H_W}$, defines a representation $\pi_M$ of the conformal net of $W$ on $\cal H_M$.
\item If $M$ is a simple $W$-module, then $\pi_M$ is irreducible.
\item Let $V_1$ and $V_2$ be strongly local simple unitary vertex operator algebras which both contain $W$ as a unitary subalgebra with the same conformal vector, the scalar products of $V_1$ and $V_2$ restricting to the same scalar product on $W$. Let $M_1 \subset V_1$ and $M_2 \subset V_2$ be simple $W$-submodules. Then the representations $\pi_{M_1}$ and $\pi_{M_2}$ are unitarily equivalent if and only if the $W$-modules $M_1$ and $M_2$ are isomorphic.
\end{enumerate}
\end{lemma}

\begin{proof}
1. Let $P$ be the orthogonal projection onto $\cal H_M$. Since the homogeneous spaces are finite-dimensional and orthogonal, $P$ maps $V$ onto $M$ and $(1-P)V = M^\perp \cap V$ is the sum of the orthogonal complements of the spaces $M(n) \ov{\mathrm{def}}{=} M \cap V_n$ in $V_n$. Let $w \in W$ be homogeneous. By \cite[(95) p.~43]{CKLW18}, the adjoint of a mode $w_n$ on $V$ is a finite linear combination of modes of the vectors $L_1^j\theta w$, which lie in $W$. Hence the modes of $w$ preserve $M$ and $M^\perp \cap V$. Consequently, the operator $Y_0(w,f)$ with domain $V$ satisfies $PY_0(w,f)b = Y_0(w,f)Pb$ for any $b \in V$. Since $Y(w,f)$ is the closure of $Y_0(w,f)$ \cite[p.~49]{CKLW18} and $PV \subset V$, we obtain $PY(w,f) \subset Y(w,f)P$. The same holds for the adjoint, so that $P$ commutes with the von Neumann algebra generated by the operators $Y(w,f)$, $\operatorname{supp} f \subset I$, that is, with $\cal A_W(I)$. If $\xi \in \cal H_M$ belongs to the domain of $Y(w,f)$, and if $b_k \in V$ satisfies $b_k \to \xi$ and $Y(w,f)b_k \to Y(w,f)\xi$, then $Pb_k \in M$, $Pb_k \to \xi$ and $Y(w,f)Pb_k = PY(w,f)b_k \to Y(w,f)\xi$. Hence $M$ is a core for $Y_M(w,f)$. The compatibility of the maps $\pi_M$ with the inclusions of intervals holds by construction.

 Since the normal isomorphism $r_I$ maps the operator $Y(w,f)$ affiliated with $\cal A_W(I)$ to its restriction to $\cal H_W$, the representation $\pi_M$ maps the smeared vertex operator $Y(w,f)|_{\cal H_W}$ of $W$ to $Y_M(w,f)$.

2. We first show that $A \subset L_0|_{\cal H_M}$ for the operator $A$ defined below. Choose real functions $f_1, f_2 \in C^\infty(S^1)$ with support in proper intervals and $f_1+f_2 = 1$. The vector $\omega$ is quasi-primary and fixed by $\theta$, so that $Y(\omega,f_j)$ is symmetric by \cite[(113) p.~50]{CKLW18}. Hence the operator $A = Y_M(\omega,f_1)+Y_M(\omega,f_2)$ with domain $D(Y_M(\omega,f_1)) \cap D(Y_M(\omega,f_2))$ is symmetric. On $M$, it coincides with $Y_0(\omega,1) = L_0$. The restriction of $L_0$ to $\cal H_M$ is self-adjoint and $M$ is a core for it. Thus the closure of $A$ is a symmetric extension of a self-adjoint operator, so it is equal to it, and $A \subset L_0|_{\cal H_M}$.

Let $T \in \pi_M(\cal A_W)'$ be a projection. By part 1, $TY_M(w,f) \subset Y_M(w,f)T$ whenever $\operatorname{supp} f$ is contained in a proper interval. For $\xi \in M$, we obtain $T\xi \in D(A)$ and $AT\xi = TA\xi = TL_0\xi$. Hence $L_0T\xi = TL_0\xi$, and $T$ maps each eigenspace $M(n)$ of $L_0|_{\cal H_M}$ into itself. Now let $w \in W$ be homogeneous and $\xi \in M(k)$. The vector $Y_0(w,f)\xi = \sum_n \hat f_n w_n\xi$ has component $\hat f_nw_n\xi$ in $M(k-n)$. Comparing the components of $TY_M(w,f)\xi = Y_M(w,f)T\xi$, we obtain $\hat f_nTw_n\xi = \hat f_nw_nT\xi$. For any $n$ there exists $f$ supported in a proper interval with $\hat f_n \neq 0$, for instance $f(z) = z^ng(z)$ with $g \geq 0$ nonzero and supported in a proper interval. Hence $T|_M$ is an endomorphism of the simple $W$-module $M$, so it is a scalar, and $T \in \{0,1\}$.

3. The restrictions to $W$ of the normalized invariant bilinear forms of $V_1$ and $V_2$ are normalized invariant bilinear forms on $W$, hence they coincide \cite{Li94}. Since the scalar products also coincide on $W$, the PCT operators of $V_1$ and $V_2$ restrict to the same antilinear map $\theta$ of $W$, and the adjoint formula \cite[(95) p.~43]{CKLW18} is the same in $V_1$ and $V_2$. If $\phi \co M_1 \to M_2$ is an isomorphism of $W$-modules, then $(u,v) \mapsto (\phi u|\phi v)$ is an invariant scalar product on $M_1$. The invariant sesquilinear forms on $M_1$ correspond to the module maps from the simple module $\ovl{M_1}$, with conjugate scalars and the action of $w$ given by that of $\theta w$, to the simple contragredient module $M_1'$. Hence they form a space of dimension at most one. Hence $\phi$ is a positive multiple of a unitary map $u \co M_1 \to M_2$. It preserves the grading and extends to a unitary $u \co \cal H_{M_1} \to \cal H_{M_2}$. Since $M_k$ is a core for $Y_{M_k}(w,f)$, we obtain $uY_{M_1}(w,f)u^* = Y_{M_2}(w,f)$ for any $w \in W$ and any $f$ with support in $I$. By part 1, the normal representations $x \mapsto u\pi_{M_1}(x)u^*$ and $\pi_{M_2}$ of $\cal A_W(I)|_{\cal H_W}$ map the closed operator $Y(w,f)|_{\cal H_W}$ to the same closed operator. Hence they coincide on the partial isometries and on the spectral projections of the moduli of the polar decompositions of these operators, which generate $\cal A_W(I)|_{\cal H_W}$. Thus $u$ intertwines $\pi_{M_1}$ and $\pi_{M_2}$.

Conversely, let $T \co \cal H_{M_1} \to \cal H_{M_2}$ be a unitary intertwiner. The operator $\left(\begin{smallmatrix} 0 & 0 \\ T & 0\end{smallmatrix}\right)$ commutes with $(\pi_{M_1} \oplus \pi_{M_2})(\cal A_W(I))$, with which the closed operator $Y_{M_1}(w,f) \oplus Y_{M_2}(w,f)$ is affiliated by part 1. Hence $TY_{M_1}(w,f) \subset Y_{M_2}(w,f)T$. The argument of part 2, applied to the operators $A_k = Y_{M_k}(\omega,f_1)+Y_{M_k}(\omega,f_2) \subset L_0|_{\cal H_{M_k}}$, shows that $T$ maps $M_1(n)$ into $M_2(n)$ and that $T|_{M_1} \co M_1 \to M_2$ is a nonzero morphism of $W$-modules. Since $M_1$ and $M_2$ are simple, it is an isomorphism.
\end{proof}

\begin{remark}\normalfont
In the framework of bounded localized vertex operators, the statements of Lemma \ref{lem:submodules} are counterparts of \cite[Proposition 5.4 p.~50, Proposition 5.5 p.~51 and Theorem 5.6 p.~52]{Ten19a}.
\end{remark}
%%%%%%%%%%%%%%%%%%%%%%%%%%%%%%%%%%%%%%%%%%%%%%%%%%%%%%%%%%%%%%%%%%%%%%%%%%%%%
\section{Unitarity of fusion products and complete rationality}
\label{sec:rationality}

In Section \ref{subsec:krein}, $V$ is a general simple regular unitary vertex operator algebra. From Section \ref{subsec:rationality-U} on, $U \subset V'$ is a Haagerup algebra which is a vertex subalgebra of $V'$ containing the conformal vector of $V'$. The main example is $U = S$ under $\mathrm{(H_J)}$.

%%%%%%%%%%%%%%%%%%%%%%%%%%%%%%%%%%%%%%%%%%%%%%%%%%%%%%%%%%%%%%%%%%%%%%%%%%%%%
\subsection{Positivity of transport forms and unitarity of fusion products}
\label{subsec:krein}

Let $V$ be a simple regular unitary vertex operator algebra and let $M$, $N$ be unitary $V$-modules. Throughout this subsection, modules are ordinary grading-restricted modules. By regularity they are finite direct sums of irreducible ordinary modules. Their conformal weights are rational \cite[Theorem 11.3 p.~48]{DLM00}, in particular real, and their $L_0$-eigenspaces are finite-dimensional. All forms and module homomorphisms are initially defined on the algebraic finite-energy subspaces. Tener defined the transport form $\langle \cdot,\cdot\rangle_{z,w}$ on $M \ot N$, for $(z,w)$ in the domain $\cal R$ of \cite[(3.5) p.~24]{Ten24}, by means of the annihilation operator $\cal Y^-_M$ of $M$ and of the vertex operator of $N$ \cite[Definition 3.1 p.~25]{Ten24}. A transport module for $M$ and $N$ is a unitary $V$-module $K$ with a dominant intertwining operator $\cal Y$ of type $\binom{K}{M N}$ belonging to the space $\cal I_{\mathrm{Hilb}}$ of \cite[Definition 2.8 p.~12]{Ten24}, such that $\langle a_1 \ot b_1,a_2 \ot b_2\rangle_{z,w} = \langle \cal Y(a_1,z)b_1,\cal Y(a_2,w)b_2\rangle$ on $\cal R$ \cite[Definition 3.3 p.~26]{Ten24}. Tener proved that the fusion product $M \boxtimes N$ is a transport module as soon as the transport form is positive semidefinite at one point $(z,z)$ of $\cal R$, under the hypothesis that $M \boxtimes N$ admits a unitary structure \cite[Proposition 3.16 p.~36]{Ten24}, and he noted that the unitarizability of $M \boxtimes N$ is open in general \cite[Remark 3.18 p.~37]{Ten24}. We show that this hypothesis holds in the setting of bounded localized vertex operators when the representation $\pi^N$ exists. We denote by $P = M \boxtimes N$ the fusion product, by $\cal Y_\boxtimes$ its universal intertwining operator, by $P(s)$ the eigenspace of $L_0$ for the eigenvalue $s$ and by $P_s$ the projection onto $P(s)$ defined on the algebraic completion of $P$.

\begin{lemma}
\label{lem:moment-positivity}
Let $\Sigma$ be a countable subset of $[s_0,\infty)$ for some $s_0 \in \R$, let $\epsi > 0$ and let $(c_s)_{s \in \Sigma}$ be complex numbers with $\sum_s |c_s|\e^{\epsi s} < \infty$. If $\sum_s |p(s)|^2c_s \geq 0$ for every polynomial $p \in \Cbb[X]$, then $c_s \geq 0$ for every $s \in \Sigma$.
\end{lemma}

\begin{proof}
Let $\nu = \sum_s |c_s|\delta_s$, a finite positive measure on $[s_0,\infty)$. Let $f \in \L^2(\nu)$ be orthogonal to all polynomials. For $\operatorname{Im}\zeta > -\frac{\epsi}{2}$ and $s \geq s_0$, we have $|\e^{\i\zeta s}| \leq C_\zeta\e^{\frac{\epsi}{2}s}$ for a constant $C_\zeta$ which is locally bounded in $\zeta$, and $\int |f|\e^{\frac{\epsi}{2}s} \d\nu \leq \|f\|_2(\int \e^{\epsi s} \d\nu)^{\frac12} < \infty$. Hence the function $F(\zeta) = \int f(s)\e^{\i\zeta s} \d\nu(s)$ is holomorphic on the half-plane $\operatorname{Im}\zeta > -\frac{\epsi}{2}$, and all its derivatives at $0$ vanish. Thus $F = 0$, the Fourier transform of the finite complex measure $f\nu$ vanishes, and $f = 0$ in $\L^2(\nu)$. So the polynomials are dense in $\L^2(\nu)$. Fix $t \in \Sigma$ and choose polynomials $p_n$ with $\|p_n-1_{\{t\}}\|_2 \to 0$. Using $\big||\alpha|^2-|\beta|^2\big| \leq |\alpha-\beta|(|\alpha|+|\beta|)$, we obtain
\[
\Big|\sum_s |p_n(s)|^2c_s - c_t\Big| \leq \int \big||p_n|^2-1_{\{t\}}\big| \d\nu \leq \|p_n-1_{\{t\}}\|_2\big(\|p_n\|_2+\|1_{\{t\}}\|_2\big) \xrightarrow[n \to \infty]{} 0.
\]
Since the sums on the left are nonnegative, $c_t \geq 0$.
\end{proof}

\begin{prop}
\label{prop:fusion-unitary}
Let $V$ be a simple regular unitary vertex operator algebra with bounded localized vertex operators, and let $M$ and $N$ be ordinary unitary $V$-modules such that the representation $\pi^N$ of \cite[Definition 2.18 p.~20]{Ten24} exists. Then $M \boxtimes N$ admits an invariant scalar product for which $(M \boxtimes N,\cal Y_\boxtimes)$ is a transport module for $M$ and $N$.
\end{prop}

\begin{proof}
By regularity the ordinary modules $M$ and $N$ are finite direct sums of simple modules. These decompositions can be chosen orthogonal by invariance of their scalar products. The annihilation operator $\cal Y^-_M$ pairs only matching summands of $\ovl M$ and $M$, by its definition as the adjoint of the creation operator. The vertex operator of $N$ preserves each summand, and $\pi^N$ restricts to the representation of each summand by \cite[Proposition 5.5]{Ten19a}. Consequently the transport kernel is the orthogonal direct sum of the kernels for the pairs of summands. The tensor product distributes over these finite direct sums and the sum of the corresponding universal intertwining operators is universal. It therefore suffices to prove the proposition for simple $M$ and $N$.

\emph{A nondegenerate invariant Hermitian form.} The unitary forms identify $\ovl M$ with $M'$ and $\ovl N$ with $N'$. For any ordinary module $X$, the complex conjugate space $\ovl X$ is a module via $Y_{\ovl X}(v,x)\ovl a=\ovl{Y_X(\theta v,x)a}$, with coefficientwise conjugation and $\theta$ the PCT operator of $V$. Conjugation of intertwining operators \cite[(1.25)]{Gui22} acts on coefficients and uses no scalar product on the target. Hence the conjugate of $\cal Y_\boxtimes$ is universal and $\ovl P\cong\ovl M\boxtimes\ovl N$. Rigidity identifies $P'$ with $N'\boxtimes M'$, and braiding gives $M'\boxtimes N'\cong P'$ \cite{Hua08}.

The space $\cal F\cong\operatorname{Hom}_V(\ovl P,P')$ of invariant sesquilinear forms on $P$ therefore contains a nondegenerate form. Here a form $B$ is invariant if
\[
B(Y_P(v,x)a,b)=B(a,Y_P(\e^{xL_1}(-x^{-2})^{L_0}\theta v,x^{-1})b).
\]
The involution $B\mapsto B^*$, with $B^*(a,b)=\ovl{B(b,a)}$, preserves $\cal F$ and its fixed points are a real form of $\cal F$. Since $L_0$ is symmetric for $B$, its real eigenspaces are mutually orthogonal. The radical in the first variable is a submodule. If it is nonzero it contains a simple submodule whose lowest weight is among the finitely many lowest weights of the simple summands of $P$. It follows that nondegeneracy is detected by the product of the determinants of $B$ on these finitely many eigenspaces. This polynomial is nonzero on $\cal F$, so it cannot vanish identically on its real form. Choose an invariant nondegenerate Hermitian form $Q_0$ on $P$ in this way. It is nondegenerate on every $P(s)$.

\emph{An algebraic adjoint.} By associativity and the universal property of $\cal Y_\boxtimes$, there is an intertwining operator $\cal Z$ of type $\binom{N}{\ovl M\ P}$ such that
\[
Y^N(\cal Y^-_M(c,z_1-z_2)a,z_2)b
  =\cal Z(c,z_1)\cal Y_\boxtimes(a,z_2)b
\]
as germs of matrix coefficients in the overlap $|z_1|>|z_2|>|z_1-z_2|>0$. Take positive real base points with $0<z_1-z_2<z_2<z_1$ and real logarithms. Convergence and associativity here are for ordinary modules paired with their graded duals \cite{Hua05}. The product uses $P$ only as an ordinary module.

Write $Q_P=Q_0$ and $Q_N=\langle\cdot,\cdot\rangle_N$, linear in the first variable. The maps $h_X:X\to(\ovl X)'$ given by $h_X(x)(\ovl y)=Q_X(x,y)$, for $X=P,N$, are isomorphisms of modules by invariance and graded nondegeneracy. For the conjugate intertwining operator $\ovl{\cal Z}$ of type $\binom{\ovl N}{M\ \ovl P}$, let $A_r(\ovl{\cal Z})$ denote its contragredient as in \cite[(3.87) and Proposition 3.46]{HLZII}. It uses the factor $\e^{xL_1}\e^{(2r+1)\pi\i L_0}x^{-2L_0}$ on the first input. Define the intertwining operator of type $\binom{P}{M\ N}$ by
\[
\cal Z^\dagger = h_P^{-1}A_{-1}(\ovl{\cal Z})(\mathrm{id}_M,h_N).
\]
This definition is graded and applies to the possibly indefinite form $Q_0$. It is a bijection of the corresponding spaces of intertwining operators. With conjugate logarithms at $w$ and $\ovl w$, put $c=\e^{\ovl w L_1}\e^{\i\pi L_0}(\ovl w)^{-2L_0}\ovl{a_2}$, as in \cite[(3.4) and (3.6)]{Ten24}. The plus sign in the rotation factor is obtained by conjugating the factor $\e^{-\i\pi L_0}$ in $A_{-1}$. For $\xi\in P(s)$ the defining graded pairing gives
\[
\langle \cal Z(c,\ovl w^{-1})\xi,b_2\rangle_N
  =Q_0\bigl(\xi,P_s\cal Z^\dagger(a_2,w)b_2\bigr).
\]
By universality $\cal Z^\dagger=\Gamma\cal Y_\boxtimes$ for a $V$-endomorphism $\Gamma$ of $P$. Its adjoint $\Lambda=\Gamma^\sharp$ is taken grade by grade for $Q_0$. Invariance of $Q_0$ makes $\Lambda$ a $V$-endomorphism, and the last pairing is $Q_0(\Lambda\xi,P_s\cal Y_\boxtimes(a_2,w)b_2)$.

\emph{Nondegeneracy.} Let $u\in\operatorname{Hom}_V(\ovl M\boxtimes P,N)$ correspond to the intertwining operator of which $\cal Y_\boxtimes$ is the algebraic adjoint just defined. By the graded pairing, $\cal Z$ corresponds to $u\circ(\mathrm{id}_{\ovl M}\boxtimes\Lambda)$. At the positive real base points above, the associator of the $P(z)$-tensor products identifies the iterate with the product. Transporting the tensor products to $P(1)$ along the positive real axis gives, with the associativity convention of \cite[(3.21)--(3.22)]{Ten24},
\[
u\circ(\mathrm{id}_{\ovl M}\boxtimes\Lambda)
 = \ell_N\circ(y\boxtimes\mathrm{id}_N)\circ\alpha^{-1}.
\]
Here $\alpha:(\ovl M\boxtimes M)\boxtimes N\to\ovl M\boxtimes P$ is the associator, $\ell_N$ is the left unitor corresponding to $Y^N$, and $y:\ovl M\boxtimes M\to V$ corresponds to $\cal Y^-_M$. Since $M$ is simple, rigidity gives $y=d_M\,\mathrm{ev}_M$, where $\mathrm{ev}_M$ is the evaluation for the dual $\ovl M\cong M'$ and $d_M\in\Cbb$. Creation is nonzero on the vacuum, so annihilation is nonzero and $d_M\ne 0$. Define
\[
X=(\mathrm{id}_M\boxtimes u)\circ\alpha'\circ
   (\mathrm{coev}_M\boxtimes\mathrm{id}_P)\circ\ell_P^{-1},
\]
where $\alpha':(M\boxtimes\ovl M)\boxtimes P\to M\boxtimes(\ovl M\boxtimes P)$. Naturality of the associators and the triangle identity for evaluation and coevaluation give $X\Lambda=d_M\,\mathrm{id}_P$. Thus $\Lambda$ is injective. It preserves the finite-dimensional spaces $P(s)$, so it is invertible. Therefore $Q=Q_0(\Lambda\cdot,\cdot)$ is invariant and nondegenerate and makes distinct eigenspaces of $L_0$ orthogonal.

\emph{The transport identity.} Set $z_1=\ovl w^{-1}$ and $z_2=z$. The product converges absolutely on $\cal R$ because $|z_1|>|z_2|$. Pairing with the algebraic vector $b_2$ uses the graded dual of $N$. Substituting $\xi=P_s\cal Y_\boxtimes(a_1,z)b_1$ into the adjoint identity and summing gives the product expansion
\begin{equation}
\label{eq:krein-transport}
\langle a_1\ot b_1,a_2\ot b_2\rangle_{z,w}
 = \sum_s Q\bigl(P_s\cal Y_\boxtimes(a_1,z)b_1,
                P_s\cal Y_\boxtimes(a_2,w)b_2\bigr).
\end{equation}
To see that the equality holds throughout $\cal R$, start near a point $z=w\in(2^{-1/2},1)$, where the product and the iterate agree as germs with the chosen real branches. For homogeneous inputs, the summand on the right is a coefficient times $z^{s-h_1}\ovl w^{s-h_2}$, where $h_i$ is the sum of the weights of $a_i$ and $b_i$. The common branch convention of \cite[Section 3.1]{Ten24} makes the sum single-valued, since a common change of logarithm has phase $\e^{2\pi\i(h_2-h_1)}=1$ for the simple modules $M$ and $N$. After factoring out the input powers, the series is in the variable $t=z\ovl w$ with the logarithm on the right half-plane. Absolute convergence at $z=w=\sqrt{\rho}$, $1/2<\rho<1$, gives normal convergence on compact subsets with $|t|<\rho$, since the weights are bounded below. Finally $\cal R$ is connected: in coordinates $(t,z)$ it is given by $\operatorname{Re} t>1/2$ and $|t|<|z|<1$. Thus the germ identity continues to the entire domain, with the single-valued iterate of \cite[Definition 3.1]{Ten24} on the left, which converges absolutely on $\cal R$ by \cite[Proposition 3.11 p.~33]{Ten24}.

\emph{Positivity of the transport form.} By \cite[Proposition 3.11 p.~33]{Ten24}, $\cal Y^-_M$ belongs to the space $\cal I_{\mathrm{Hilb}}$ of \cite[Definition 2.8 p.~12]{Ten24}. By \cite[Proposition 4.2 p.~43]{Ten24}, whose hypotheses do not involve $P$, there is $r < 1$ such that the form $\langle\cdot,\cdot\rangle_{z,z}$ on $M \ot N$ is positive semidefinite whenever $(z,z) \in \cal R$ and $|z| > r$. Fix a real number $z$ with $\max(r,2^{-\frac12}) < z < 1$, so that $(z,z) \in \cal R$. For $x \in M \ot N$, put $\xi = \cal Y_\boxtimes(x,z)$, $\xi_s = P_s\xi$ and $c_s(x) = Q(\xi_s,\xi_s)$. By \eqref{eq:krein-transport}, $\sum_s c_s(x) = \langle x,x\rangle_{z,z} \geq 0$.

\emph{Invariance under $L_0$.} For an intertwining operator $\cal Y$, we have $[L_0,\cal Y(a,x)] = x\frac{\d}{\d x}\cal Y(a,x)+\cal Y(L_0a,x)$ and $\frac{\d}{\d x}\cal Y(a,x) = \cal Y(L_{-1}a,x)$. Hence
\[
L_0\cal Y_\boxtimes(a,z)b = \cal Y_\boxtimes(L_0a+zL_{-1}a,z)b+\cal Y_\boxtimes(a,z)L_0b,
\]
where $L_0$ acts componentwise on the algebraic completion $\prod_s P(s)$ of $P$, which contains the vectors $\cal Y_\boxtimes(x,z)$, $x \in M \ot N$, defined by linearity. With the operator $T_z \ov{\mathrm{def}}{=} (L_0+zL_{-1}) \ot \mathrm{id}_N+\mathrm{id}_M \ot L_0$ of $M \ot N$, we obtain $p(L_0)\cal Y_\boxtimes(x,z) = \cal Y_\boxtimes(p(T_z)x,z)$ for any polynomial $p$ and any $x \in M \ot N$. Put $x_p \ov{\mathrm{def}}{=} p(T_z)x$. Since $P_s p(L_0)\xi = p(s)\xi_s$, we get $\sum_s |p(s)|^2c_s(x) = \langle x_p,x_p\rangle_{z,z} \geq 0$.

\emph{Exponential decay.} Let $1 < \lambda < z^{-1}$. If $x = \sum_k a_k \ot b_k$ with homogeneous $a_k$ and $b_k$ of weights $h_k$ and $h'_k$, put $x' = \sum_k \lambda^{h_k+h'_k}a_k \ot b_k$. Then $P_s\cal Y_\boxtimes(x',\lambda z) = \lambda^s\xi_s$, since $\lambda^{L_0}\cal Y_\boxtimes(a,z)b = \cal Y_\boxtimes(\lambda^{L_0}a,\lambda z)\lambda^{L_0}b$. The absolute convergence of \eqref{eq:krein-transport} at the point $(\lambda z,\lambda z) \in \cal R$ gives $\sum_s\lambda^{2s}|c_s(x)| < \infty$. The weights of $P$ are bounded below.

\emph{Conclusion.} By Lemma \ref{lem:moment-positivity}, $c_s(x) \geq 0$ for all $s$ and all $x$. For homogeneous $a \in M$ and $b \in N$ of weights $h$ and $h'$, we have $P_s\cal Y_\boxtimes(a,z)b = z^{-n-1}a_{(n)}b$ with $n = h+h'-s-1$. By the commutator formula, the coefficients $a_{(n)}b$ span a submodule of $P$, which is equal to $P$ by uniqueness in the universal property of $(P,\cal Y_\boxtimes)$ and semisimplicity. That is, $\cal Y_\boxtimes$ is dominant. Hence every $v \in P(s)$ can be written $v = \sum_k (a_k)_{(n_k)}b_k = P_s\cal Y_\boxtimes(x,z)$ with $x = \sum_k z^{n_k+1}a_k \ot b_k$, and $Q(v,v) = c_s(x) \geq 0$. By polarization, $Q$ is a positive semidefinite Hermitian form. It is invariant and nondegenerate, hence positive definite, and $P$ is a unitary $V$-module for $Q$. Finally, \eqref{eq:krein-transport} reads $\langle x,y\rangle_{z,w} = Q(\cal Y_\boxtimes(x,z),\cal Y_\boxtimes(y,w))$ on $\cal R$, and $\cal Y_\boxtimes$ is dominant and lies in $\cal I_{\mathrm{Hilb}}$ by \cite[Proposition 3.11 p.~33]{Ten24}, so that $(P,\cal Y_\boxtimes)$ is a transport module.
\end{proof}

\begin{remark}\normalfont
\label{rem:krein}
The proof of \cite[Proposition 3.16 p.~36]{Ten24} uses a Hilbert completion of an already unitary fusion product. Here nondegeneracy comes from the algebraic adjoint and rigidity, and positivity follows from Lemma \ref{lem:moment-positivity}. The only integration used in Proposition \ref{prop:fusion-unitary} is that of the second input, through \cite[Proposition 4.2]{Ten24}. In particular, if $\pi^M$ exists, the proposition constructs the transport module for $(\ovl M,M)$ in the order required by \cite[Theorem 4.19 p.~55]{Ten24}. This supplies directly the transport module that the proof of \cite[Theorem 4.20 p.~55]{Ten24} obtains from \cite[Corollary 4.5 p.~44]{Ten24}, which assumes, through \cite[Theorem 4.4 p.~44]{Ten24}, that $\ovl M\boxtimes M$ admits a unitary structure. This assumption is not listed in the statement of \cite[Theorem 4.20]{Ten24}, see Remark \ref{rem:delta}.
\end{remark}

%%%%%%%%%%%%%%%%%%%%%%%%%%%%%%%%%%%%%%%%%%%%%%%%%%%%%%%%%%%%%%%%%%%%%%%%%%%%%
\subsection{Unitarity of the modules and complete rationality}
\label{subsec:rationality-U}

\begin{lemma}
\label{lem:V-decomposition}
The $U$-module $V'$ decomposes as $V' = U \oplus M$ with $M \cong M_b$. Consequently, for each isotropic line $\Lambda$, the vertex operator algebra $H_\Lambda$ of Lemma \ref{lem:cocycle} is a holomorphic extension of $U$ containing $V' = U \oplus M_b$, and the facts recalled in Section \ref{subsec:haagerup-algebras} apply to $U \subset V' \subset H_\Lambda$. In particular, the isotropic lines $\Lambda_a$ and $\Lambda_{c_0}$ are defined, and for $w \in \Lambda_x \setminus \{0\}$, $x \in \{a,c_0\}$, the $U$-module $V'_w$ is simple and isomorphic to $M_x$.
\end{lemma}

\begin{proof}
The vertex operator algebra $V'$ is a simple extension of $U$, and $V' \cong (V')'$. By Lemma \ref{lem:intermediate-etale} applied to $W = B = V'$ and by \cite[Theorem 3.4 p.~16]{HKL15}, the algebra $A$ in $\Rep(U)$ attached to $V'$ is a connected \'etale algebra, and the category $\Rep(U)^0_A$ is isomorphic to $\Rep(V')$, which has nine simple objects \cite{Don93}. Hence $\operatorname{FPdim} \Rep(U)^0_A \geq 9$, and by \eqref{eq:dims-local} and \eqref{eq:dim-center}, we have $\operatorname{FPdim}(A)^2 \leq \frac{9(3d+2)^2}{9}$, that is, $\operatorname{FPdim}(A) \leq 3d+2$. As in the proof of Lemma \ref{lem:possible-objects}, $A$ is a direct sum of simple objects with trivial twist containing $0$ once. By \eqref{eq:center-dimensions}, the only such objects of Frobenius--Perron dimension at most $3d+2$ are $0$ and $0 \oplus b$. Since $V' \neq U$, we obtain $A \cong 0 \oplus b$. The other assertions follow, since $H_\Lambda \cong V_{E_8}$ is holomorphic and contains $V'$ by Lemma \ref{lem:cocycle}. Here the part of the proof of \cite[Lemma 4.3]{Arh26} which shows that exactly two nonzero classes restrict to $M_{c_0}$ uses only the decomposition $V' = U \oplus M_b$, the twists and the modular transformation of the one-point functions, so that it applies to the given embedding $U \subset V'$.
\end{proof}

Assume now that $U$ is a unitary subalgebra of $V'$ for a normalized invariant scalar product, and equip each $H_\Lambda$ with the scalar product of part 2 of Lemma \ref{lem:unitary-extension}. Let $N_b$ be the orthogonal complement of $U$ in $V'$, which is a $U$-submodule isomorphic to $M_b$ by Lemma \ref{lem:V-decomposition} and part 1 of Lemma \ref{lem:submodules}. With the notation of Lemma \ref{lem:submodules}, we put
\[
\pi_0 \ov{\mathrm{def}}{=} \pi_U,\quad \pi_b \ov{\mathrm{def}}{=} \pi_{N_b},\quad \pi_x \ov{\mathrm{def}}{=} \pi_{V'_w} \text{ for } w \in \Lambda_x \setminus \{0\},\ x \in \{a,c_0\},
\]
where $V'_w$ is viewed as a $U$-submodule of $H_{\Lambda_x}$. By part 3 of Lemma \ref{lem:submodules}, $\pi_x$ does not depend on the choice of $w$ up to unitary equivalence, and $\pi_0$ is the vacuum representation of $\cal A_U$.

\begin{cor}
\label{cor:modules-unitary}
Let $U \subset V'$ be a Haagerup algebra containing the conformal vector of $V'$ which is a unitary subalgebra for a normalized invariant scalar product on $V'$. Then every simple $U$-module admits a scalar product which is invariant with respect to the PCT operator of $U$ induced by this scalar product. In particular, under $\mathrm{(H_J)}$, every simple $S$-module admits a scalar product which is invariant for the PCT operator of $S$ induced by $H_\chi$, which is the assertion $\mathrm{(H_U)}$ of the introduction.
\end{cor}

\begin{proof}
Fix an isotropic line $\Lambda$ and put $H = H_\Lambda$, with the scalar product of part 2 of Lemma \ref{lem:unitary-extension}. By part 1 of Proposition \ref{prop:lattice-nets}, $H$ has bounded localized vertex operators. As in the proof of Theorem \ref{thm:strong-locality}, $U$ is a unitary subalgebra of $H$ containing its conformal vector, so that $U$ has bounded localized vertex operators, for the same system of generalized annuli \cite[Proposition 4.16 p.~43]{Ten19a}. The $U$-module $N_b$ is a unitary $U$-submodule of $H$ isomorphic to $M_b$, and the representation $\pi^{N_b}$ exists by \cite[Theorem 5.6 p.~52]{Ten19a}, applied to the $H$-module $H$, for which $\pi^H$ is the vacuum representation. The strongly rational algebra $U$ is regular \cite[Theorem 4.5]{ABD04}. By Proposition \ref{prop:fusion-unitary}, $N_b \boxtimes N_b$ admits an invariant scalar product. In the normalized $S$-matrix of \cite[(2.5a)--(2.5b) p.~20]{EvG11}, put $x = \frac{13-3\sqrt{13}}{26}$ and $y = \frac3{\sqrt{13}}$. Then $x(1-x) = \frac1{13}$ and $1-2x = y$. Outside the columns $0$ and $b$, the squares of the entries in the rows $0$ and $b$ agree. Since all simple objects are self-dual, $S^2 = 1$, and the Verlinde formula gives
\[
N_{bb}^k = \delta_{0k}+\frac{1-2x}{3}\left(\frac{S_{k0}}{x}-\frac{S_{kb}}{1-x}\right).
\]
For $k = 0$ the added term vanishes. For $k \in \{b,a,c_0,c_1,c_2\}$ it is $\frac{(1-2x)^2}{9x(1-x)} = 1$, and for $k = d_l$ it is $\frac{(1-2x)y}{9x(1-x)} = 1$. Thus $b \ot b$ is the direct sum of the twelve simple objects, each with multiplicity one. Hence every simple $U$-module is isomorphic to a submodule of $N_b \boxtimes N_b$, and the restriction of the scalar product to this submodule is invariant. The last assertion follows, since under $\mathrm{(H_J)}$ the algebra $S$ is a Haagerup algebra which contains the conformal vector of $V'$ and is a unitary subalgebra of $V'$ \cite[Theorem 10.10 and Theorem 10.24]{Arh26}.
\end{proof}

\begin{thm}
\label{thm:complete-rationality}
Let $U \subset V'$ be a Haagerup algebra containing the conformal vector of $V'$ which is a unitary subalgebra for a normalized invariant scalar product on $V'$.
\begin{enumerate}
\item The conformal net $\cal A_U$ is completely rational. For $x \in \{a,c_0\}$, the index $[\cal A_{H_{\Lambda_x}} : \cal A_U]$ is finite.
\item The representations $\pi_0$, $\pi_b$, $\pi_a$, $\pi_{c_0}$ are irreducible and pairwise inequivalent. The vacuum representation of $\cal A_{V'}$ restricts to $\cal A_U$ as $\pi_0 \oplus \pi_b$, and the vacuum representation of $\cal A_{H_{\Lambda_x}}$ restricts to $\cal A_U$ as $\pi_0 \oplus \pi_b \oplus 2\pi_x$.
\item Put $\delta \ov{\mathrm{def}}{=} d(\pi_b)$. Then $\delta \geq 1$, $[\cal A_{V'} : \cal A_U] = 1+\delta$, $d(\pi_a) = d(\pi_{c_0}) = 1+\delta$, $[\cal A_{H_{\Lambda_x}} : \cal A_U] = 3(1+\delta)$ and $\mu_{\cal A_U} = 9(1+\delta)^2$.
\item For $x \in \{a,c_0\}$, the dual canonical endomorphism $\Theta_x$ of $\cal A_U \subset \cal A_{H_{\Lambda_x}}$ is a Lagrangian algebra in $\Rep(\cal A_U)$, and the fusion category $\cal F_x \ov{\mathrm{def}}{=} \Rep(\cal A_U)_{\Theta_x}$ is a unitary fusion category of rank six and Frobenius--Perron dimension $3(1+\delta)$ such that $\Rep(\cal A_U)$ is braided equivalent to $Z(\cal F_x)$. The category $\cal F_x$ contains a pointed fusion subcategory of dimension three, corresponding to the intermediate net $\cal A_{V'}$.
\end{enumerate}
\end{thm}

\begin{proof}
1. Fix $x$ and put $H=H_{\Lambda_x}$, with the scalar product of Lemma \ref{lem:unitary-extension}. By Proposition \ref{prop:lattice-nets}, $H$ is strongly local, has bounded localized vertex operators, and its nets in the sense of \cite{Ten19a} and \cite{CKLW18} are isomorphic. By the proof of Theorem \ref{thm:strong-locality}, $U$ is a unitary conformal subalgebra of $H$. The algebra $U$ is regular by \cite[Theorem 4.5]{ABD04} and inherits bounded localized vertex operators by \cite[Proposition 4.16]{Ten19a}.

Let $M$ be a simple $U$-summand of $H$. The representations $\pi^M$ and $\pi^{\ovl M}$ of $\cal A_U^{\mathrm T}$ in the sense of \cite[Definition 2.18]{Ten24} exist by \cite[Theorem 5.6 and Lemma 5.7]{Ten19a}, applied to $M$ and to its image under the PCT operator of $H$, which realizes $\ovl M$. With $p_M$ the orthogonal projection onto $M$, skew symmetry gives $p_MY^H(m,z)v=\cal Y_M^+(m,z)v$ for $m\in M$, $v\in U$. The common conformal vector makes the power shift in \cite[Theorem 5.15]{Ten19a} zero, so $\cal Y_M^+\in\cal I_{\mathrm{loc}}$. Proposition \ref{prop:fusion-unitary}, with inputs $\ovl M,M$, makes $\ovl M\boxtimes M$ a transport module, and rigidity and semisimplicity give a copy of $U$ in it. Now \cite[Theorem 4.19]{Ten24} gives finite index for each such sector $\pi^M$.

The restriction $\pi_H$ of the vacuum of $\cal A_H^{\mathrm T}$ to $\cal A_U^{\mathrm T}$ is a finite sum of these sectors and has finite index. For a complementary interval $I'$ the inclusion
\[
\pi_H(\cal A_U^{\mathrm T}(I))
 \subset \cal A_H^{\mathrm T}(I)
 \subset \pi_H(\cal A_U^{\mathrm T}(I'))'
\]
then gives finite index for the subnet by multiplicativity of index, as in the last step of \cite[Corollary 4.24]{Ten24}. The ambient net is completely rational by part 1 of Proposition \ref{prop:lattice-nets}, so part 1 of Proposition \ref{prop:nets} gives complete rationality of $\cal A_U^{\mathrm T}$. By \cite[Lemma 7.6 p.~72]{Ten19a}, this net is isomorphic to the net $\cal A_U$ of \cite{CKLW18}, which is therefore completely rational.

Now consider the inclusion $\cal A_U \subset \cal A_H$ of nets of \cite{CKLW18}. By Theorem \ref{thm:strong-locality}, its coset is trivial and $\cal A_U(I)' \cap \cal A_H(I) = \Cbb$. In particular, $\cal A_H$ is an irreducible local extension of $\cal A_U$ in the sense of \cite[p.~9]{KL04}. Since $\cal A_U$ is completely rational, the index $[\cal A_H : \cal A_U]$ is finite by part 2 of Proposition \ref{prop:nets}. The same argument applies to the inclusion $\cal A_U \subset \cal A_{V'}$.

2. As a $U$-module, $H = U \oplus N_b \oplus V'_w \oplus V'_{-w}$ for $w \in \Lambda_x \setminus \{0\}$, and the summands are pairwise orthogonal. Indeed, the lattice sectors are orthogonal for $H_1$, the automorphism $\tilde h = \exp(y(0))$, $y \in (V')_1$, of Lemma \ref{lem:unitary-extension} preserves each $V'_w$, and $N_b$ is the orthogonal complement of $U$ in $V'$. Hence $\cal H_H = \cal H_U \oplus \cal H_{N_b} \oplus \cal H_{V'_w} \oplus \cal H_{V'_{-w}}$, and by part 1 of Lemma \ref{lem:submodules}, the vacuum representation of $\cal A_H$ restricts to $\pi_0 \oplus \pi_b \oplus \pi_{V'_w} \oplus \pi_{V'_{-w}}$. By Lemma \ref{lem:V-decomposition}, the summands are simple $U$-modules, isomorphic to $U$, $M_b$, $M_x$, $M_x$. These are pairwise non-isomorphic except for the last two. By parts 2 and 3 of Lemma \ref{lem:submodules}, we obtain the decomposition $\pi_0 \oplus \pi_b \oplus 2\pi_x$ into irreducible representations, with $\pi_0$, $\pi_b$, $\pi_x$ pairwise inequivalent. Applying part 3 of Lemma \ref{lem:submodules} to $V_1 = H_{\Lambda_a}$, $V_2 = H_{\Lambda_{c_0}}$ and to the non-isomorphic $U$-modules $M_a$ and $M_{c_0}$, we see that $\pi_a$ and $\pi_{c_0}$ are inequivalent. The same argument applied to $\cal H_{V'} = \cal H_U \oplus \cal H_{N_b}$ gives the restriction of the vacuum representation of $\cal A_{V'}$.

3. By part 1 and part 3 of Proposition \ref{prop:nets}, the dual canonical endomorphism of $\cal A_U \subset \cal A_{V'}$ is $\theta = \pi_0 \oplus \pi_b$, so that $[\cal A_{V'}:\cal A_U] = d(\theta) = 1+\delta$. Here $\delta \geq 1$, as for any statistical dimension. By part 1 of Proposition \ref{prop:nets} applied to the two finite index inclusions $\cal A_U \subset \cal A_{V'}$ and $\cal A_U \subset \cal A_H$, and by Proposition \ref{prop:lattice-nets}, we have
\[
\mu_{\cal A_U} = [\cal A_{V'}:\cal A_U]^2\mu_{\cal A_{V'}} = 9(1+\delta)^2 \quad \text{and} \quad \mu_{\cal A_U} = [\cal A_H:\cal A_U]^2\mu_{\cal A_H} = [\cal A_H:\cal A_U]^2.
\]
Hence $1+\delta+2d(\pi_x) = d(\Theta_x) = [\cal A_H : \cal A_U] = 3(1+\delta)$, so that $d(\pi_x) = 1+\delta$.

4. By part 3 of Proposition \ref{prop:nets} and part 1, $\Theta_x$ is a commutative $\Q$-system with
\[
\dim\operatorname{Hom}(\mathrm{id},\Theta_x) = \dim \cal A_U(I)' \cap \cal A_H(I) = 1,
\]
hence a connected \'etale algebra in the unitary modular tensor category $\Rep(\cal A_U)$. By \eqref{eq:mu-index}, $\operatorname{FPdim}\Rep(\cal A_U) = \mu_{\cal A_U} = d(\Theta_x)^2$, so that $\Theta_x$ is Lagrangian. By part 1 of Proposition \ref{prop:DMNO}, $\Rep(\cal A_U)$ is braided equivalent to $Z(\cal F_x)$, and \eqref{eq:dims-local} gives $\operatorname{FPdim}(\cal F_x) = \frac{9(1+\delta)^2}{3(1+\delta)} = 3(1+\delta)$. The category $\cal F_x$ of $\Theta_x$-modules is equivalent to the category of $\alpha$-induced sectors of the inclusion $\cal A_U(I) \subset \cal A_H(I)$, which is a unitary fusion category \cite[Section 4.6 p.~54]{BKLR15}, see also \cite[Proposition 4.5 p.~21]{Bis17a}. By Lemma \ref{lem:rank} and part 2, its rank is $\dim\operatorname{Hom}(\Theta_x,\Theta_x) = 1+1+4 = 6$. By \cite[Proposition 4.11 p.~49]{BKLR15}, the intermediate subfactor $\cal A_U(I) \subset \cal A_{V'}(I) \subset \cal A_H(I)$ corresponds to a projection in $\operatorname{Hom}(\Theta_x,\Theta_x)$ whose range, with the restricted multiplication, is the $\Q$-system $\theta$ of $\cal A_U \subset \cal A_{V'}$. Hence $\theta$ is a connected \'etale subalgebra of $\Theta_x$. By Lemma \ref{lem:rank}, under the equivalence $\Rep(\cal A_U) \cong Z(\cal F_x)$, the algebra $\Theta_x$ corresponds to the algebra $I(1)$, and by part 3 of Proposition \ref{prop:DMNO}, $\theta$ corresponds to a fusion subcategory $\cal P_x$ of $\cal F_x$ with $\operatorname{FPdim}(\theta)\operatorname{FPdim}(\cal P_x) = \operatorname{FPdim}(\cal F_x)$, that is, $\operatorname{FPdim}(\cal P_x) = 3$. A fusion category of prime Frobenius--Perron dimension is pointed \cite[Corollary 8.30 p.~37]{ENO05}.
\end{proof}

\begin{prop}
\label{prop:hypergroup-U}
Under the hypotheses of Theorem \ref{thm:complete-rationality}, there is a proper action of the hypergroup $K = \{k_0,k_1\}$ with
\[
k_1k_1 = \frac{1}{\delta}k_0+\Big(1-\frac{1}{\delta}\Big)k_1
\]
on the net $\cal A_{V'}$ such that $\cal A_U = \cal A_{V'}^K$. The weight of $K$ is $1+\delta = [\cal A_{V'} : \cal A_U]$.
\end{prop}

\begin{proof}
By part 1 of Theorem \ref{thm:complete-rationality}, the net $\cal A_U$ is completely rational and $[\cal A_{V'} : \cal A_U] < \infty$. By part 3 of Proposition \ref{prop:bischoff}, there is a canonical hypergroup $K$ acting properly on $\cal A_{V'}$ with $\cal A_U = \cal A_{V'}^K$, and its weight is $[\cal A_{V'} : \cal A_U] = 1+\delta$. By part 2 of Theorem \ref{thm:complete-rationality}, the dual canonical endomorphism of $\cal A_U(I) \subset \cal A_{V'}(I)$ is $\theta = \pi_0 \oplus \pi_b$ with $\pi_b$ irreducible and inequivalent to $\pi_0$. Hence $\langle \gamma,\gamma\rangle = \langle\theta,\theta\rangle = 2$, where $\gamma$ is the canonical endomorphism, and by \cite[Section 4.6 p.~25]{Bis17b}, the canonical hypergroup is $K = \{k_0,k_1\}$ with the stated multiplication, where $\delta = [\cal A_{V'}:\cal A_U]-1$. The elements $k_0$, $k_1$ are denoted $c_0$, $c_1$ in \cite{Bis17b}.
\end{proof}

\begin{prop}
\label{prop:delta-big}
Under the hypotheses of Theorem \ref{thm:complete-rationality}, we have $\delta > 1$. Consequently, the fusion category $\cal F_x$ has exactly three invertible simple objects and three non-invertible simple objects.
\end{prop}

\begin{proof}
Assume that $\delta = 1$. By Proposition \ref{prop:hypergroup-U}, there is a proper action of a hypergroup $K = \{k_0,k_1\}$ on $\cal A_{V'}$ with $\cal A_U = \cal A_{V'}^K$ and $k_1k_1 = \delta^{-1}k_0+(1-\delta^{-1})k_1 = k_0$. Thus the vacuum preserving unital completely positive maps $\sigma_I = \phi^I_{k_1}$ of $\cal A_{V'}(I)$ satisfy $\sigma_I^2 = \mathrm{id}$. By the Kadison--Schwarz inequality and the positivity of $\sigma_I$, we have $x^*x = \sigma_I^2(x^*x) \geq \sigma_I(\sigma_I(x)^*\sigma_I(x)) \geq x^*x$. Since $\sigma_I$ is injective and positive, we obtain $\sigma_I(x^*x) = \sigma_I(x)^*\sigma_I(x)$ for any $x$, so that $\sigma_I$ is multiplicative by the multiplicative domain argument, hence an automorphism. It is nontrivial since $\cal A_U \neq \cal A_{V'}$. Since $\sigma_I$ preserves the vacuum state and $\Omega$ is cyclic for $\cal A_{V'}(I)$, the formula $u_I x\Omega = \sigma_I(x)\Omega$ defines a unitary $u_I$ with $u_I x u_I^* = \sigma_I(x)$ and $u_I\Omega = \Omega$. The maps $\sigma_I$ are compatible \cite[Proposition 4.21 p.~23]{Bis17b}, so that $u_I = u_J$ if $I \subset J$. Since any two intervals contain subintervals which are contained in a common interval, $g \ov{\mathrm{def}}{=} u_I$ does not depend on $I$, and $g$ is an automorphism of the net $\cal A_{V'}$ in the sense of \cite[(41) p.~21]{CKLW18}. By \cite[Theorem 6.9 p.~53]{CKLW18}, $g$ is a unitary automorphism of $V'$ of order two. By \cite[Proposition 7.6 and Theorem 7.5 p.~58]{CKLW18}, we obtain $U = (V')^g$, so that $U_1$ is the fixed point subalgebra of the automorphism $g$ of order at most two of $(V')_1 \cong \mathfrak e_6 \oplus \mathfrak{sl}_3$. The automorphism $g$ preserves the two simple ideals $\mathfrak e_6$ and $\mathfrak{sl}_3$ of $(V')_1$, which are not isomorphic, and it is nontrivial on $(V')_1$ since $V'$ is generated by $(V')_1$. The fixed point subalgebra of an involution of $\mathfrak e_6$ is isomorphic to $\mathfrak{sp}_8$, $\mathfrak{sl}_6 \oplus \mathfrak{sl}_2$, $\mathfrak{so}_{10} \oplus \Cbb$ or $\mathfrak f_4$, of dimension $36$, $38$, $46$ or $52$, and that of an involution of $\mathfrak{sl}_3$ is isomorphic to $\mathfrak{so}_3$ or $\mathfrak{gl}_2$, of dimension $3$ or $4$ \cite[Chapter 8]{Kac90}, \cite[Chapter X]{Hel78}. Taking into account the case where $g$ is trivial on one of the two ideals, of dimensions $78$ and $8$, we obtain $\dim U_1 \geq \min(36+3,36+8,78+3) = 39$, which contradicts $\dim U_1 = 6$ \cite[Theorem 7.1]{Arh26}. A similar argument, with the rank of the fixed point subalgebra, appears in \cite[pp.~7--8]{Hua26}.

The invertible objects of $\cal F_x$ form a group containing the group of order three of $\cal P_x$. Since $\cal F_x$ has rank six, the order of this group is three or six by Lagrange's theorem. If this order were six, then $\cal F_x$ would be pointed of dimension $6 = 3(1+\delta)$, so that $\delta = 1$.
\end{proof}

\begin{cor}
\label{cor:S-rational}
Assume $\mathrm{(H_J)}$. Then Theorem \ref{thm:complete-rationality} holds for $U = S$, with the unitary structure of \cite[Theorem 10.10]{Arh26}. In particular, $\cal A_S$ is completely rational, and $\Rep(\cal A_S)$ is braided equivalent to the Drinfeld center of a unitary fusion category of rank six with a pointed fusion subcategory of dimension three.
\end{cor}

\begin{proof}
By \cite[Theorem 10.24]{Arh26}, $\mathrm{(H_J)}$ implies that $S$ is a Haagerup algebra. It contains the conformal vector of $V'$ and is a unitary subalgebra for $H_\chi$ by \cite[Theorem 10.10]{Arh26}.
\end{proof}

\begin{remark}\normalfont
\label{rem:delta}
The proof of Theorem \ref{thm:complete-rationality} uses the transport modules $\ovl M\boxtimes_t M$ of Proposition \ref{prop:fusion-unitary} for the simple summands $M$ of $H_\Lambda$, together with \cite[Theorem 4.19]{Ten24}. Corollary \ref{cor:modules-unitary} uses the proposition for a different purpose. Under $\mathrm{(H_J)}$, the fusion rules of $Z(\cal H_3)$ show that $M_b\boxtimes M_b$ contains every simple module, so the invariant scalar product on this product also unitarizes the six modules $M_{d_l}$. This use of the fusion rules is separate from the index argument.

For the algebra of \cite{Hua26}, the complete rationality of the net and the unitarity of all simple modules are announced in \cite[Theorem 1 p.~1 and Section 6 pp.~8--9]{Hua26}. In the outline of \cite[Section 6 p.~8]{Hua26}, the finiteness of the index of the subnet generated by the bounded insertions is obtained from \cite{Ten24} and the regularity of the algebra, without reference to a specific statement, and the unitarity of the simple modules is derived afterwards from the twelve sectors of the resulting completely rational net, through \cite[Theorem 9.5 and Remark 9.8]{HeT25}. The finite index statements \cite[Theorem B and Corollary C p.~6, Theorem 4.20 p.~55 and Corollaries 4.21 and 4.24 pp.~55--56]{Ten24} do not list a unitarity hypothesis, but their proofs use \cite[Corollary 4.5 p.~44]{Ten24}, which requires a unitary structure on $\ovl M \boxtimes M$, see also \cite[Remark 3.18 p.~37]{Ten24}. For $M = M_b$, this module contains every simple module, including the six modules $M_{d_l}$. Proposition \ref{prop:fusion-unitary} supplies the transport module for algebras with bounded localized vertex operators as soon as $\pi^M$ exists. Hence this step of the outline can be completed as in the proofs of Corollary \ref{cor:modules-unitary} and Theorem \ref{thm:complete-rationality}, provided that the algebra of \cite{Hua26} is realized as a unitary subalgebra of $V'$ containing its conformal vector, as in \cite[Remark 7.14]{ArhIII}. The determination of the fusion rules and of the associator in \cite[Section 5 pp.~6--8]{Hua26} only uses the pseudo-unitarity of the category of modules over the extension algebra, and it does not depend on this issue. The application of the classification of $C^*$-generalized Haagerup categories \cite[Theorem 5.9 p.~36]{Izu18} in \cite[Section 6 p.~9]{Hua26} uses the unitarity of this category, which is deduced there from the complete unitarity of the algebra. In Section \ref{sec:recognition}, we apply this classification to the category $\cal F_x$ of the net, which is unitary by construction, see Remark \ref{rem:huang-recognition}.

The value $\delta = d^2$ is obtained in Theorem \ref{thm:Fx-haagerup}.
\end{remark}
%%%%%%%%%%%%%%%%%%%%%%%%%%%%%%%%%%%%%%%%%%%%%%%%%%%%%%%%%%%%%%%%%%%%%%%%%%%%%
\section{A recognition theorem for the Haagerup categories}
\label{sec:recognition}

In this section we prove Theorem \ref{thm:recognition}, which is Theorem \ref{thm:intro-D} of the introduction, and we deduce from it that each category $\cal F_x$ of Theorem \ref{thm:complete-rationality} is equivalent to $\cal H_2$ or to $\cal H_3$ (Theorem \ref{thm:Fx-haagerup}). The strategy for the fusion rules, which combines the induction to the Drinfeld center, Frobenius--Schur indicators and lower bounds for sums of roots of unity, is the one outlined by Huang \cite[Sections 5 and 6 pp.~7--9]{Hua26}, and we compare our proofs with this outline in Remark \ref{rem:huang-recognition}.

Throughout this section, $\cal F$ is a unitary fusion category. We equip it with the spherical structure for which the dimension $d(X)$ of every simple object $X$ is its Frobenius--Perron dimension \cite[Proposition 8.23 p.~34]{ENO05}, and we put $\dim \cal F \ov{\mathrm{def}}{=} \sum_{X \in \Irr(\cal F)} d(X)^2$, where $\Irr(\cal F)$ is the set of isomorphism classes of simple objects. We denote by $F \co Z(\cal F) \to \cal F$ the forgetful functor and by $I \co \cal F \to Z(\cal F)$ its adjoint. Since both categories are semisimple, we have $\dim\operatorname{Hom}(Z,I(V)) = \dim\operatorname{Hom}(I(V),Z) = \dim\operatorname{Hom}(F(Z),V)$ for any objects $Z$ of $Z(\cal F)$ and $V$ of $\cal F$. We denote by $\theta_Z$ the twist of a simple object $Z$ of $Z(\cal F)$, which is a root of unity \cite[p.~9 and Theorem 5.5 p.~13]{NS07}, and by $K(\cal F)$ the Grothendieck ring of $\cal F$. For an object $V$ of $\cal F$ and $n \geq 1$, the Frobenius--Schur indicator $\nu_n(V)$ is the trace of the linear map $E_V^{(n)}$ of $\operatorname{Hom}(1,V^{\ot n})$ defined in \cite[p.~8]{NS07b}. We have $(E^{(n)}_V)^n = \mathrm{id}$ \cite[Theorem 5.1 p.~12]{NS07b}. In particular, $\nu_1(V) = \dim\operatorname{Hom}(1,V)$, and if $\operatorname{Hom}(1,V^{\ot n})$ is one-dimensional, then $\nu_n(V)$ is an $n$-th root of unity. By \cite[Theorem 4.1 p.~10]{NS07}, we have
\begin{equation}
\label{eq:NS-formula}
\nu_n(V) = \frac{1}{\dim \cal F}\sum_{Z \in \Irr(Z(\cal F))} \theta_Z^n  d(Z) \dim \operatorname{Hom}(F(Z),V).
\end{equation}
We equip $Z(\cal F)$ with the braiding of \cite[Theorem 4.1 p.~10]{NS07}. For the reverse braiding, the twists are replaced by their complex conjugates \cite[Remark 4.3 p.~10]{NS07}. Since the right-hand sides of \eqref{eq:moment-0}--\eqref{eq:moment-2} below are real, the only effect of this replacement on Lemma \ref{lem:moments} and Proposition \ref{prop:twisted} is that $\nu_3(\alpha)$ is replaced by $\ovl{\nu_3(\alpha)}$ in the identities $\zeta^3 = \theta_{Z_i}^3 = \nu_3(\alpha)$.

\begin{thm}
\label{thm:recognition}
Let $\cal F$ be a unitary fusion category with six simple objects up to isomorphism, three of which are invertible, and with noncommutative Grothendieck ring.
\begin{enumerate}
\item There exist an invertible object $\alpha$ of order three and a simple object $\rho$ such that the simple objects of $\cal F$ are $1$, $\alpha$, $\alpha^2$, $\rho$, $\alpha\rho$, $\alpha^2\rho$ up to isomorphism, the objects $\alpha^j\rho$ are self-dual, and
\[
\rho\alpha \cong \alpha^2\rho,\quad \rho^2 \cong 1 \oplus \rho \oplus \alpha\rho \oplus \alpha^2\rho.
\]
In particular, $\cal F$ has the fusion rules \eqref{eq:H2-rules} of $\cal H_2$ and $\cal H_3$, and $d(\rho) = d$.
\item If the fusion subcategory of $\cal F$ generated by $\alpha$ is equivalent to $\mathrm{Vec}_{\Z_3}$ with trivial associator, then $\cal F$ is equivalent to $\cal H_2$ or to $\cal H_3$.
\end{enumerate}
Conversely, $\cal H_2$ and $\cal H_3$ satisfy the hypotheses of part 2.
\end{thm}

%%%%%%%%%%%%%%%%%%%%%%%%%%%%%%%%%%%%%%%%%%%%%%%%%%%%%%%%%%%%%%%%%%%%%%%%%%%%%
\subsection{The fusion rules}
\label{subsec:recognition-rules}

\begin{lemma}
\label{lem:ring}
Let $\cal F$ be as in Theorem \ref{thm:recognition}. Then there exist an invertible object $\alpha$ of order three, a simple object $\rho$ and an integer $m \geq 1$ such that the simple objects of $\cal F$ are $1$, $\alpha$, $\alpha^2$, $\rho$, $\alpha\rho$, $\alpha^2\rho$, up to isomorphism, and
\[
\rho\alpha \cong \alpha^2\rho,\quad (\alpha^j\rho)^* \cong \alpha^j\rho \ (j \in \{0,1,2\}),\quad \rho^2 \cong 1 \oplus m(\rho \oplus \alpha\rho \oplus \alpha^2\rho).
\]
Consequently, $d(\rho) = \frac{3m+\sqrt{9m^2+4}}{2}$ and $d(\rho)^2 = 1+3m d(\rho)$.
\end{lemma}

\begin{proof}
The invertible simple objects form a group of order three, generated by some $\alpha$. Let $\cal X$ be the set of the isomorphism classes of the three non-invertible simple objects. The group acts on $\cal X$ by left and by right multiplication.

The algebra $K(\cal F) \ot \Cbb$ is semisimple. Indeed, consider the basis $\Irr(\cal F)$, the involution $X \mapsto X^*$ extended antilinearly, and the linear form $\tau$ which gives the coefficient of $1$. We have $\tau(X^*Y) = \delta_{X,Y}$ for simple $X$, $Y$, so that $\tau(x^*x) > 0$ for $x \neq 0$. Hence $\langle x,y\rangle \ov{\mathrm{def}}{=} \tau(y^*x)$ defines a scalar product on $K(\cal F) \ot \Cbb$, for which left multiplication is a faithful $*$-representation. Thus $K(\cal F) \ot \Cbb$ is a finite-dimensional $C^*$-algebra, hence semisimple. Since $K(\cal F) \ot \Cbb$ has dimension six and is noncommutative, it is isomorphic to $\Cbb \oplus \Cbb \oplus \M_2(\Cbb)$, which has two characters.

Assume that $\alpha X \cong X$ for every $X \in \cal X$. For $X \in \cal X$, we have $X^* \in \cal X$, so that $(X\alpha)^* \cong \alpha^2 X^* \cong X^*$ and $X\alpha \cong X$. Hence, for $X, Y \in \cal X$, the objects $XY$ and $\alpha XY$ are isomorphic, and the invertible summands of $XY$ form a multiple of $1 \oplus \alpha \oplus \alpha^2$. For a primitive cube root of unity $\omega$, the linear map $\chi_\omega \co K(\cal F) \to \Cbb$ with $\chi_\omega(\alpha^j) = \omega^j$ and $\chi_\omega(X) = 0$ for $X \in \cal X$ is therefore a character. The two characters $\chi_\omega$, $\chi_{\omega^2}$ and the Frobenius--Perron character are three distinct characters, a contradiction. Since the orbits of a group of order three have one or three elements, the left action of the group on $\cal X$ is free, and $\cal X = \{\rho,\alpha\rho,\alpha^2\rho\}$ for any $\rho \in \cal X$. If $\rho\alpha \cong \rho$, then $\alpha^2\rho^* \cong (\rho\alpha)^* \cong \rho^*$, which contradicts the freeness. Hence the right action is free as well.

Write $\rho\alpha \cong \alpha^s\rho$. We have $s \neq 0$, since the right action is free. If $s = 1$, then $\alpha$ commutes with every simple object, and $(\alpha^i\rho)(\alpha^j\rho) \cong \alpha^{i+j}\rho^2 \cong (\alpha^j\rho)(\alpha^i\rho)$, so that $K(\cal F)$ is commutative, a contradiction. Hence $s = 2$, and $\rho\alpha^i \cong \alpha^{2i}\rho$ for any $i$. Write $\rho^* \cong \alpha^k\rho$. Then $(\alpha^j\rho)^* \cong \rho^*\alpha^{-j} \cong \alpha^k\rho\alpha^{-j} \cong \alpha^{k+j}\rho$. Applying this twice, $\alpha^j\rho \cong (\alpha^j\rho)^{**} \cong \alpha^{2k+j}\rho$, so that $k = 0$ by freeness, and every non-invertible simple object is self-dual.

Since $\rho$ is self-dual, $\dim\operatorname{Hom}(\alpha^i,\rho^2) = \dim\operatorname{Hom}(\alpha^i\rho,\rho) = \delta_{i,0}$, so that $\rho^2 \cong 1 \oplus \bigoplus_j n_j\alpha^j\rho$ for some integers $n_j \geq 0$. Moreover $\rho^2\alpha \cong \rho\alpha^2\rho \cong \alpha\rho^2$. The coefficients of $\alpha^k\rho$ in $\alpha\rho^2$ and in $\rho^2\alpha$ are $n_{k-1}$ and $n_{k-2}$, so that $n_0 = n_1 = n_2 \ov{\mathrm{def}}{=} m$. If $m = 0$, then $\rho$ is invertible, a contradiction. Taking Frobenius--Perron dimensions gives $d(\rho)^2 = 1+3m d(\rho)$.
\end{proof}

By Lemma \ref{lem:ring}, $\cal F$ is a quadratic category with $(\Z_3,-1,m)$ in the sense of \cite[Definition 2.1 p.~6]{Izu18}. We put $d_m \ov{\mathrm{def}}{=} d(\rho)$ and $\Sigma_\rho \ov{\mathrm{def}}{=} \rho \oplus \alpha\rho \oplus \alpha^2\rho$. Then $d_m - d_m^{-1} = 3m$, $d_m+d_m^{-1} = \sqrt{9m^2+4}$ and
\begin{equation}
\label{eq:dim-F-m}
\dim \cal F = 3+3d_m^2 = 3(2+3m d_m).
\end{equation}
For $m = 1$, we recover $d_1 = d$ and the fusion rules \eqref{eq:H2-rules}.

%%%%%%%%%%%%%%%%%%%%%%%%%%%%%%%%%%%%%%%%%%%%%%%%%%%%%%%%%%%%%%%%%%%%%%%%%%%%%
\subsection{The Drinfeld center}
\label{subsec:recognition-center}

\begin{lemma}
\label{lem:recognition-center}
Let $\cal F$ be as in Theorem \ref{thm:recognition}, with the notation of Lemma \ref{lem:ring}. The simple objects of $Z(\cal F)$ are, up to isomorphism, of the following three kinds.
\begin{enumerate}
\item The summands of $I(1) \cong 1 \oplus B \oplus 2X$, where $B$ and $X$ are simple and non-isomorphic, $F(B) \cong 1 \oplus m\Sigma_\rho$, $F(X) \cong 2\cdot 1 \oplus m\Sigma_\rho$, $d(B) = d_m^2$, $d(X) = 1+d_m^2$ and $\theta_B = \theta_X = 1$.
\item The summands of $I(\alpha) \cong I(\alpha^2) \cong Z_1 \oplus Z_2 \oplus Z_3$, where the $Z_i$ are simple and pairwise non-isomorphic, and $F(Z_i) \cong \alpha \oplus \alpha^2 \oplus k_i\Sigma_\rho$ with integers $k_i \geq 0$ such that $k_1+k_2+k_3 = 3m$.
\item Simple objects $Y_j$, $j \in J$, with $F(Y_j) \cong n_j\Sigma_\rho$ for some integers $n_j \geq 1$.
\end{enumerate}
\end{lemma}

\begin{proof}
Let $Z$ be a simple object of $Z(\cal F)$. The half-braiding gives $\alpha F(Z) \cong F(Z)\alpha$, so that $F(Z)$ is invariant under conjugation by $\alpha$. By Lemma \ref{lem:ring}, $\alpha(\alpha^j\rho)\alpha^{-1} \cong \alpha^{j+1}\rho\alpha^2 \cong \alpha^{j+2}\rho$, so that conjugation by $\alpha$ permutes $\rho$, $\alpha\rho$, $\alpha^2\rho$ cyclically. Hence the multiplicities of $\rho$, $\alpha\rho$ and $\alpha^2\rho$ in $F(Z)$ are equal. Consequently, $Z$ is a summand of $I(1)$, of $I(\alpha)$ or of $I(\alpha^2)$, or $F(Z)$ is a multiple of $\Sigma_\rho$.

By \cite[Proposition 5.4 p.~21]{ENO05}, we have $F(I(V)) \cong \bigoplus_{Y} Y \ot V \ot Y^*$, the sum running over $\Irr(\cal F)$. For $V = 1$, Lemma \ref{lem:ring} gives $F(I(1)) \cong 3\cdot 1 \oplus 3\rho^2 \cong 6\cdot 1 \oplus 3m\Sigma_\rho$. For $V = \alpha$ and $Y = \alpha^j\rho$, we have $Y^* \cong \rho\alpha^{-j}$ and $Y\alpha Y^* \cong \alpha^j\rho\alpha\rho\alpha^{-j} \cong \alpha^{j+2}\rho^2\alpha^{-j} \cong \alpha^2\rho^2$, since $\rho^2\alpha \cong \alpha\rho^2$ by the proof of Lemma \ref{lem:ring}. Hence $F(I(\alpha)) \cong 3\alpha \oplus 3\alpha^2 \oplus 3m\Sigma_\rho$, and the same computation gives $F(I(\alpha^2)) \cong 3\alpha \oplus 3\alpha^2 \oplus 3m\Sigma_\rho$. By adjunction, we obtain
\begin{align*}
\dim \operatorname{End}(I(\alpha)) &= \dim \operatorname{End}(I(\alpha^2)) = \dim\operatorname{Hom}(I(\alpha),I(\alpha^2)) = 3,\\
\dim\operatorname{Hom}(I(1),I(\alpha)) &= \dim\operatorname{Hom}(I(1),I(\alpha^2)) = 0.
\end{align*}
The first equalities show that $I(\alpha)$ and $I(\alpha^2)$ are both direct sums of three pairwise non-isomorphic simple objects, and that the three summands $Z_1$, $Z_2$, $Z_3$ of $I(\alpha)$ are summands of $I(\alpha^2)$, so that $I(\alpha^2) \cong I(\alpha)$. The second ones show that $F(Z_i)$ does not contain $1$, and that $F(Z)$ contains neither $\alpha$ nor $\alpha^2$ for any summand $Z$ of $I(1)$. Since $\dim\operatorname{Hom}(F(Z_i),\alpha) = \dim\operatorname{Hom}(F(Z_i),\alpha^2) = 1$, we get $F(Z_i) \cong \alpha \oplus \alpha^2 \oplus k_i\Sigma_\rho$, and $\sum_i k_i = 3m$ by comparison with $F(I(\alpha))$.

By \cite[Theorem 2.13 p.~6]{Ost15}, the simple summands of $I(1)$ are in bijection with the irreducible representations $E$ of $K(\cal F)$, the summand $A_E$ attached to $E$ has multiplicity $\dim E$ and dimension $\frac{\dim \cal F}{f_E}$, where $f_E$ is the formal codegree of $E$. The unit object of $Z(\cal F)$ has multiplicity one in $I(1)$, so that it is $A_E$ for a one-dimensional $E$ with $f_E = \dim \cal F$. A character $\chi$ of $K(\cal F)$ satisfies $\chi(\alpha)\chi(\rho) = \chi(\rho\alpha) = \chi(\alpha)^2\chi(\rho)$ and $\chi(\rho)^2 = 1+m\chi(\rho)(1+\chi(\alpha)+\chi(\alpha)^2)$. The second relation gives $\chi(\rho) \neq 0$, so that $\chi(\alpha) = 1$ and $\chi(\rho) \in \{d_m,-d_m^{-1}\}$. The formal codegrees of these two characters are $\sum_Y \chi(Y)\chi(Y^*) = 3+3\chi(\rho)^2$, that is, $\dim\cal F$ and $3+3d_m^{-2}$. The corresponding summands of $I(1)$ are the unit object and a simple object $B$ of dimension $\frac{\dim\cal F}{3+3d_m^{-2}} = d_m^2$ by \eqref{eq:dim-F-m}. The summand $X$ attached to the two-dimensional representation has multiplicity two, and $d(I(1)) = \dim \cal F = 1+d_m^2+2d(X)$ gives $d(X) = 1+d_m^2$. By adjunction, $1$ has multiplicity one in $F(B)$ and two in $F(X)$. Writing $F(B) \cong 1 \oplus s\Sigma_\rho$ and $F(X) \cong 2\cdot 1 \oplus s'\Sigma_\rho$, the dimensions give $1+3s d_m = d_m^2$ and $2+3s'd_m = 1+d_m^2$, so that $s = s' = m$ by Lemma \ref{lem:ring}. Finally $\theta_{I(1)} = \mathrm{id}_{I(1)}$ by \cite[Theorem 2.5 p.~4]{Ost15}.
\end{proof}

\begin{lemma}
\label{lem:moments}
In the situation of Lemma \ref{lem:recognition-center}, put $\Delta \ov{\mathrm{def}}{=} 9m^2+4$ and $\epsi \ov{\mathrm{def}}{=} \nu_2(\rho)$. Then $\epsi \in \{1,-1\}$ and $k_1 = k_2 = k_3 = m$. The twists of $Z_1$, $Z_2$, $Z_3$ are $\zeta$, $\zeta\omega$, $\zeta\omega^2$ for some root of unity $\zeta$ with $\zeta^3 = \nu_3(\alpha)$ and a primitive cube root of unity $\omega$. Moreover
\begin{align}
\label{eq:moment-0}
\sum_{j \in J} n_j^2 &= 4m^2+2,\\
\label{eq:moment-1}
\sum_{j \in J} n_j^2 \theta_{Y_j} &= -\frac{m^2+m\sqrt{\Delta}}{2},\\
\label{eq:moment-2}
\sum_{j \in J} n_j^2 \theta_{Y_j}^2 &= \frac{-m^2+(2\epsi-m)\sqrt{\Delta}}{2}.
\end{align}
\end{lemma}

\begin{proof}
Since $\rho$ is simple and self-dual, $\operatorname{Hom}(1,\rho^{\ot 2})$ is one-dimensional, so that $\epsi \in \{1,-1\}$. Similarly, $\operatorname{Hom}(1,\alpha^{\ot 3})$ is one-dimensional, so that $|\nu_3(\alpha)| = 1$. Moreover $\nu_1(\alpha) = \nu_1(\rho) = 0$. Put $\ell_i \ov{\mathrm{def}}{=} d(Z_i) = 2+3k_i d_m > 0$, so that $\sum_i \ell_i = \dim\cal F$ by \eqref{eq:dim-F-m}. By \eqref{eq:NS-formula} and Lemma \ref{lem:recognition-center}, we have
\[
\nu_n(\alpha) = \frac{1}{\dim\cal F}\sum_{i=1}^3 \theta_{Z_i}^n \ell_i.
\]
For $n = 3$, this gives $\sum_i \theta_{Z_i}^3\ell_i = \nu_3(\alpha)\sum_i \ell_i$ with $|\nu_3(\alpha)| = 1$. Since the $\ell_i$ are positive and the $\theta_{Z_i}$ have modulus one, we obtain $\theta_{Z_1}^3 = \theta_{Z_2}^3 = \theta_{Z_3}^3 = \nu_3(\alpha)$. Write $\theta_{Z_i} = \zeta u_i$ with $\zeta = \theta_{Z_1}$ and $u_i^3 = 1$. For $n = 1$, we obtain $\sum_i u_i\ell_i = 0$. If two of the $u_i$ coincide, then $\sum_i u_i\ell_i$ is a combination with positive coefficients of two cube roots of unity whose quotient is not real, or a positive multiple of a cube root of unity, which is nonzero. Hence the $u_i$ are the three distinct cube roots of unity, and $1+\omega+\omega^2 = 0$ and the linear independence of $1$ and $\omega$ over $\mathbb{R}$ give $\ell_1 = \ell_2 = \ell_3$, that is, $k_i = m$ for all $i$. In particular $d(Z_i) = 2+3m d_m = 1+d_m^2$, and $\sum_i\theta_{Z_i}^n = 0$ for $n \in \{1,2\}$.

Since $d(Z) = d(F(Z))$ is the Frobenius--Perron dimension of $Z$ for every simple object $Z$ of $Z(\cal F)$, we have $\sum_Z d(Z)^2 = (\dim\cal F)^2 = 9(1+d_m^2)^2$ by \cite[Proposition 8.12 p.~32]{ENO05}, and $d(Y_j) = 3n_j d_m$. By Lemma \ref{lem:recognition-center}, we obtain
\[
9d_m^2\sum_j n_j^2 = 9(1+d_m^2)^2-1-d_m^4-(1+d_m^2)^2-3(1+d_m^2)^2 = 4d_m^4+10d_m^2+4,
\]
and $\sum_j n_j^2 = \frac{4(d_m^2+d_m^{-2})+10}{9} = 4m^2+2$, since $d_m^2+d_m^{-2} = (d_m-d_m^{-1})^2+2 = 9m^2+2$.

By \eqref{eq:NS-formula} for $V = \rho$ and $n \in \{1,2\}$, using $\nu_1(\rho) = 0$, $\nu_2(\rho) = \epsi$, Lemma \ref{lem:recognition-center} and $\sum_i\theta_{Z_i}^n = 0$, we get
\[
m d_m^2+m(1+d_m^2)+3d_m\sum_j n_j^2\theta_{Y_j} = 0,\quad m d_m^2+m(1+d_m^2)+3d_m\sum_j n_j^2\theta_{Y_j}^2 = \epsi\dim\cal F.
\]
Since $1+2d_m^2 = 3+6m d_m$ and $d_m^{-1} = d_m-3m$, the first equality gives $\sum_j n_j^2\theta_{Y_j} = -m(d_m^{-1}+2m) = m^2-m d_m$, which is \eqref{eq:moment-1} since $d_m = \frac{3m+\sqrt\Delta}{2}$. By \eqref{eq:dim-F-m}, the second one gives $\sum_j n_j^2\theta_{Y_j}^2 = \epsi(d_m+d_m^{-1})+m^2-m d_m = \epsi\sqrt\Delta+m^2-m d_m$, which is \eqref{eq:moment-2}.
\end{proof}

The identities \eqref{eq:moment-0}, \eqref{eq:moment-1} and \eqref{eq:moment-2} are those of \cite[(25) p.~7]{Hua26}. For $m = 1$ and the Haagerup categories, the objects $Y_j$ are $d_1, \ldots, d_6$, with $n_j = 1$ and twists $\e^{\frac{2\pi\i 6l^2}{13}}$ by \eqref{eq:center-dimensions} and \eqref{eq:center-twists}, since $d(Y_j) = 3n_jd$. Since $13 \equiv 1 \bmod 4$, and since $12 \equiv -1$ is a square modulo $13$ whereas $6$ is not, the quadratic Gauss sums give $\sum_{l=1}^6\e^{\frac{2\pi\i 6l^2}{13}} = \frac{-1-\sqrt{13}}{2}$ and $\sum_{l=1}^6\e^{\frac{2\pi\i 12l^2}{13}} = \frac{-1+\sqrt{13}}{2}$, which are \eqref{eq:moment-1} and \eqref{eq:moment-2} for $m = 1$ and $\epsi = 1$.

%%%%%%%%%%%%%%%%%%%%%%%%%%%%%%%%%%%%%%%%%%%%%%%%%%%%%%%%%%%%%%%%%%%%%%%%%%%%%
\subsection{Sums of roots of unity}
\label{subsec:roots-of-unity}

For $n \geq 1$, we denote by $r(n)$ the sum of the primitive $n$-th roots of unity, that is, the Möbius function of $n$. Thus $r(n) \in \{0,1,-1\}$, $r(n) = 0$ if $4$ divides $n$, $r(1) = 1$, $r(2) = -1$, and $r(kl) = r(k)r(l)$ if $k$ and $l$ are coprime. We denote by $\varphi$ the Euler function.

Let $D > 1$ be a squarefree integer. For a root of unity $\zeta$, choose $L \geq 1$ such that $\zeta^L = 1$ and $\sqrt D \in \mathbb{Q}(\e^{\frac{2\pi\i}{L}})$, put $G_L \ov{\mathrm{def}}{=} \operatorname{Gal}(\mathbb{Q}(\e^{\frac{2\pi\i}{L}})/\mathbb{Q}(\sqrt D))$ and
\[
\langle \zeta \rangle_D \ov{\mathrm{def}}{=} \frac{1}{|G_L|}\sum_{\sigma \in G_L}\sigma(\zeta) = A_D(\zeta)+B_D(\zeta)\sqrt D,\quad A_D(\zeta), B_D(\zeta) \in \mathbb{Q}.
\]
This does not depend on $L$, since the restriction maps $G_{L'} \to G_L$ for $L \mid L'$ are surjective with fibers of the same cardinality. If $\zeta_1, \ldots, \zeta_N$ are roots of unity with $\sum_q \zeta_q = a+b\sqrt D$ for rational numbers $a$, $b$, then averaging over $G_L$ for a common $L$ gives
\begin{equation}
\label{eq:average-sum}
\sum_{q=1}^N A_D(\zeta_q) = a,\quad \sum_{q=1}^N B_D(\zeta_q) = b.
\end{equation}

\begin{lemma}
\label{lem:galois-average}
Let $\zeta$ be a root of unity of order $n$.
\begin{enumerate}
\item If $40$ does not divide $n$, then $\langle\zeta\rangle_{10} = \frac{r(n)}{\varphi(n)}$.
\item If $n = 40k$ with $\gcd(k,40) > 1$, then $\langle\zeta\rangle_{10} = 0$.
\item If $n = 40k$ with $\gcd(k,40) = 1$, then $A_{10}(\zeta) = 0$ and $|B_{10}(\zeta)| = \frac{|r(k)|}{8\varphi(k)}$.
\end{enumerate}
Consequently, $|A_{10}(\zeta)|+8|B_{10}(\zeta)| \leq 1$, with equality if and only if $\zeta \in \{1,-1\}$ or $n = 40$. If $\zeta_1, \ldots, \zeta_N$ are roots of unity with $\sum_q\zeta_q = a+b\sqrt{10}$ for rational numbers $a$, $b$, then $N \geq |a|+8|b|$.
\end{lemma}

\begin{proof}
The field $\mathbb{Q}(\sqrt{10})$ has discriminant $40$, so that $\sqrt{10} \in \mathbb{Q}(\zeta)$ if and only if $40$ divides $n$. If $40 \nmid n$, then $\mathbb{Q}(\zeta) \cap \mathbb{Q}(\sqrt{10}) = \mathbb{Q}$, the group $G_L$ maps onto $\operatorname{Gal}(\mathbb{Q}(\zeta)/\mathbb{Q})$, and $\langle\zeta\rangle_{10}$ is the average of the primitive $n$-th roots of unity, which is $\frac{r(n)}{\varphi(n)}$.

Assume that $n = 40k$, and take $L = n$. We identify $\operatorname{Gal}(\mathbb{Q}(\e^{\frac{2\pi\i}{n}})/\mathbb{Q})$ with $(\Z/n\Z)^\times$. The subgroup $G_n$ has index two, and we denote by $\chi$ the character of $(\Z/n\Z)^\times$ with kernel $G_n$. It factors through $(\Z/40\Z)^\times$, since the kernel of the map $(\Z/n\Z)^\times \to (\Z/40\Z)^\times$ fixes $\mathbb{Q}(\e^{\frac{2\pi\i}{40}})$, which contains $\sqrt{10}$. Put $\tau_n \ov{\mathrm{def}}{=} \sum_{u \in (\Z/n\Z)^\times}\chi(u)\e^{\frac{2\pi\i u}{n}}$. Writing $\zeta = \e^{\frac{2\pi\i v}{n}}$ with $v \in (\Z/n\Z)^\times$, we obtain
\[
\langle\zeta\rangle_{10} = \frac{2}{\varphi(n)}\sum_{u \in G_n}\e^{\frac{2\pi\i uv}{n}} = \frac{1}{\varphi(n)}\sum_{u \in (\Z/n\Z)^\times}(1+\chi(u))\e^{\frac{2\pi\i uv}{n}} = \frac{r(n)+\chi(v)\tau_n}{\varphi(n)} = \frac{\chi(v)\tau_n}{\varphi(n)},
\]
since $r(n) = 0$. For $\sigma \notin G_n$ we have $\sigma(\tau_n) = -\tau_n$, so that $\tau_n \in \mathbb{Q}\sqrt{10}$ and $A_{10}(\zeta) = 0$.

If $\gcd(k,40) > 1$, choose a prime $p$ dividing $k$ and $40$. Then $p$ divides $\frac np$ and $40$ divides $\frac np$. Hence $u \mapsto u+\frac np$ permutes $(\Z/n\Z)^\times$ and preserves $\chi$, and it multiplies $\e^{\frac{2\pi\i u}{n}}$ by $\e^{\frac{2\pi\i}{p}} \neq 1$, so that $\tau_n = 0$. This proves part 2.

Assume that $\gcd(k,40) = 1$. By the Chinese remainder theorem, $\e^{\frac{2\pi\i}{n}} = \e^{\frac{2\pi\i x}{40}}\e^{\frac{2\pi\i y}{k}}$ with $x \in (\Z/40\Z)^\times$ and $y \in (\Z/k\Z)^\times$, and $(\Z/n\Z)^\times \cong (\Z/40\Z)^\times \times (\Z/k\Z)^\times$. Hence
\[
\tau_n = \Big(\sum_{u \in (\Z/40\Z)^\times}\chi(u)\e^{\frac{2\pi\i xu}{40}}\Big)\Big(\sum_{u \in (\Z/k\Z)^\times}\e^{\frac{2\pi\i yu}{k}}\Big) = \chi(x)\tau_{40} r(k).
\]
Let $s_+$ and $s_-$ be the sums of the primitive $40$-th roots of unity in the two cosets of $G_{40}$, so that $\tau_{40} = \pm(s_+-s_-)$. They are fixed by $G_{40}$ and exchanged by the other elements of the Galois group. By \cite[proof of Theorem 6.2, Case 2, p.~18]{Lar14}, applied with $c = 10$, we have $\{s_+,s_-\} = \{\sqrt{10},-\sqrt{10}\}$. Consequently, $|\tau_{40}| = 2\sqrt{10}$, and $|B_{10}(\zeta)|\sqrt{10} = \frac{2\sqrt{10} |r(k)|}{\varphi(40)\varphi(k)} = \frac{2\sqrt{10} |r(k)|}{16\varphi(k)}$, which is part 3.

By parts 1, 2 and 3, $|A_{10}(\zeta)|+8|B_{10}(\zeta)|$ is $\frac{|r(n)|}{\varphi(n)}$, $0$ or $\frac{|r(k)|}{\varphi(k)}$ in the three cases. This is at most $1$, with equality if and only if $n \in \{1,2\}$ in the first case and $k = 1$ in the third case, since $k$ is then odd. The last assertion follows from \eqref{eq:average-sum}, since $|a|+8|b| \leq \sum_q(|A_{10}(\zeta_q)|+8|B_{10}(\zeta_q)|) \leq N$.
\end{proof}

For a general squarefree integer $D > 1$ and rational numbers $a$, $b$ with $2a, 2b \in \Z$, a sum of $N$ roots of unity equal to $a+b\sqrt D$ satisfies
\begin{equation}
\label{eq:larson}
N \geq |b|\varphi(2D).
\end{equation}
This follows from \cite[Theorem 1.3 p.~2 and Theorem 6.2 p.~18]{Lar14}, applied to the sum of $2N$ roots of unity obtained by repeating each root twice, which is $2a+2b\sqrt D$ with $2a, 2b \in \Z$.

%%%%%%%%%%%%%%%%%%%%%%%%%%%%%%%%%%%%%%%%%%%%%%%%%%%%%%%%%%%%%%%%%%%%%%%%%%%%%
\subsection{The multiplicity}
\label{subsec:multiplicity}

\begin{prop}
\label{prop:m-equals-one}
In Lemma \ref{lem:ring}, we have $m = 1$.
\end{prop}

\begin{proof}
Write $\Delta = 9m^2+4 = t^2D$ with $t \geq 1$ and $D$ squarefree. Since $(3m)^2 < \Delta < (3m+1)^2$, we have $D > 1$. Since $\Delta \equiv 1 \bmod 3$, the integer $t$ is prime to $3$ and $D \equiv 1 \bmod 3$. An odd prime $p$ dividing $D$ divides $(3m)^2+4$ and does not divide $3m$, so that $-1$ is a square modulo $p$ and $p \equiv 1 \bmod 4$. Hence every prime factor of $D$ is $2$ or congruent to $1$ modulo $4$, and $D \notin \{2,5\}$ since $D \equiv 1 \bmod 3$. Write $D = 2^e\prod_i p_i$ with distinct odd primes $p_i$ and $e \in \{0,1\}$. Then
\[
\frac{\varphi(2D)}{\sqrt D} = 2^{\frac e2}\prod_i\frac{p_i-1}{\sqrt{p_i}},
\]
with $\frac{p-1}{\sqrt p} \geq \frac{4}{\sqrt 5} > 1$ for $p \geq 5$ and $\frac{p-1}{\sqrt p} \geq \frac{12}{\sqrt{13}} > 3$ for $p \geq 13$. If some $p_i$ is at least $13$, the product is at least $\frac{12}{\sqrt{13}} > 3$, since the other factors are at least $1$. Otherwise every $p_i$ equals $5$, so that $D \in \{2,5,10\}$, and $D = 10$. Hence $\frac{\varphi(2D)}{\sqrt D} > 3$ unless $D = 10$.

Consider the family of $N = 4m^2+2$ roots of unity in which each twist $\theta_{Y_j}$ of Lemma \ref{lem:moments} is repeated $n_j^2$ times. By \eqref{eq:moment-0} and \eqref{eq:moment-1}, its sum is $a+b\sqrt D$ with $a = -\frac{m^2}{2}$ and $b = -\frac{mt}{2}$. Note that $t\sqrt D = \sqrt\Delta > 3m$.

Assume that $D \neq 10$. By \eqref{eq:larson}, we have $4m^2+2 \geq \frac{mt}{2}\varphi(2D) = \frac m2  t\sqrt D \frac{\varphi(2D)}{\sqrt D} > \frac{9m^2}{2}$, so that $m^2 < 4$ and $m = 1$.

Assume that $D = 10$. By Lemma \ref{lem:galois-average}, we have $4m^2+2 \geq \frac{m^2}{2}+4mt > \frac{m^2}{2}+\frac{12m^2}{\sqrt{10}}$, so that $m^2 < 7$. Since $D = 13$ for $m = 1$, we get $m = 2$ and $t = 2$. Then $N = 18 = |a|+8|b|$. By \eqref{eq:average-sum}, we have $\sum_q(-A_{10}(\zeta_q)-8B_{10}(\zeta_q)) = 18$, and each term is at most $1$ by Lemma \ref{lem:galois-average}. Hence $-A_{10}(\zeta_q)-8B_{10}(\zeta_q) = 1$ for every root $\zeta_q$ of the family. Since this quantity is $-1$ for $\zeta_q = 1$, Lemma \ref{lem:galois-average} shows that each $\zeta_q$ is equal to $-1$ or is a primitive $40$-th root of unity, and that $A_{10}(\zeta_q) = 0$ in the second case. Since $\sum_q A_{10}(\zeta_q) = a = -2$ by \eqref{eq:average-sum}, the root $-1$ occurs exactly twice in the family. Hence the family of the squares $\zeta_q^2$ consists of $1$, twice, and of $16$ primitive $20$-th roots of unity, and $A_{10}(\zeta_q^2) = \frac{r(20)}{\varphi(20)} = 0$ for the latter by part 1 of Lemma \ref{lem:galois-average}. By \eqref{eq:moment-2}, the sum of this family is $\frac{-4+(2\epsi-2)\sqrt{40}}{2} = -2+2(\epsi-1)\sqrt{10}$. Comparing the rational parts by \eqref{eq:average-sum}, we obtain $-2 = \sum_q A_{10}(\zeta_q^2) = 2$, a contradiction. This excludes $m = 2$.
\end{proof}

%%%%%%%%%%%%%%%%%%%%%%%%%%%%%%%%%%%%%%%%%%%%%%%%%%%%%%%%%%%%%%%%%%%%%%%%%%%%%
\subsection{Pointed parts and proof of the recognition theorem}
\label{subsec:recognition-proof}

\begin{lemma}
\label{lem:pointed-trivial}
Let $\cal P$ be a fusion category with three simple objects, all invertible, so that $\cal P$ is equivalent to $\mathrm{Vec}^\varpi_{\Z_3}$ for some $3$-cocycle $\varpi$.
\begin{enumerate}
\item If a left $\cal P$-module category $\cal M$ contains a simple object $M$ with $g \ot M \cong M$ for every invertible object $g$ of $\cal P$, then $\varpi$ is a coboundary.
\item If $\cal P$ has two inequivalent indecomposable left module categories, then $\varpi$ is a coboundary.
\end{enumerate}
\end{lemma}

\begin{proof}
1. The full subcategory of $\cal M$ whose objects are the multiples of $M$ is stable under the action of $\cal P$, and it is equivalent to $\mathrm{Vec}$. A structure of $\cal P$-module category on $\mathrm{Vec}$ is the same thing as a fiber functor $\cal P \to \mathrm{Vec}$ \cite[Example (ii) pp.~6--7]{Ost03}. The monoidal structure of a fiber functor on $\mathrm{Vec}^\varpi_{\Z_3}$ consists of nonzero scalars $j(g,h)$, after a choice of basis vectors of the lines attached to the simple objects, and its compatibility with the associators says that $\varpi$ is the coboundary of $j$.

2. Let $\cal M$ be an indecomposable $\cal P$-module category. The invertible objects permute the isomorphism classes of simple objects of $\cal M$. The full subcategory whose objects are the direct sums of simple objects in one orbit is a module subcategory and a direct summand of $\cal M$, so that the action is transitive, with orbits of cardinality one or three. If a simple object is fixed, then $\varpi$ is a coboundary by part 1. Otherwise the action is free. Let $M$ be a simple object of $\cal M$, and let $\underline{\operatorname{Hom}}(M,M)$ be its internal Hom, characterized by $\operatorname{Hom}_{\cal P}(Y,\underline{\operatorname{Hom}}(M,M)) \cong \operatorname{Hom}_{\cal M}(Y \ot M,M)$ \cite[Definition 11 p.~9]{Ost03}. For $Y = g$ invertible, this space is nonzero if and only if $g \cong 1$, so that $\underline{\operatorname{Hom}}(M,M) \cong 1$ as an object. Its algebra structure is then isomorphic to the trivial algebra $1$, and by \cite[Theorem 1 and its proof pp.~10--11]{Ost03}, $\cal M$ is equivalent to the category of modules over $1$, that is, to the regular module category $\cal P$. Hence, if $\varpi$ is not a coboundary, every indecomposable $\cal P$-module category is equivalent to $\cal P$.
\end{proof}

We now prove Theorem \ref{thm:recognition}.

\begin{proof}
Part 1 follows from Lemma \ref{lem:ring} and Proposition \ref{prop:m-equals-one}, with $d(\rho) = d_1 = d$.

We next prove the converse statement. The categories $\cal H_2$ and $\cal H_3$ are unitary fusion categories with the fusion rules \eqref{eq:H2-rules} \cite[p.~14, Theorem 1.4 p.~4 and Corollary 3.24 p.~24]{GrS12}, the category $\cal H_3$ being the category of bimodules over a $\Q$-system in $\cal H_2$ \cite[Lemma 3.14 p.~21 and Definition 3.19 p.~23]{GrS12}. They have six simple objects, three of which are invertible, and their Grothendieck rings are isomorphic and noncommutative since $\xi\alpha \cong \alpha^2\xi \not\cong \alpha\xi$ in $\cal H_2$. The pointed part of $\cal H_2$ has trivial associator by \cite[Lemma 2.21 p.~15]{GrS12}. For $\cal H_3$, let $A \ov{\mathrm{def}}{=} 1 \oplus \alpha \oplus \alpha^2$ with this algebra structure, so that $\cal H_3$ is the category of $A$-$A$-bimodules in $\cal H_2$. The category $\cal M$ of left $A$-modules in $\cal H_2$ is a left $\cal H_3$-module category for the relative tensor product $\ot_A$ and a right $\cal H_2$-module category, and the two actions commute. Since the $\Q$-system $A$ is separable, the category $\cal M$ is semisimple \cite[Proposition 2.7 p.~5]{DMNO13}. For objects $Y$, $Y'$ of $\cal H_2$, we have $\operatorname{Hom}_A(A \ot Y,A \ot Y') \cong \operatorname{Hom}(Y,A \ot Y')$, and every simple left $A$-module is a quotient, hence a summand, of some $A \ot Y$ with $Y$ simple \cite[Lemma 4 p.~9]{Ost03}. By \eqref{eq:H2-rules}, we deduce that the simple objects of $\cal M$ are $A \cong A \ot \alpha \cong A \ot \alpha^2$ and $A \ot \xi \cong A \ot \alpha\xi \cong A \ot \alpha^2\xi$, up to isomorphism. Moreover $(A \ot \xi) \ot \xi \cong A \oplus 3(A \ot \xi)$ is not simple. Hence $A$ is, up to isomorphism, the only simple object $M$ of $\cal M$ such that $M \ot \xi$ is simple. For an invertible object $g$ of $\cal H_3$, the functor $g \ot_A -$ is an autoequivalence of $\cal M$ commuting with the right action of $\cal H_2$, so that $(g \ot_A A) \ot \xi \cong g \ot_A (A \ot \xi)$ is simple and $g \ot_A A \cong A$. By part 1 of Lemma \ref{lem:pointed-trivial}, the pointed part of $\cal H_3$ has trivial associator.

We finally prove part 2. By part 1, $\cal F$ is a quadratic category with $(\Z_3,-1,1)$ \cite[Definition 2.1 p.~6]{Izu18}. Its invariant $c_{0,3}(\cal F) \in H^3(\Z_3,\Cbb^\times)$ is the class of the associator of its pointed part \cite[p.~9]{Izu18}, which vanishes by assumption. The invariant $c_{1,2}(\cal F)$ belongs to $H^1(\Z_2,H^2(\Z_3,\Cbb^\times))$ \cite[pp.~10--11]{Izu18}, which is zero since $H^2(\Z_3,\Cbb^\times) = 0$. Hence $\cal F$ is a generalized Haagerup category with group $\Z_3$ \cite[Definition 2.7 p.~12]{Izu18} which is a $C^*$-fusion category. For $G = \Z_3$, the group $\Gamma = (H^2(G,\T) \times G/2G) \rtimes \operatorname{Aut}(G)$ of \cite[Theorem 5.9 p.~36]{Izu18} is $\operatorname{Aut}(\Z_3)$, and since $G$ is odd one may take the signs $\epsi_h(g)$ of \cite{Izu18} equal to $1$ \cite[p.~56]{Izu18}. Hence the number of $\Gamma$-orbits of gauge equivalence classes of solutions in \cite[Theorem 5.9 p.~36]{Izu18} is at most the number of solutions up to group automorphism in \cite[Theorem 9.1 p.~57]{Izu18}, and there are at most two $C^*$-generalized Haagerup categories with group $\Z_3$ up to equivalence. By the converse statement, $\cal H_2$ and $\cal H_3$ are two such categories, and they are not equivalent by part 1 of Proposition \ref{prop:GS}. Hence $\cal F$ is equivalent to $\cal H_2$ or to $\cal H_3$.
\end{proof}

%%%%%%%%%%%%%%%%%%%%%%%%%%%%%%%%%%%%%%%%%%%%%%%%%%%%%%%%%%%%%%%%%%%%%%%%%%%%%
\subsection{The associator of the pointed part and the twists}
\label{subsec:twisted}

\begin{lemma}
\label{lem:nu3}
Let $\cal P$ be a fusion category with three simple objects, all invertible, equipped with the spherical structure for which all simple objects have dimension $1$, and let $\alpha$ be a non-trivial simple object. Then $\nu_3(\alpha) = 1$ if the associator of $\cal P$ is trivial, and $\nu_3(\alpha)$ is a primitive cube root of unity otherwise.
\end{lemma}

\begin{proof}
The category $\cal P$ is pseudo-unitary, since $\dim \cal P = \operatorname{FPdim} \cal P = 3$. Hence its spherical structure with positive dimensions is unique \cite[Proposition 8.23 p.~34]{ENO05}. By \cite[Corollary 4.4 p.~11]{NS07b}, Frobenius--Schur indicators are invariant under tensor equivalences which preserve the pivotal structures, and every tensor equivalence between such pointed categories preserves the spherical structures with positive dimensions. Hence we may assume that $\cal P$ is the category $\cal C(\Z/3\Z,\omega_t)$ of \cite[p.~30]{NS07} with the pivotal structure of dimension $1$ described there, where $t \in \{0,1,2\}$ and the class of $\omega_t$ is trivial if and only if $t = 0$, and that $\alpha$ is the object $V_1$ or $V_2$. By \cite[p.~31]{NS07}, we have $\nu_3(V_1) = \e^{\frac{2\pi\i t}{3}}$. Since $V_2 \cong V_1^*$, we also have $\nu_3(V_2) = \nu_3(V_1)$ \cite[Remark 4.2 p.~10]{NS07}. Hence $\nu_3(\alpha) = \e^{\frac{2\pi\i t}{3}}$, which is equal to $1$ if and only if $t = 0$.
\end{proof}

\begin{prop}
\label{prop:twisted}
Let $\cal F$ be as in Theorem \ref{thm:recognition}, and let $Z_1$, $Z_2$, $Z_3$ be the simple summands of $I(\alpha)$ (Lemma \ref{lem:recognition-center}). Then $\theta_{Z_i}^3 = \nu_3(\alpha)$ for $i \in \{1,2,3\}$.
\begin{enumerate}
\item If the fusion subcategory generated by $\alpha$ has trivial associator, then $\theta_{Z_i}^3 = 1$, and $\cal F$ is equivalent to $\cal H_2$ or to $\cal H_3$.
\item Otherwise the twists $\theta_{Z_i}$ are primitive ninth roots of unity. In particular, $Z(\cal F)$ is not braided equivalent to $Z(\cal H_3)$.
\end{enumerate}
Consequently, $\cal F$ is equivalent to $\cal H_2$ or to $\cal H_3$ if and only if no simple object of $Z(\cal F)$ has a twist of order nine.
\end{prop}

\begin{proof}
The equality $\theta_{Z_i}^3 = \nu_3(\alpha)$ is proved in Lemma \ref{lem:moments}. The indicator $\nu_3(\alpha)$ is computed in the fusion subcategory $\cal P$ generated by $\alpha$, since $\operatorname{Hom}(1,\alpha^{\ot 3})$ and the map $E^{(3)}_\alpha$ only involve objects of $\cal P$ \cite[p.~31]{NS07}, and the spherical structure of $\cal F$ restricts to the spherical structure of $\cal P$ with dimensions $1$. Part 1 then follows from Lemma \ref{lem:nu3} and part 2 of Theorem \ref{thm:recognition}. In the situation of part 2, $\nu_3(\alpha)$ is a primitive cube root of unity by Lemma \ref{lem:nu3}, so that $\theta_{Z_i}^9 = 1 \neq \theta_{Z_i}^3$. A braided equivalence between pseudo-unitary modular tensor categories preserves the spherical structures with positive dimensions \cite[Proposition 8.23 p.~34]{ENO05}, hence the twists. By \eqref{eq:center-twists}, the twists of $Z(\cal H_3)$ have order $1$, $3$ or $13$, which proves part 2. The last assertion follows from parts 1 and 2, since $Z(\cal H_2)$ and $Z(\cal H_3)$ are braided equivalent to $Z(\cal H_3)$ (Section \ref{subsec:haagerup-categories}).
\end{proof}

%%%%%%%%%%%%%%%%%%%%%%%%%%%%%%%%%%%%%%%%%%%%%%%%%%%%%%%%%%%%%%%%%%%%%%%%%%%%%
\subsection{Application to the conformal nets}
\label{subsec:recognition-nets}

\begin{lemma}
\label{lem:trivial-associator}
Under the hypotheses of Theorem \ref{thm:complete-rationality}, let $x \in \{a,c_0\}$. Then the pointed fusion subcategory $\cal P_x$ of $\cal F_x$ is equivalent to $\mathrm{Vec}_{\Z_3}$ with trivial associator.
\end{lemma}

\begin{proof}
By part 4 of Theorem \ref{thm:complete-rationality}, part 1 of Proposition \ref{prop:DMNO} and Lemma \ref{lem:rank}, there is a braided equivalence $\Rep(\cal A_U) \cong Z(\cal F_x)$ under which $\Theta_x$ corresponds to $I(1)$. By the proof of part 4 of Theorem \ref{thm:complete-rationality}, the connected \'etale subalgebra $\theta = \pi_0 \oplus \pi_b$ of $\Theta_x$ corresponds to the \'etale subalgebra $A_x$ of $I(1)$ attached to $\cal P_x$ by \cite[Theorem 4.10 p.~18]{DMNO13}. By part 5 of Proposition \ref{prop:DMNO}, the category of local modules $Z(\cal F_x)^0_{A_x}$ is braided equivalent to $Z(\cal P_x)$. On the other hand, the category of local $\theta$-modules $\Rep(\cal A_U)^0_\theta$ is braided equivalent to $\Rep(\cal A_{V'})$ \cite[Proposition 6.4 p.~33]{BKL15}, see also \cite[Theorem 3.1 p.~2]{Mug10}. Hence there is a braided equivalence $\Psi_x \co \Rep(\cal A_{V'}) \to Z(\cal P_x)$. If the braiding conventions of \cite[Proposition 6.4 p.~33]{BKL15} and of part 5 of Proposition \ref{prop:DMNO} differ, then $\Psi_x$ is a braided equivalence onto $Z(\cal P_x)^{\mathrm{rev}}$. This does not affect the argument below, since an algebra which is commutative for a braiding is commutative for the reverse braiding, so that $Z(\cal P_x)$ and $Z(\cal P_x)^{\mathrm{rev}}$ have the same Lagrangian algebras.

The net $\cal A_{V'}$ is a covariant subnet of $\cal A_{H_{\Lambda_a}}$ and of $\cal A_{H_{\Lambda_{c_0}}}$ by Theorem \ref{thm:strong-locality}. Let $y \in \{a,c_0\}$. The inclusion $\cal A_{V'} \subset \cal A_{H_{\Lambda_y}}$ is irreducible by Theorem \ref{thm:strong-locality}, and of finite index by part 2 of Proposition \ref{prop:lattice-nets}. By part 3 of Proposition \ref{prop:nets}, its dual canonical endomorphism $\Theta'_y$, which is the restriction of the vacuum representation of $\cal A_{H_{\Lambda_y}}$ to $\cal A_{V'}$, is a connected \'etale algebra in $\Rep(\cal A_{V'})$. By part 3 of Proposition \ref{prop:nets} and part 2 of Proposition \ref{prop:lattice-nets}, we have $d(\Theta'_y) = [\cal A_{H_{\Lambda_y}} : \cal A_{V'}] = 3$ and $\mu_{\cal A_{V'}} = 9 = d(\Theta'_y)^2$, so that $\Theta'_y$ is Lagrangian by \eqref{eq:mu-index}. The restriction of $\Theta'_y$ to $\cal A_U$ is the restriction of the vacuum representation of $\cal A_{H_{\Lambda_y}}$, that is, $\pi_0 \oplus \pi_b \oplus 2\pi_y$ by part 2 of Theorem \ref{thm:complete-rationality}. Since $\pi_a$ and $\pi_{c_0}$ are inequivalent, the algebras $\Theta'_a$ and $\Theta'_{c_0}$ are not isomorphic. Hence $\Psi_x(\Theta'_a)$ and $\Psi_x(\Theta'_{c_0})$ are two non-isomorphic Lagrangian algebras in $Z(\cal P_x)$. By part 2 of Proposition \ref{prop:DMNO}, which gives a bijection between the isomorphism classes of Lagrangian algebras of $Z(\cal P_x)$ and the equivalence classes of indecomposable $\cal P_x$-module categories, the category $\cal P_x$ has two inequivalent indecomposable module categories. We conclude with part 2 of Lemma \ref{lem:pointed-trivial}.
\end{proof}

\begin{thm}
\label{thm:Fx-haagerup}
Let $U \subset V'$ be as in Theorem \ref{thm:complete-rationality}, and let $x \in \{a,c_0\}$. Then the fusion category $\cal F_x$ is equivalent to $\cal H_2$ or to $\cal H_3$, and $\delta = d^2$. Consequently, $\Rep(\cal A_U)$ is equivalent to $Z(\cal H_3)$ as a unitary modular tensor category. In particular, $\mathrm{(H_J)}$ implies $\mathrm{(N_0)}$.
\end{thm}

\begin{proof}
By part 4 of Theorem \ref{thm:complete-rationality} and Proposition \ref{prop:delta-big}, the category $\cal F_x$ is a unitary fusion category with six simple objects, three of which are invertible. By part 2 of Theorem \ref{thm:complete-rationality}, we have $\Theta_x \cong \pi_0 \oplus \pi_b \oplus 2\pi_x$ with pairwise inequivalent irreducible summands, so that $\operatorname{Hom}(\Theta_x,\Theta_x) \cong \Cbb \oplus \Cbb \oplus \M_2(\Cbb)$. Since $\Theta_x$ corresponds to $I(1)$ (Lemma \ref{lem:rank}), and since $K(\cal F_x) \ot \Cbb$ is isomorphic to $\operatorname{End}(I(1))$ for the spherical category $\cal F_x$ \cite[Corollary 2.16 p.~7]{Ost15}, see also \cite[Corollary 3.4 p.~16]{BiJ22}, the Grothendieck ring of $\cal F_x$ is noncommutative. The fusion subcategory $\cal P_x$ is pointed of dimension three, so that its three simple objects are invertible, and $\cal F_x$ has three invertible simple objects by Proposition \ref{prop:delta-big}. Hence $\cal P_x$ is the fusion subcategory generated by the invertible objects, and it has trivial associator by Lemma \ref{lem:trivial-associator}. By Theorem \ref{thm:recognition}, $\cal F_x$ is equivalent to $\cal H_2$ or to $\cal H_3$. Comparing $\operatorname{FPdim}(\cal F_x) = 3(1+\delta)$ with \eqref{eq:global-dimension}, we obtain $\delta = d^2$. By part 4 of Theorem \ref{thm:complete-rationality} and Section \ref{subsec:haagerup-categories}, $\Rep(\cal A_U)$ is braided equivalent to $Z(\cal F_x)$, hence to $Z(\cal H_3)$, and such an equivalence is naturally isomorphic to a unitary braided equivalence \cite[Theorem 5 p.~2]{Reu23}. Under $\mathrm{(H_J)}$, these results apply to $U = S$ by Corollary \ref{cor:S-rational}.
\end{proof}

\begin{remark}\normalfont
\label{rem:huang-recognition}
For the algebra of \cite{Hua26}, the fusion rules of part 1 of Theorem \ref{thm:recognition}, the triviality of the associator of the pointed part and the identification with $\cal H_2$ or $\cal H_3$ are stated in \cite[Sections 5 and 6 pp.~7--9]{Hua26}, with an outline of proof, and the proofs of this section follow that outline. The outline leaves the following points implicit or only sketched.
\begin{enumerate}
\item The free actions of the invertible objects (Lemma \ref{lem:ring}), which are stated in \cite[p.~7]{Hua26}, and the self-duality of the non-invertible simple objects, which is implicit in \cite[(24) p.~7]{Hua26}.
\item The structure of $Z(\cal F)$ (Lemma \ref{lem:recognition-center}) and the equalities $k_1 = k_2 = k_3 = m$, which are obtained from $\nu_1(\alpha)$ and $\nu_3(\alpha)$ (Lemma \ref{lem:moments}). The identities \cite[(25) p.~7]{Hua26}, attributed there to the first two Frobenius--Schur indicators, follow from them.
\item The restrictions on the squarefree part $D$ of $9m^2+4$ in the proof of Proposition \ref{prop:m-equals-one}. For $D \neq 10$, the second inequality displayed in \cite[p.~7]{Hua26} amounts to $\varphi(2D) \geq \frac{12\sqrt D}{\sqrt{13}}$. It holds for the values of $D$ that can occur, but not for every $D \neq 10$, for instance not for $D \in \{2,5,7\}$. The value $D = 7$ is excluded since the odd prime factors of $D$ are congruent to $1$ modulo $4$, and the values $D \in \{2,5\}$ are excluded by the congruence $D \equiv 1 \bmod 3$.
\item The bound of Lemma \ref{lem:galois-average} for $D = 10$, which is asserted in \cite[p.~7]{Hua26}, and the analysis of the case $m = 2$, which is sketched in \cite[pp.~7--8]{Hua26}.
\item The triviality of the associator of the pointed part. It is derived in \cite[p.~7]{Hua26} from the transversality of the two isotropic lines, by induction from the second one, and it is used in \cite[pp.~8--9]{Hua26}. Lemma \ref{lem:trivial-associator} uses the same input, the second lattice extension, through the module categories of $\cal P_x$.
\item The converse statement of Theorem \ref{thm:recognition}, which is needed to deduce from \cite[Theorems 5.9 and 9.1 pp.~36 and 57]{Izu18} that the category is $\cal H_2$ or $\cal H_3$.
\end{enumerate}
The use of \cite{Hua26} in Section \ref{sec:rationality} is different, see Remark \ref{rem:delta}.
\end{remark}
%%%%%%%%%%%%%%%%%%%%%%%%%%%%%%%%%%%%%%%%%%%%%%%%%%%%%%%%%%%%%%%%%%%%%%%%%%%%%
\section{Consequences for the net of \texorpdfstring{$S$}{S}, complete unitarity and \texorpdfstring{$\mathrm{(N)}$}{(N)}}
\label{sec:N}

The proofs in this section take $\mathrm{(H_J)}$ as a hypothesis. For the explicit $S$ it is supplied by \cite[Theorem 6.1]{ArhIII}, and there is no circularity. The proof of \cite[Theorem 6.1]{ArhIII} uses no statement of the present article which assumes $\mathrm{(H_J)}$. It uses the results of Sections \ref{sec:preliminaries}, \ref{sec:lagrangian} and \ref{sec:recognition} which do not concern Haagerup algebras, Section \ref{sec:net} and Proposition \ref{prop:fusion-unitary}, and it repeats, with \textup{(I1)} and \textup{(F1)}--\textup{(F3)} of \cite[Theorems 2.1 and 2.2]{ArhIII} in place of $\mathrm{(H_J)}$, the arguments of Sections \ref{subsec:rationality-U} and \ref{subsec:recognition-nets}, of the proof of Lemma \ref{lem:N0}, of the categorical part of the proof of Proposition \ref{prop:diagonal} and of the proof of Theorem \ref{thm:S-cu}. The proofs of Lemma \ref{lem:N0}, of parts 1--3 of Theorem \ref{thm:consequences}, of Corollaries \ref{cor:bischoff} and \ref{cor:orbifold}, of Proposition \ref{prop:hypergroup} and of Remark \ref{rem:near-group} use $\mathrm{(H_J)}$ only through Corollary \ref{cor:S-rational} and the assertion $\mathrm{(N_0)}$ below, except for the last assertion of Lemma \ref{lem:N0}, which involves $\Rep(S)$. We use the notation of Section \ref{sec:rationality} for $U = S$. In particular, $\cal A_S$ is completely rational by Corollary \ref{cor:S-rational}. We consider the two assertions
\begin{itemize}
\item[$\mathrm{(N_0)}$] the unitary modular tensor category $\Rep(\cal A_S)$ is braided equivalent to $Z(\cal H_3)$,
\item[$\mathrm{(N)}$] there exists a braided equivalence $\Phi \co \Rep(\cal A_S) \to \Rep(S)$ with $\Phi(\pi_0) \cong S$ and $\Phi(\pi_M) \cong M$ for $M \in \{M_b,M_a,M_{c_0}\}$.
\end{itemize}
Under $\mathrm{(H_J)}$, the category $\Rep(S)$ is braided equivalent to $Z(\cal H_3)$, so that $\mathrm{(N)}$ implies $\mathrm{(N_0)}$. By Theorem \ref{thm:Fx-haagerup}, $\mathrm{(N_0)}$ holds under $\mathrm{(H_J)}$. We prove $\mathrm{(N)}$ under $\mathrm{(H_J)}$ in Theorem \ref{thm:S-cu}, at the end of this section. Its proof uses Theorem \ref{thm:strong-locality}, Lemma \ref{lem:submodules}, Propositions \ref{prop:diagonal} and \ref{prop:diagonal-unitary} and the results of Gui \cite{Gui26,Gui25}, but neither Theorem \ref{thm:Fx-haagerup} nor any statement of this section which assumes $\mathrm{(N)}$ or the complete unitarity of $S$. Since $\mathrm{(N)}$ implies $\mathrm{(N_0)}$ under $\mathrm{(H_J)}$, Theorem \ref{thm:S-cu} gives a second, computer-assisted, proof of $\mathrm{(N_0)}$, independent of Section \ref{sec:recognition}.

%%%%%%%%%%%%%%%%%%%%%%%%%%%%%%%%%%%%%%%%%%%%%%%%%%%%%%%%%%%%%%%%%%%%%%%%%%%%%
\subsection{The representation category and the holomorphic extensions}

\begin{lemma}
\label{lem:N0}
Assume $\mathrm{(H_J)}$, and let $\Psi \co \Rep(\cal A_S) \to Z(\cal H_3)$ be a braided equivalence, which exists by Theorem \ref{thm:Fx-haagerup}. Then $\Psi(\pi_b) \cong b$, $\{\Psi(\pi_a),\Psi(\pi_{c_0})\} = \{a,c_0\}$ up to isomorphism, and $\delta = d^2 = 3d+1$. Moreover, $\Psi$ is unique up to natural isomorphism and naturally isomorphic to a unitary braided equivalence. Consequently, $\mathrm{(N)}$ holds if and only if $\Psi(\pi_a) \cong a$.
\end{lemma}

\begin{proof}
The uniqueness follows from Section \ref{subsec:haagerup-categories}, and the unitarity from \cite[Theorem 5 p.~2]{Reu23}. In the unitary category $\Rep(\cal A_S)$, the Frobenius--Perron dimension of an irreducible sector is its statistical dimension, and $\operatorname{FPdim}\Rep(\cal A_S) = \mu_{\cal A_S}$ by \eqref{eq:mu-index}. By the existence of $\Psi$, \eqref{eq:dim-center} and part 3 of Theorem \ref{thm:complete-rationality}, we obtain $9(1+\delta)^2 = 9(3d+2)^2$, that is, $\delta = 3d+1 = d^2$. The simple object $\Psi(\pi_b)$ has dimension $3d+1$, so that $\Psi(\pi_b) \cong b$ by \eqref{eq:center-dimensions}. For $x \in \{a,c_0\}$, the object $\Psi(\pi_x)$ is a simple summand of dimension $3d+2$ of the Lagrangian algebra $\Psi(\Theta_x)$, hence it belongs to $\{a,c_0\}$ by Lemma \ref{lem:possible-objects}. The objects $\Psi(\pi_a)$ and $\Psi(\pi_{c_0})$ are not isomorphic since $\pi_a$ and $\pi_{c_0}$ are inequivalent by Theorem \ref{thm:complete-rationality}. Finally, under $\mathrm{(H_J)}$, the braided equivalence $\Rep(S) \cong Z(\cal H_3)$ is unique up to natural isomorphism, and it maps $M_z$ to $z$. A braided equivalence $\Phi \co \Rep(\cal A_S) \to \Rep(S)$ is therefore isomorphic to the composition of $\Psi$ with the inverse of this equivalence, and $\Phi(\pi_b) \cong M_b$ by the first part. Hence $\mathrm{(N)}$ amounts to $\Psi(\pi_a) \cong a$.
\end{proof}

\begin{thm}
\label{thm:consequences}
Assume $\mathrm{(H_J)}$.
\begin{enumerate}
\item $\Rep(\cal A_S)$ is equivalent to $Z(\cal H_3)$ as a unitary modular tensor category. We have $d(\pi_b) = 3d+1$, $d(\pi_a) = d(\pi_{c_0}) = 3d+2$, $[\cal A_{V'} : \cal A_S] = 3d+2$, $[\cal A_{E_8} : \cal A_S] = 3(3d+2)$ for both lattice extensions, and $\mu_{\cal A_S} = 9(3d+2)^2$.
\item Up to isomorphism of extensions, $\cal A_S$ has exactly three irreducible local extensions $\cal B$ with $\mu_{\cal B} = 1$. Their dual canonical endomorphisms $\Theta_1$, $\Theta_2$, $\Theta_3$ satisfy $\Rep(\cal A_S)_{\Theta_i} \cong \cal H_i$, $\Theta_1 \cong \pi_0 \oplus \pi_b \oplus \pi_a \oplus \pi_{c_0}$ and
\[
\{\Theta_2,\Theta_3\} = \{\Theta_a,\Theta_{c_0}\},\quad \Theta_x \cong \pi_0 \oplus \pi_b \oplus 2\pi_x.
\]
The extensions attached to $\Theta_a$ and $\Theta_{c_0}$ are the nets $\cal A_{H_{\Lambda_a}}$ and $\cal A_{H_{\Lambda_{c_0}}}$ of the lattice extensions, both isomorphic to $\cal A_{E_8}$.
\item The only intermediate nets between $\cal A_S$ and the extension $\cal B_1$ attached to $\Theta_1$ are $\cal A_S$ and $\cal B_1$. For $x \in \{a,c_0\}$, the only intermediate nets between $\cal A_S$ and $\cal A_{H_{\Lambda_x}}$ are $\cal A_S$, $\cal A_{V'}$ and $\cal A_{H_{\Lambda_x}}$.
\item Assume $\mathrm{(N)}$. Then $\Theta_2 = \Theta_a$ and $\Theta_3 = \Theta_{c_0}$, that is, the net $\cal A_{H_{\Lambda_a}}$ carries $\cal H_2$ and the net $\cal A_{H_{\Lambda_{c_0}}}$ carries $\cal H_3$, as the vertex operator algebras $H_{\Lambda_a}$ and $H_{\Lambda_{c_0}}$ do.
\end{enumerate}
\end{thm}

\begin{proof}
1. This follows from Lemma \ref{lem:N0} and part 3 of Theorem \ref{thm:complete-rationality}.

2. Let $\cal B$ be an irreducible local extension of $\cal A_S$, that is, $\cal A_S(I)' \cap \cal B(I) = \Cbb$. Then $\cal A_S(S^1)' \cap \cal B(I) = \Cbb$, so that $[\cal B : \cal A_S] < \infty$ by part 2 of Proposition \ref{prop:nets}, and $\cal B$ is completely rational with $\mu_{\cal B} = \mu_{\cal A_S}/d(\Theta_{\cal B})^2$ by part 1 of the same proposition, where $\Theta_{\cal B}$ is the dual canonical endomorphism. Hence $\mu_{\cal B} = 1$ if and only if the connected \'etale algebra $\Theta_{\cal B}$ is Lagrangian in $\Rep(\cal A_S)$. Irreducible local extensions of $\cal A_S$ up to isomorphism correspond to commutative $\Q$-systems $\Theta$ with $\dim\operatorname{Hom}(\mathrm{id},\Theta) = 1$ up to unitary equivalence \cite[Theorem 4.9 p.~25]{LoR95}, \cite[Section 5.2.1 p.~77]{BKLR15}, see also \cite[Lemma 15 p.~14]{Lon03}. By \cite[Theorem 3.2 p.~11 and Theorem 3.9 p.~15]{CGGH23}, the latter correspond to connected \'etale algebras up to isomorphism. By Lemma \ref{lem:N0} and Propositions \ref{prop:three-lagrangian} and \ref{prop:which-line}, there are exactly three Lagrangian algebras in $\Rep(\cal A_S)$, with the stated categories of modules and objects. Their objects are pairwise non-isomorphic. By part 2 of Theorem \ref{thm:complete-rationality}, $\Theta_x \cong \pi_0 \oplus \pi_b \oplus 2\pi_x$ is the dual canonical endomorphism of $\cal A_S \subset \cal A_{H_{\Lambda_x}}$, and by part 4 of the same theorem it is Lagrangian. The nets $\cal A_{H_{\Lambda_x}}$ are isomorphic to $\cal A_{E_8}$ by Theorem \ref{thm:strong-locality}.

3. By \cite[Proposition 4.11 p.~49]{BKLR15}, an intermediate net gives, at a fixed interval, a projection in $\operatorname{Hom}(\Theta_i,\Theta_i)$ whose range is a connected \'etale subalgebra of $\Theta_i$, and distinct intermediate nets give distinct subalgebras. As in the proof of part 3 of Theorem \ref{thm:lagrangian-extensions}, using Lemma \ref{lem:rank}, the connected \'etale subalgebras of $\Theta_i$ correspond to the fusion subcategories of $\cal H_i$, which are given by part 4 of Proposition \ref{prop:GS}. Equivalently, one can use the Galois correspondence of part 1 of Proposition \ref{prop:bischoff}. For $\Theta_x$, the subalgebra of Frobenius--Perron dimension $3d+2$ is the dual canonical endomorphism $\pi_0 \oplus \pi_b$ of the intermediate net $\cal A_{V'}$.

4. By Lemma \ref{lem:N0}, the braided equivalence $\Psi$ maps $\pi_a$ to $a$ and $\pi_{c_0}$ to $c_0$, so that it maps $\Theta_a$ to a Lagrangian algebra isomorphic, as an object, to $0 \oplus b \oplus 2a$, that is, to $L_2$ by Proposition \ref{prop:which-line}. Similarly $\Psi(\Theta_{c_0}) \cong L_3$.
\end{proof}

\begin{remark}\normalfont
Part 2 alone shows that the two lattice nets carry the categories $\cal H_2$ and $\cal H_3$ in some order. Under $\mathrm{(N)}$, the order is the same as for the vertex operator algebras, by Proposition \ref{prop:which-line}, and $\mathrm{(N)}$ holds under $\mathrm{(H_J)}$ by Theorem \ref{thm:S-cu}, see Corollary \ref{cor:N-consequences}.
\end{remark}

%%%%%%%%%%%%%%%%%%%%%%%%%%%%%%%%%%%%%%%%%%%%%%%%%%%%%%%%%%%%%%%%%%%%%%%%%%%%%
\subsection{The Haagerup subfactor and the generalized orbifold}

\begin{cor}
\label{cor:bischoff}
Assume $\mathrm{(H_J)}$. Then the Haagerup subfactor arises from the completely rational net $\cal A_S$ in the sense of \cite[Definition 4.1 p.~20]{Bis17a}, that is, the hyperfinite type $\mathrm{III}_1$ subfactor with the standard invariant of the Haagerup subfactor is of the form described there.
\end{cor}

\begin{proof}
The Haagerup subfactor $N \subset M$ has finite index and finite depth, and its quantum double $D(N \subset M)$ is the Drinfeld center of either of its even parts $\cal H_1$ and $\cal H_2$, which is braided equivalent to $Z(\cal H_3)$. By part 1 of Theorem \ref{thm:consequences}, $\Rep(\cal A_S)$ is braided equivalent to $D(N \subset M)$, and the conclusion follows from part 2 of Proposition \ref{prop:bischoff}.
\end{proof}

\begin{prop}
\label{prop:hypergroup}
Assume $\mathrm{(H_J)}$. Then there is a proper action of the hypergroup $K = \{k_0,k_1\}$ with
\[
k_1k_1 = \frac{1}{\delta}k_0+\Big(1-\frac{1}{\delta}\Big)k_1
\]
on the net $\cal A_{V'}$ such that $\cal A_S = \cal A_{V'}^K$. The weight of $K$ is $1+\delta = [\cal A_{V'} : \cal A_S]$.
\end{prop}

\begin{proof}
This is Proposition \ref{prop:hypergroup-U} for $U = S$, which applies by Corollary \ref{cor:S-rational}.
\end{proof}

\begin{cor}
\label{cor:orbifold}
Assume $\mathrm{(H_J)}$. Then the hypergroup of Proposition \ref{prop:hypergroup} has parameter $w = \delta = d^2 = \frac{11+3\sqrt{13}}{2}$. With $\tilde k_1 \ov{\mathrm{def}}{=} d k_1 = w^{1/2}k_1$, we have $\tilde k_1\tilde k_1 = k_0+3\tilde k_1$, and $\Rep(\cal A_{V'}^K) = \Rep(\cal A_S)$ is braided equivalent to the quantum double of the Haagerup subfactor. The net $\cal A_{V'}$ is the conformal net $\cal A_{E_6,1} \ot \cal A_{SU(3),1}$, denoted $E_{6,1} \times SU(3)_1$ in \cite[Table 1 p.~25]{Bis17b}, of the lattice $A_2 \times E_6$. Hence, under these hypotheses, $\cal A_S$ solves \cite[Problem 4.27 p.~27]{Bis17b}.
\end{cor}

\begin{proof}
By Lemma \ref{lem:N0}, $\delta = d^2 = 3d+1 = \frac{11+3\sqrt{13}}{2}$. Then
\[
\tilde k_1\tilde k_1 = d^2k_1k_1 = k_0+(d^2-1)k_1 = k_0+\frac{d^2-1}{d}\tilde k_1 = k_0+3\tilde k_1,
\]
since $d^2-1 = 3d$. The other assertions follow from Proposition \ref{prop:hypergroup}, Corollary \ref{cor:bischoff}, part 1 of Theorem \ref{thm:consequences} and part 2 of Proposition \ref{prop:lattice-nets}.
\end{proof}

\begin{remark}\normalfont
The vertex algebraic counterpart of Problem 4.27 of \cite{Bis17b} is formulated in \cite[Claim 4.4 and (4.54) p.~86]{GaR26}. The existence of the hypergroup action is a general fact \cite[Theorem 4.22 p.~24]{Bis17b}. The new input of Corollary \ref{cor:orbifold} is the identification of $\Rep(\cal A_S)$ in Theorem \ref{thm:Fx-haagerup}. The value of $w$ is consistent with \cite[Corollary 4.20 p.~23 and Lemma 4.26 p.~27]{Bis17b}. Indeed, by Theorem \ref{thm:Fx-haagerup}, the canonical hypergroup of the inclusion $\cal A_S \subset \cal A_{E_8}$ is attached to a fusion category with the fusion rules \eqref{eq:H2-rules} \cite[Theorem 5.16 p.~31]{Bis17b}, the intermediate net $\cal A_{V'}$ is the fixed point net of its maximal subgroup $\Z_3$, and \cite[Lemma 4.26 p.~27]{Bis17b} with $n = 3$ gives $w = \frac{2+9+3\sqrt{13}}{2}$.
\end{remark}

\begin{remark}\normalfont
\label{rem:near-group}
Bischoff deduced from the modular invariant of \cite[(2.4) p.~20]{EvG11} that the unitary fusion category attached to the inclusion $\cal A_S \subset \cal A_{V'}$ should be a near-group category of type $\Z_3^2+9$ \cite[Remark 5.14 and Proposition 5.15 p.~31]{Bis17b}. Under $\mathrm{(H_J)}$, this holds for the category $\Rep(\cal A_S)_\theta$ of modules over the dual canonical endomorphism $\theta = \pi_0 \oplus \pi_b$. Indeed, by part 1 of Theorem \ref{thm:consequences}, we may compute in $Z(\cal H_3)$ with the algebra $A = 0 \oplus b$. For simple objects $z$ and $z'$, the free module adjunction gives $\operatorname{Hom}_A(z \ot A,z' \ot A) \cong \operatorname{Hom}(z,z' \ot A)$, and the fusion rules of $Z(\cal H_3)$, computed from the $S$-matrix of Section \ref{subsec:center} by the Verlinde formula,
\[
b \ot b \cong \bigoplus_z z,\quad y \ot b \cong y \oplus \bigoplus_{z \neq 0} z,\quad d_l \ot b \cong \bigoplus_{z \notin \{0,d_l\}} z,
\]
for $y \in \{a,c_0,c_1,c_2\}$ and $1 \leq l \leq 6$, give the following. The free modules $d_l \ot A$ are simple and pairwise isomorphic, to an object $\rho$. The module $b \ot A$ is isomorphic to $A \oplus \rho$. For $y \in \{a,c_0,c_1,c_2\}$, the module $y \ot A$ is isomorphic to $\rho \oplus g_y \oplus g'_y$, where the eight simple modules $g_y$, $g'_y$ are pairwise non-isomorphic and not isomorphic to $A$ or $\rho$. Since every simple $A$-module is a summand of a free module, the category $\Rep(\cal A_S)_\theta$ has ten simple objects. Since the free module functor $z \mapsto z \ot A$ is a tensor functor, it preserves Frobenius--Perron dimensions. Hence $\rho$ has dimension $3d$ and $\operatorname{FPdim}(g_y)+\operatorname{FPdim}(g'_y) = 2$. Since the Frobenius--Perron dimension of a simple object is at least one \cite[Proposition 3.3.4 p.~53]{EGNO15}, and a simple object of Frobenius--Perron dimension one is invertible \cite[Corollary 3.3.10 p.~54]{EGNO15}, the nine objects $A$, $g_y$, $g'_y$ are invertible. They form a group of order nine. By \eqref{eq:dims-local}, the category $\Rep(\cal A_S)^0_\theta$ of local modules has Frobenius--Perron dimension nine. Since $\operatorname{FPdim}(\rho)^2 = 9d^2 > 9$, it cannot contain $\rho$, so its simple objects are these nine invertible objects. By \cite[Proposition 6.4 p.~33]{BKL15} and the proof of Lemma \ref{lem:trivial-associator}, it is braided equivalent to $\Rep(\cal A_{V'})$ and to $Z(\cal P_x) \cong Z(\mathrm{Vec}_{\Z_3})$, whose group of invertible objects is $\Z_3 \times \Z_3$. For an invertible object $g$, the object $\rho \ot g$ is simple of dimension $3d$, so that $\rho \ot g \cong \rho$, and similarly $\rho^* \cong \rho$. Hence $\operatorname{Hom}(\rho \ot \rho,g) \cong \operatorname{Hom}(\rho,g \ot \rho)$ is one-dimensional, and comparing dimensions with $9d^2 = 9+9 \cdot 3d$ gives $\rho \ot \rho \cong \bigoplus_g g \oplus 9\rho$. Thus $\Rep(\cal A_S)_\theta$ is a near-group category of type $\Z_3^2+9$, of Frobenius--Perron dimension $9(3d+2)$.
\end{remark}

\begin{remark}\normalfont
By \cite[Theorem 4.22 p.~24]{Bis17b} and the computation of \cite[Section 4.6 p.~25]{Bis17b}, the action of $k_1$ is implemented on $\cal H_{V'} = \cal H_S \oplus \cal H_{N_b}$ by the operator $T = P_S-\delta^{-1}P_{N_b}$, where $P_S$ and $P_{N_b}$ are the orthogonal projections. In other words, the stochastic map $\phi_{k_1}$ on $\cal A_{V'}(I)$ is determined by $\phi_{k_1}(x)\Omega = Tx\Omega$. The restriction of $T$ to $V' = S \oplus N_b$ is the linear map $\mathrm{id}_S \oplus (-\delta^{-1})\mathrm{id}_{N_b}$, which is not an automorphism of $V'$. Compare with the hypergroup symmetries of \cite[Section 4.3 pp.~83--88]{GaR26}.
\end{remark}

%%%%%%%%%%%%%%%%%%%%%%%%%%%%%%%%%%%%%%%%%%%%%%%%%%%%%%%%%%%%%%%%%%%%%%%%%%%%%
\subsection{The third holomorphic extension}

We say that $S$ is completely unitary if it is completely unitary in the sense of \cite[Definition 4.2 p.~18]{CGGH23}, that is, every $S$-module is unitarizable and the invariant sesquilinear forms on the tensor products of $S$-modules of \cite{Gui19} are positive definite.

\begin{prop}
\label{prop:third-net}
Assume $\mathrm{(H_J)}$ and that $S$ is completely unitary. Let $B_1$ be the holomorphic extension of $S$ attached to $L_1$ by Theorem \ref{thm:lagrangian-extensions}.
\begin{enumerate}
\item The vertex operator algebra $B_1$ has a unique unitary structure extending the unitary structure of $S$ given by $H_\chi$. With this structure, it is strongly local and its conformal net $\cal A_{B_1}$ is isomorphic to $\cal A_{E_8}$.
\item The net $\cal A_{B_1}$ is an irreducible local extension of $\cal A_S$, and its vacuum representation restricts to $\cal A_S$ as $\pi_0 \oplus \pi_b \oplus \pi_a \oplus \pi_{c_0}$. In particular, $\cal A_{H_{\Lambda_a}}$, $\cal A_{H_{\Lambda_{c_0}}}$ and $\cal A_{B_1}$ are three pairwise non-isomorphic extensions of $\cal A_S$ with $\mu$-index one, each isomorphic to $\cal A_{E_8}$ as a conformal net.
\item The net $\cal A_{B_1}$ is the extension $\cal B_1$ of part 2 of Theorem \ref{thm:consequences}.
\end{enumerate}
\end{prop}

\begin{proof}
1. By \cite[Theorem 4.7 p.~19]{CGGH23}, $B_1$ is a unitary extension of $S$, with a unique unitary structure extending the given one on $S$. By part 1 of Theorem \ref{thm:lagrangian-extensions}, $B_1 \cong V_{E_8}$. The argument of part 1 of Proposition \ref{prop:lattice-nets} applies to any unitary structure on a vertex operator algebra isomorphic to $V_{E_8}$, so that $B_1$ is strongly local, has bounded localized vertex operators and $\cal A_{B_1} \cong \cal A_{E_8}$.

2. By Lemma \ref{lem:trivial-coset}, the coset of $\cal A_S$ in $\cal A_{B_1}$ is trivial and $\cal A_S(I)' \cap \cal A_{B_1}(I) = \Cbb$. As in the proof of Theorem \ref{thm:complete-rationality}, the index is finite. Since $L_1 \cong 0 \oplus b \oplus a \oplus c_0$ is multiplicity free, the isotypic components of the $S$-module $B_1$ are simple modules isomorphic to $S$, $M_b$, $M_a$, $M_{c_0}$. They are pairwise orthogonal, since the orthogonal projection onto the closure of a simple submodule maps $B_1$ into this closure and commutes with the modes of $S$, as in the proof of part 1 of Lemma \ref{lem:submodules}. By part 3 of Lemma \ref{lem:submodules}, applied to $V_1 = B_1$ and $V_2 \in \{H_{\Lambda_a},H_{\Lambda_{c_0}}\}$, the corresponding representations are $\pi_0$, $\pi_b$, $\pi_a$, $\pi_{c_0}$. The three extensions are pairwise non-isomorphic since their dual canonical endomorphisms are not equivalent.

3. By part 2 of Theorem \ref{thm:consequences}, $\Theta_1$ is the unique Lagrangian $\Q$-system isomorphic to $\pi_0 \oplus \pi_b \oplus \pi_a \oplus \pi_{c_0}$.
\end{proof}

\begin{remark}\normalfont
Under complete unitarity of $S$, the comparison of \cite[Theorem 4.10 p.~21]{CGGH23} between vertex algebra extensions of $S$ and net extensions of $\cal A_S$ requires an equivalence $\Rep^u(S) \cong \Rep^f(\cal A_S)$ of unitary modular tensor categories. If $S$ is completely unitary, such an equivalence follows from $\mathrm{(N_0)}$, $\mathrm{(H_J)}$ and \cite[Theorem 5 p.~2]{Reu23}, and it does not involve the correspondence $\pi_M \mapsto M$ of $\mathrm{(N)}$. Proposition \ref{prop:third-net} only uses the unitarity of the extension $B_1$ and Lemma \ref{lem:submodules}. Both the complete unitarity of $S$ and this equivalence hold under $\mathrm{(H_J)}$ (Theorem \ref{thm:S-cu} and Corollary \ref{cor:N-consequences}).
\end{remark}

%%%%%%%%%%%%%%%%%%%%%%%%%%%%%%%%%%%%%%%%%%%%%%%%%%%%%%%%%%%%%%%%%%%%%%%%%%%%%
\subsection{A diagonal holomorphic extension of \texorpdfstring{$S \ot S$}{S x S}}
\label{subsec:diagonal}

The following statement only uses $\mathrm{(H_J)}$. It realizes $S$ as its own commutant in $V_{E_8 \oplus E_8}$.

\begin{prop}
\label{prop:diagonal}
Assume $\mathrm{(H_J)}$.
\begin{enumerate}
\item There is a braided equivalence $F \co \Rep(S)^{\mathrm{rev}} \to \Rep(S)$, and it is unique up to natural isomorphism. The induced involution $\pi$ of the set of simple objects satisfies $\pi(0) = 0$, $\pi(b) = b$, $\pi(a) = a$, $\pi(c_0) = c_0$, $\pi(c_1) = c_2$, $\pi(d_1) = d_5$, $\pi(d_2) = d_3$ and $\pi(d_4) = d_6$.
\item The $S \ot S$-module $W \ov{\mathrm{def}}{=} \bigoplus_z M_z \ot M_{\pi(z)}$ carries the structure of a simple vertex operator algebra of CFT type extending $S \ot S$, with the same conformal vector. It is holomorphic, it is isomorphic to $V_{E_8 \oplus E_8}$, and $\operatorname{Com}_W(S \ot 1) = 1 \ot S$, $\operatorname{Com}_W(1 \ot S) = S \ot 1$.
\end{enumerate}
\end{prop}

\begin{proof}
1. For a fusion category $\cal C$, let $\cal C^{\ot\mathrm{op}}$ be the category $\cal C$ with the reversed tensor product. A half-braiding $\gamma$ of an object $X$ in $\cal C^{\ot\mathrm{op}}$ gives the half-braiding $\gamma^{-1}$ of $X$ in $\cal C$, and the functor $(X,\gamma) \mapsto (X,\gamma^{-1})$, with the monoidal structure given by the braiding of $Z(\cal C)$, is a braided equivalence from $Z(\cal C^{\ot\mathrm{op}})$ onto $Z(\cal C)^{\mathrm{rev}}$. If $\cal C$ is a unitary fusion category, the functor $f \mapsto f^*$ is a linear monoidal equivalence from the complex conjugate category $\ovl{\cal C}$, which has the same objects and the complex conjugate spaces of morphisms, to the opposite category $\cal C^{\mathrm{op}}$, and the duality functor $X \mapsto X^\vee$ is a monoidal equivalence from $\cal C^{\mathrm{op}}$ to $\cal C^{\ot\mathrm{op}}$. Hence $\cal C^{\ot\mathrm{op}} \cong \ovl{\cal C}$ and $Z(\cal C)^{\mathrm{rev}} \cong Z(\ovl{\cal C})$. Apply this to the even parts $\cal H_1$ and $\cal H_2$ of the Haagerup subfactor. The categories $\ovl{\cal H_1}$ and $\ovl{\cal H_2}$ are the even parts of the complex conjugate subfactor, whose planar algebra is the complex conjugate of the Haagerup planar algebra, with the same index and the same principal graphs. By \cite[Theorem 4.1 p.~21]{JMS14}, it is the Haagerup planar algebra or its dual, whose even parts are $\cal H_1$ and $\cal H_2$, of ranks four and six. Comparing ranks, we obtain $\ovl{\cal H_1} \cong \cal H_1$ and $\ovl{\cal H_2} \cong \cal H_2$. Hence $Z(\cal H_3)^{\mathrm{rev}} \cong Z(\cal H_1)^{\mathrm{rev}} \cong Z(\ovl{\cal H_1}) \cong Z(\cal H_1) \cong Z(\cal H_3)$ by Section \ref{subsec:haagerup-categories}, and $\mathrm{(H_J)}$ gives $F$. Its uniqueness follows from the triviality of the braided autoequivalences of $Z(\cal H_3)$ (Section \ref{subsec:haagerup-categories}). Since $F$ is also a braided equivalence $\Rep(S) \to \Rep(S)^{\mathrm{rev}}$, the composite $F \circ F$ is a braided autoequivalence of $\Rep(S)$. For the same reason, it is isomorphic to the identity functor, so that $\pi$ is an involution. The equivalence $F$ preserves Frobenius--Perron dimensions, and it maps the canonical twist of $\Rep(S)^{\mathrm{rev}}$, which is $\theta^{-1}$, to the canonical twist $\theta$ of $\Rep(S)$ \cite[Proposition 8.23 p.~34]{ENO05}. By \eqref{eq:center-dimensions} and \eqref{eq:center-twists}, we obtain $\pi(0) = 0$, $\pi(b) = b$, $\pi(\{a,c_0\}) = \{a,c_0\}$ and $\pi(c_1) = c_2$. Moreover $\pi(d_l) = d_{l'}$ with $6l'^2 \equiv -6l^2 \bmod 13$, that is, $l' \equiv \pm 5l \bmod 13$, since $5^2 \equiv -1 \bmod 13$. Finally, under $\Rep(S) \cong Z(\cal H_3)$, the Lagrangian algebra $L_2 \cong 0 \oplus b \oplus 2a$ of Corollary \ref{cor:categorical} is also a connected \'etale algebra of $\Rep(S)^{\mathrm{rev}}$, since an algebra is commutative for a braiding if and only if it is commutative for the reverse braiding. The two tensor structures on the category of modules over a connected \'etale algebra defined by a braiding and by its reverse are opposite to each other \cite[p.~9]{DMNO13}, so that $(\Rep(S)^{\mathrm{rev}})_{L_2}$ is equivalent to $\cal H_2^{\ot\mathrm{op}} \cong \ovl{\cal H_2} \cong \cal H_2$. Hence $F(L_2)$ is a Lagrangian algebra of $\Rep(S)$ whose category of modules is equivalent to $\cal H_2$. By Proposition \ref{prop:three-lagrangian}, $F(L_2) \cong L_2$, so that $0 \oplus b \oplus 2\pi(a) \cong 0 \oplus b \oplus 2a$, and $\pi(a) = a$, $\pi(c_0) = c_0$.

2. Every simple $S$-module is self-dual (Section \ref{subsec:center}), and $F$ inverts the twists. By \cite[Main Theorem 3 pp.~4--5 and Theorem 5.10 p.~51]{CKM22}, applied to the braid-reversed equivalence $F$, the space $W$ carries the structure of a simple conformal vertex algebra extending $S \ot S$. It is $\Z$-graded since $\theta_z\theta_{\pi(z)} = 1$, and of CFT type since $S$ is positive. Under the braided equivalence $\Rep(S \ot S) \cong \Rep(S) \boxtimes \Rep(S)$ \cite[Theorem 5.5 p.~48]{CKM22}, $W$ corresponds to the haploid commutative algebra $\bigoplus_z M_z \boxtimes M_{\pi(z)}$ with trivial twist \cite[Theorem 3.2 p.~10]{HKL15}, which is a connected \'etale algebra by \cite[Theorem 5.12 p.~52 and its proof]{CKM22} and \cite[Proposition 2.7 p.~5]{DMNO13}. Its Frobenius--Perron dimension $\sum_z d(M_z)d(M_{\pi(z)}) = \operatorname{FPdim}\Rep(S)$ is the square root of $\operatorname{FPdim}\Rep(S \ot S)$. It is therefore Lagrangian, so that $\Rep(W) \cong \mathrm{Vec}$ by \cite[Theorem 3.4 p.~16]{HKL15} and \cite[Remark 4.3 p.~16]{DMNO13}, and $W$ is strongly rational by \cite[Corollary 5.13 p.~53]{CKM22}. By \cite[Theorem 2 p.~2]{DoM04b}, $W$ is isomorphic to $V_{E_8 \oplus E_8}$ or to $V_{D_{16}^+}$.
The weight-one space $W_1$ is the sum of $(S \ot S)_1 = (\mathfrak g' \oplus \mathfrak g'') \oplus (\mathfrak g' \oplus \mathfrak g'')$ and of the spaces $R_z \ot R_{\pi(z)}$ with $h_z+h_{\pi(z)} = 1$, that is, $z \in \{c_1,c_2,d_1,d_5,d_4,d_6\}$, where $R_z$ is the lowest weight space of $M_z$. Write $\mathfrak g'_1 = \mathfrak g' \ot 1$ and $\mathfrak g'_2 = 1 \ot \mathfrak g'$. By \cite[Theorem 8.5 and Table 2]{Arh26}, the $\mathfrak g'$-weights occurring in these spaces $R_z$ are at most $2$, so that $W_1$ contains no irreducible $\mathfrak g'_k$-module of dimension $4$, and the subspace of $(\mathfrak g'_1 \oplus \mathfrak g'_2)$-invariant vectors of $W_1$ has dimension $6+2 \cdot 18+2 \cdot 15+2 \cdot 3 = 78$. Assume that $W \cong V_{D_{16}^+}$, so that $W_1 \cong \mathfrak{so}(E)$ with $\dim E = 32$ and $W_1 \cong \Lambda^2 E$. The affine subalgebra generated by $\mathfrak g'_k$ has level $3$, so that the Dynkin index of $\mathfrak g'_k \subset \mathfrak{so}(E)$ is $3$. Since the vector representation of $\mathfrak{so}(E)$ has Dynkin index $2$, the restriction of $E$ to $\mathfrak g'_k$ has Dynkin index $6$. Hence $E|_{\mathfrak g'_k}$ is $6V(1) \oplus 20V(0)$ or $V(2) \oplus 2V(1) \oplus 25V(0)$, where $V(n)$ is the irreducible module of highest weight $n$. In the second case $\Lambda^2 E \supset V(2) \ot V(1) \supset V(3)$, which is excluded. Thus $E|_{\mathfrak g'_1 \oplus \mathfrak g'_2} = pV(1,1) \oplus (6-2p)(V(1,0) \oplus V(0,1)) \oplus (8+4p)V(0,0)$ for some $p \in \{0,1,2,3\}$. The multiplicity spaces of $V(1,1)$ and $V(0,0)$ carry nondegenerate symmetric forms and those of $V(1,0)$ and $V(0,1)$ nondegenerate alternating forms, so that the centralizer of $\mathfrak g'_1 \oplus \mathfrak g'_2$ in $\mathfrak{so}(E)$ is $\mathfrak{so}_p \oplus \mathfrak{so}_{8+4p} \oplus \mathfrak{sp}_{6-2p} \oplus \mathfrak{sp}_{6-2p}$, of dimension $70$, $86$, $127$ or $193$. This contradicts the value $78$, so that $W \cong V_{E_8 \oplus E_8}$.
A vector of $W$ annihilated by the nonnegative modes of $S \ot 1$ is annihilated by $L^S(0) \ot 1$. Since $h_z > 0$ for $z \neq 0$, it lies in $S_0 \ot S = \mathbf 1 \ot S$. This gives the commutants.
\end{proof}

\begin{remark}\normalfont
The vertex operator algebra $W \cong V_{E_8 \oplus E_8}$ is generated by $W_1 \cong \mathfrak e_8 \oplus \mathfrak e_8$, so that its unitary structures are determined by compact real forms of $W_1$, and $S \ot S$ is generated by its weight-one and weight-two spaces. By Corollary \ref{cor:modules-unitary}, every module $M_z \ot M_{\pi(z)}$ is unitarizable. If $S$ is completely unitary, then so is $S \ot S$ \cite[Proposition 3.31 p.~48]{Gui22}, and $W$ is a unitary extension of $S \ot S$ \cite[Theorem 4.7 p.~19]{CGGH23}. Without complete unitarity, it suffices to find the PCT operator of a unitary structure on $W$ which preserves $(S \ot S)_1$ and $\omega_{C \ot 1}$, see the proof of Proposition \ref{prop:diagonal-unitary}. Such an operator is found in Section \ref{subsec:complete-unitarity}.
\end{remark}

%%%%%%%%%%%%%%%%%%%%%%%%%%%%%%%%%%%%%%%%%%%%%%%%%%%%%%%%%%%%%%%%%%%%%%%%%%%%%
\subsection{Complete unitarity of \texorpdfstring{$S$}{S} and the comparison \texorpdfstring{$\mathrm{(N)}$}{(N)}}
\label{subsec:complete-unitarity}

We solve the existence problem of the preceding remark, and we deduce the complete unitarity of $S$ and $\mathrm{(N)}$ from the results of Gui on cosets \cite{Gui26,Gui25}. We assume $\mathrm{(H_J)}$ and keep the notation of Proposition \ref{prop:diagonal}. Let $A$ be the vertex subalgebra of $S$ generated by $S_1 = \mathfrak g' \oplus \mathfrak g''$, which is isomorphic to $L_{\mathfrak{sl}_2}(3,0) \ot L_{\mathfrak{sl}_2}(13,0)$ \cite[Corollary 9.2]{Arh26}, and let $C = \operatorname{Com}_S(A)$. The vector $\omega_C = \omega-\omega_A$, where $\omega_A$ is the Sugawara vector of $A$, is the conformal vector of $C$, of central charge $8-\frac95-\frac{13}{5} = \frac{18}{5}$. We write $\omega_{S \ot 1} = \omega \ot \mathbf 1$, $\omega_{A \ot 1} = \omega_A \ot \mathbf 1$, $\omega_{C \ot 1} = \omega_C \ot \mathbf 1$, and similarly $\omega_{1 \ot S}$, $\omega_{1 \ot A}$, $\omega_{1 \ot C}$, and we denote by $L_n^{S \ot 1}$, $L_n^{A \ot 1}, \ldots$ the modes of the corresponding Virasoro fields. We put $\omega_{A \ot A} = \omega_{A \ot 1}+\omega_{1 \ot A}$ and $\omega_{C \ot C} = \omega_{C \ot 1}+\omega_{1 \ot C}$, so that $\omega_W = \omega_{S \ot 1}+\omega_{1 \ot S} = \omega_{A \ot A}+\omega_{C \ot C}$ is the conformal vector of $W$. We denote by $\Rep^u(S)$ the category of unitary $S$-modules \cite[Section 3.4 p.~44]{Gui26} and $\Rep^f(\cal A_S) = \Rep(\cal A_S)$ is as in Section \ref{subsec:nets}.

Let $\Gamma = E_8 \oplus E_8 \subset \R^8 \oplus \R^8$, where $E_8$ is the lattice of the vectors of $\R^8$ whose coordinates are all in $\Z$ or all in $\frac12+\Z$ and have an even sum, with the simple roots $\alpha_1, \ldots, \alpha_8$ and the fundamental weights $\varpi_1, \ldots, \varpi_8$ of \cite[Planche VII]{Bou68}. Let $\epsi_8$ be a bimultiplicative cocycle on $E_8$ such that $\epsi_8(\alpha,\beta)\epsi_8(\beta,\alpha)^{-1} = (-1)^{\langle\alpha,\beta\rangle}$ and $\epsi_8(\alpha,\alpha) = (-1)^{\langle\alpha,\alpha\rangle/2}$, and let $V_\Gamma = M(1) \ot \Cbb_\epsi[\Gamma]$ be the lattice vertex operator algebra built with the cocycle $\epsi(\alpha,\beta) = \epsi_8(\alpha^{(1)},\beta^{(1)})\epsi_8(\alpha^{(2)},\beta^{(2)})$, where $\alpha^{(1)}$ and $\alpha^{(2)}$ are the components of $\alpha \in \Gamma$. The cocycle $\epsi$ satisfies the same two identities on $\Gamma$, and it is invariant under the exchange of the two summands of $\Gamma$. Let $\phi_1$ be the antilinear involution of $V_\Gamma$ given by the formula of Section \ref{subsec:lattice-extensions}. For a character $\kappa \co \Gamma \to \R_{>0}$, written $\kappa(\lambda) = \e^{\langle h_\kappa,\lambda\rangle}$, we put $\phi_\kappa = \phi_1 \circ \exp(h_\kappa(0))$, so that $\phi_\kappa(\e^\lambda) = \kappa(\lambda)\e^{-\lambda}$. By \cite[Theorem 4.12 p.~24]{DoL14}, $\phi_1$ is the PCT operator of a structure of unitary vertex operator algebra on $V_\Gamma$, as in part 1 of Lemma \ref{lem:unitary-extension}. By Proposition \ref{prop:diagonal}, there is an isomorphism of vertex operator algebras $W \to V_\Gamma$. We fix one, and we identify $W$ with $V_\Gamma$. The space $W_1$ is the Lie algebra $\mathfrak e_8 \oplus \mathfrak e_8$ for the bracket $[x,y] = x_{(0)}y$, with the invariant form defined by $x_{(1)}y = \langle x,y\rangle\mathbf 1$, for which the roots have norm two. The Dynkin index of a homomorphism $\psi \co \mathfrak{sl}_2 \to \mathfrak e_8$ is $\frac12\langle\psi(h),\psi(h)\rangle$, where $h$ is the standard Cartan element of $\mathfrak{sl}_2$. The level of the affine vertex subalgebra of $W$ generated by the image of a homomorphism $\mathfrak{sl}_2 \to \mathfrak e_8 \oplus \mathfrak e_8$ is the sum of the Dynkin indices of its two components.

\paragraph{An explicit subalgebra.} In each copy of $\mathfrak e_8$, the vectors $h(-1)\mathbf 1$, $h \in \Cbb^8$, and $\e^\alpha$, for the roots $\alpha = \pm e_i \pm e_j$ of $D_8 \subset E_8$, span a Lie subalgebra isomorphic to $\mathfrak{so}(16)$. We identify it with the orthogonal Lie algebra of $X = \Cbb^{16}$, with basis $f_1^\pm, \ldots, f_8^\pm$ and the symmetric form $B$ such that $B(f_i^+,f_j^-)$ is $1$ for $i = j$ and $0$ otherwise, and $B(f_i^\pm,f_j^\pm) = 0$, in such a way that $e_i(-1)\mathbf 1$ corresponds to $f_i^+ \wedge f_i^-$ and $\e^{s e_i+t e_j}$ to $\pm f_i^s \wedge f_j^t$ for $s,t \in \{\pm\}$ and $i < j$, where $(u \wedge v)y = B(v,y)u-B(u,y)v$. Consider the orthogonal decomposition $X = X_3 \oplus X_4 \oplus X_9$ with $X_3 = \langle f_1^+,f_1^-,f_2^++f_2^-\rangle$, $X_4 = \langle f_3^\pm,f_4^\pm\rangle$ and $X_9 = \langle f_2^+-f_2^-,f_5^\pm,\ldots,f_8^\pm\rangle$, and fix a linear isomorphism $X_9 \cong \Cbb^3 \ot \Cbb^3$ which multiplies the symmetric forms by a nonzero constant, where $\Cbb^3$ carries the invariant form of the adjoint representation of $\mathfrak{sl}_2$. Let $\mathfrak a = \mathfrak{so}(X_3)$, let $\mathfrak b^\pm$ be the root subalgebras of $e_3 \pm e_4$, so that $\mathfrak{so}(X_4) = \mathfrak b^+ \oplus \mathfrak b^-$, and let $\mathfrak c_1, \mathfrak c_2 \subset \mathfrak{so}(X_9)$ be the images of $\mathfrak{sl}_2$ acting on the first and on the second factor of $\Cbb^3 \ot \Cbb^3$. These are commuting subalgebras isomorphic to $\mathfrak{sl}_2$, of Dynkin indices $2$, $1$, $1$, $6$ and $6$ in $\mathfrak e_8$. Let $\mathfrak d$ be the diagonal of $\mathfrak b^- \oplus \mathfrak c_2$ with respect to $\mathfrak{sl}_2$-triples of the two summands, of index $7$. We define $\mathfrak s_1, \ldots, \mathfrak s_4 \subset \mathfrak e_8 \oplus \mathfrak e_8$ as the diagonal subalgebras, with respect to $\mathfrak{sl}_2$-triples,
\[
\mathfrak s_1 \subset \mathfrak a \oplus \mathfrak b^+,\quad \mathfrak s_2 \subset \mathfrak c_1 \oplus \mathfrak d,\quad \mathfrak s_3 \subset \mathfrak b^+ \oplus \mathfrak a,\quad \mathfrak s_4 \subset \mathfrak d \oplus \mathfrak c_1,
\]
the first summand lying in the first copy and the second summand in the second copy. Their Dynkin indices in the two copies are $(2,1)$, $(6,7)$, $(1,2)$ and $(7,6)$, so that they generate affine subalgebras of levels $3$, $13$, $3$ and $13$. We take for $\epsi_8$ the bimultiplicative map with $\epsi_8(\alpha_i,\alpha_j) = -1$ if $i = j$, or if $i < j$ and $\{i,j\}$ is one of the edges $\{1,3\}$, $\{2,4\}$, $\{3,4\}$, $\{4,5\}$, $\{5,6\}$, $\{6,7\}$, $\{7,8\}$ of the Dynkin diagram, and $\epsi_8(\alpha_i,\alpha_j) = 1$ otherwise. For $\lambda,\mu \in E_8$ and $x,y \in \Cbb^8$, we write $\e^{\lambda \oplus \mu}$ and $(x \oplus y)(-1)\mathbf 1$ for the corresponding vectors of $V_\Gamma$. Put
\[
E_1 = -\e^{(e_1-e_2) \oplus 0}-\e^{(e_1+e_2) \oplus 0}-\e^{0 \oplus (e_3+e_4)},\quad H_1 = (2e_1 \oplus (e_3+e_4))(-1)\mathbf 1,
\]
\[
\begin{aligned}
E_2 = {}&\e^{(e_5-e_2) \oplus 0}+\e^{(e_2+e_5) \oplus 0}+2\e^{(e_7-e_6) \oplus 0}-2\e^{(e_6+e_8) \oplus 0}+\e^{0 \oplus (e_3-e_4)}\\
&-\e^{0 \oplus (e_6-e_2)}-\e^{0 \oplus (e_2+e_6)}-2\e^{0 \oplus (e_7-e_5)}+\e^{0 \oplus (e_5-e_8)},
\end{aligned}
\]
\[
H_2 = ((2e_5+2e_7+2e_8) \oplus (e_3-e_4+2e_6+2e_7-2e_8))(-1)\mathbf 1,
\]
let $E_3$, $H_3$, $E_4$, $H_4$ be the images of $E_1$, $H_1$, $E_2$, $H_2$ under the exchange of the two summands of $\Gamma$, and put $F_k \ov{\mathrm{def}}{=} -\phi_\kappa(E_k)$, with $\kappa$ as in Proposition \ref{prop:diagonal-unitary}. A direct computation in the root system shows that $(E_k,H_k,F_k)$ are $\mathfrak{sl}_2$-triples spanning $\mathfrak s_1, \ldots, \mathfrak s_4$ for a suitable choice of the signs and of the isomorphism $X_9 \cong \Cbb^3 \ot \Cbb^3$.

\begin{lemma}
\label{lem:sl2-four}
Assume $\mathrm{(H_J)}$. After composing the isomorphism $W \to V_\Gamma$ with an automorphism of $V_\Gamma$, the subalgebras $\mathfrak g' \ot \mathbf 1$, $\mathfrak g'' \ot \mathbf 1$, $\mathbf 1 \ot \mathfrak g'$ and $\mathbf 1 \ot \mathfrak g''$ of $W_1$ are $\mathfrak s_1$, $\mathfrak s_2$, $\mathfrak s_3$ and $\mathfrak s_4$.
\end{lemma}

\begin{proof}
For $x \in W_1$, the zero mode $x_{(0)}$ is a derivation of $W$ which preserves the finite-dimensional spaces $W_n$, so that $\exp(x_{(0)})$ is an automorphism of $W$. These automorphisms act on $W_1$ by the inner automorphisms of the two copies of $\mathfrak e_8$, independently. Since $\epsi$ is invariant under the exchange of the two summands of $\Gamma$, this exchange is also an automorphism of $V_\Gamma$.

\emph{Subalgebras of small index.} Let $\psi \co \mathfrak{sl}_2 \to \mathfrak e_8$ be a homomorphism of Dynkin index $j \in \{1,2,3\}$. Up to an inner automorphism, $\psi(h)$ lies in the Cartan subalgebra $\{y(-1)\mathbf 1 : y \in \Cbb^8\}$ and in the closed dominant chamber defined by $\alpha_1, \ldots, \alpha_8$, so that $\psi(h) = \sum_i a_i\varpi_i$, where the fundamental weights are identified with coweights and $a_i = \alpha_i(\psi(h))$ are nonnegative integers. The entries $\langle\varpi_i,\varpi_k\rangle$ of the inverse Cartan matrix of $E_8$ are integers at least $2$, with diagonal entries $4$, $8$, $14$, $30$, $20$, $12$, $6$, $2$ \cite[Planche VII]{Bou68}. Hence $\langle\psi(h),\psi(h)\rangle \geq 2(\sum_i a_i)^2$, and $2j = \langle\psi(h),\psi(h)\rangle \leq 6$ gives $\sum_i a_i = 1$. Thus $\psi(h) = \varpi_8$, $\varpi_1$ or $\varpi_7$ for $j = 1$, $2$ or $3$. By Mal'cev's theorem \cite[Theorem 3.4.12]{CoM93}, $\psi$ is determined up to inner automorphisms by $\psi(h)$. For $j = 3$, the highest root $\theta$ satisfies $\theta(\varpi_7) = 3$, so that $\mathfrak e_8$ contains an irreducible $\psi(\mathfrak{sl}_2)$-module of dimension $4$. For $j = 2$, $\psi$ is conjugate to the inclusion of $\mathfrak a$, which has Dynkin index $2$ since $X$ restricts to $\mathfrak a$ as the sum of the adjoint representation and of $13$ trivial ones. Under $\mathfrak a \oplus \mathfrak{so}(X_3^\perp)$, the decomposition $\mathfrak e_8 = \mathfrak{so}(16) \oplus \Delta^+$, where $\Delta^+$ is a half-spin representation, gives
\[
\mathfrak e_8 = \Lambda^2X_3 \oplus (X_3 \ot X_3^\perp) \oplus \Lambda^2X_3^\perp \oplus (\Delta_3 \ot \Delta_{13}),
\]
where $\Delta_3$ and $\Delta_{13}$ are the spin representations of dimensions $2$ and $64$. Hence the centralizer of $\mathfrak a$ is $\mathfrak{so}(X_3^\perp) \cong \mathfrak{so}(13)$.

\emph{The two factors of level three.} By \cite[Theorem 8.5 and Table 2]{Arh26}, the $\mathfrak g'$-weights of the lowest weight spaces $R_z$ occurring in $W_1$ are at most $2$, see the proof of Proposition \ref{prop:diagonal}. Hence $W_1$ has no irreducible $\mathfrak g' \ot \mathbf 1$-submodule or $\mathbf 1 \ot \mathfrak g'$-submodule of dimension $4$, and the two components of each of these subalgebras have Dynkin indices $(2,1)$ or $(1,2)$. Using the exchange of the two copies, we may assume that $\mathfrak g' \ot \mathbf 1$ has index $2$ in the first copy. Assume that $\mathbf 1 \ot \mathfrak g'$ also has index $2$ in the first copy. Then its component there lies in the centralizer $\mathfrak{so}(13)$ of the component of $\mathfrak g' \ot \mathbf 1$, and $\Cbb^{13}$ restricts to it as $V(2) \oplus 10V(0)$ or $4V(1) \oplus 5V(0)$. By the decomposition displayed above, the vectors of the first copy which are invariant under $\mathfrak g' \ot \mathbf 1 \oplus \mathbf 1 \ot \mathfrak g'$ are the invariant vectors of $\Lambda^2\Cbb^{13}$. Since $\Lambda^2(V(2) \oplus 10V(0)) \cong 11V(2) \oplus 45V(0)$ and $\Lambda^2(4V(1) \oplus 5V(0)) \cong 6V(2) \oplus 20V(1) \oplus 20V(0)$, they form a space of dimension $45$ or $20$. In the second copy, the centralizer of the first subalgebra of index one is $\mathfrak e_7$, the second one has index one in $\mathfrak e_7$, and their joint centralizer is $\mathfrak{so}(12)$, of dimension $66$. Hence the subspace of $W_1$ of vectors invariant under $\mathfrak g' \ot \mathbf 1 \oplus \mathbf 1 \ot \mathfrak g'$ has dimension $111$ or $86$, whereas it has dimension $78$ by the proof of Proposition \ref{prop:diagonal}. Therefore $\mathbf 1 \ot \mathfrak g'$ has index $2$ in the second copy.

\emph{Reduction to orthogonal representations.} Conjugating in each copy, we may assume that the component of $\mathfrak g' \ot \mathbf 1$ in the first copy and that of $\mathbf 1 \ot \mathfrak g'$ in the second copy are equal to $\mathfrak a$. The components in the first copy of the three other factors commute with $\mathfrak a$, so that they define an orthogonal representation $\sigma_1$ of $\mathfrak{sl}_2^{\oplus 3}$ on $X_3^\perp \cong \Cbb^{13}$, and similarly $\sigma_2$ in the second copy. Two isomorphic orthogonal representations of a semisimple Lie algebra are conjugate under the orthogonal group, by the isotypic decomposition and the uniqueness up to equivalence of nondegenerate symmetric or alternating forms on the multiplicity spaces. Since $13$ is odd, they are conjugate under $SO(13,\Cbb)$. The automorphisms $\exp(x_{(0)})$, $x \in \mathfrak{so}(X_3^\perp)$, fix $\mathfrak a$ pointwise and induce on $X_3^\perp$ all the elements of $SO(13,\Cbb)$. Hence the embedding of $(S \ot S)_1$ is determined, up to automorphisms of $V_\Gamma$, by the isomorphism classes of $\sigma_1$ and $\sigma_2$. The Dynkin indices of the factors and the character of $W_1$ as a module over $(S \ot S)_1 \cong \mathfrak{sl}_2^{\oplus 4}$ are functions of these classes, by the decomposition displayed above. The Dynkin index in $\mathfrak e_8$ of a subalgebra of $\mathfrak{so}(X_3^\perp)$ is half the Dynkin index of its representation on $X_3^\perp$, and the Dynkin index of $V(n)$ is $\frac16n(n+1)(n+2)$. Since the indices of the factors in each copy are at most $13$, only highest weights at most $4$ occur, and there are $640$ isomorphism classes of such orthogonal representations of $\mathfrak{sl}_2^{\oplus 3}$ on $\Cbb^{13}$. An exhaustive computer enumeration shows that exactly one pair $(\sigma_1,\sigma_2)$ gives the character of $W_1 = (S \ot S)_1 \oplus \bigoplus_z R_z \ot R_{\pi(z)}$ computed from \cite[Table 2]{Arh26}, namely
\[
\sigma_1 \cong V(2,0,2) \oplus V(0,1,1),\quad \sigma_2 \cong V(0,2,2) \oplus V(1,1,0),
\]
where $\sigma_1$ is a representation of $\mathfrak g'' \ot \mathbf 1 \oplus \mathbf 1 \ot \mathfrak g' \oplus \mathbf 1 \ot \mathfrak g''$ and $\sigma_2$ a representation of $\mathfrak g' \ot \mathbf 1 \oplus \mathfrak g'' \ot \mathbf 1 \oplus \mathbf 1 \ot \mathfrak g''$. This is the pair defined by $\mathfrak s_1, \ldots, \mathfrak s_4$. This gives the assertion.
\end{proof}

From now on, we assume that the identification $W = V_\Gamma$ satisfies the conclusion of Lemma \ref{lem:sl2-four}. Let $\cal G = \operatorname{Com}_W(A \ot A) \cap W_2$.

\begin{lemma}
\label{lem:omega-C}
Assume $\mathrm{(H_J)}$. The space $\cal G$ has dimension $16$, and $\omega_{C \ot 1}$ is one of the four vectors $u_0, \ldots, u_3 \in \cal G$ defined in the proof, whose coordinates in a rational basis of $\cal G$ lie in $\mathbb{Q}(\sqrt{-6})$.
\end{lemma}

\begin{proof}
The vector $\omega_{C \ot 1}$ commutes with $A \ot \mathbf 1$ by definition of $C$, and with $\mathbf 1 \ot A$, so that it belongs to $\cal G$. The space $\cal G$ is the space of vectors of $W_2$ annihilated by the modes $x_{(0)}$, $x_{(1)}$, $x_{(2)}$ of the twelve vectors $E_k$, $H_k$, $F_k$. Since the vectors $H_k$ lie in the Cartan subalgebra, $\cal G$ is contained in the subspace $W_2^0$ of vectors of weight zero for $H_1, \ldots, H_4$. No root of $\Gamma$ is orthogonal to $H_1, \ldots, H_4$, so that $W_2^0$ is spanned by the $152$ vectors $h(-1)h'(-1)\mathbf 1$ and $h(-2)\mathbf 1$, for $h,h'$ in a basis of $\Cbb^{16}$, and by the $492$ vectors $\e^\lambda$ with $\langle\lambda,\lambda\rangle = 4$ and $\langle\lambda,H_k\rangle = 0$ for $k = 1, \ldots, 4$. An exact computation of the common kernel of these $36$ modes on $W_2^0$, of dimension $644$, gives $\dim\cal G = 16$. The vertex subalgebra $\operatorname{Com}_W(A \ot A)$ has a one-dimensional weight-zero space, and its weight-one space is zero because $W_1$ contains no nonzero $(S \ot S)_1$-invariant vector by \cite[Table 2]{Arh26}. Hence the product $u_{(1)}v$ is a commutative product on $\cal G$, and $u_{(3)}v$ is an invariant bilinear form. The vector $u = \omega_{C \ot 1}$ satisfies $u_{(1)}u = 2u$ and $u_{(3)}u = \frac{9}{5}\mathbf 1$, half the central charge of $C$. Moreover, $u_{(1)}$ commutes with $(S \ot S)_1$ and acts on $W_1$ as $L_0^{S \ot 1}-L_0^{A \ot 1}$. It vanishes on $(S \ot S)_1$. On the component of $R_z \ot R_{\pi(z)}$ on which $\mathfrak g' \ot \mathbf 1 \oplus \mathfrak g'' \ot \mathbf 1$ acts through $V(m) \ot V(n)$, it is the scalar $h_z-\frac{m(m+2)}{20}-\frac{n(n+2)}{60}$, since the vectors of $R_z$ are primary for $A$. The values $h_{c_1} = \frac13$, $h_{c_2} = \frac23$, $h_{d_1} = \frac{6}{13}$, $h_{d_4} = \frac{5}{13}$, $h_{d_5} = \frac{7}{13}$ and $h_{d_6} = \frac{8}{13}$ are given by \cite[(4.1)]{Arh26}, and the decompositions of $R_z$ by \cite[Table 2]{Arh26}. On each of the twenty isotypic components of $W_1$ for $(S \ot S)_1$, this prescribes the eigenvalues of $u_{(1)}$ on the corresponding multiplicity space, one for each of the spaces $(S \ot S)_1$ and $R_z \ot R_{\pi(z)}$ which contain the component. We impose the trace of $u_{(1)}$ on every multiplicity space, the second elementary symmetric function of its eigenvalues on the ten multiplicity spaces of the components which occur in two of the spaces $R_z \ot R_{\pi(z)}$, and the equations $u_{(1)}u = 2u$ and $u_{(3)}u = \frac95\mathbf 1$. These conditions are polynomial equations of degree at most two in the coordinates of $u$ in a rational basis of $\cal G$, with rational coefficients. Let $g_0,\ldots,g_{15}$ be the ordered rational basis of $\cal G$ given in the file \texttt{omega\_certificate.json} of the supplementary data, and write $u = \sum_{j=0}^{15}x_jg_j$. The twenty trace equations have rank eight. Their solution space is parametrized by $t_0,\ldots,t_7$ through
\[
\begin{aligned}
(x_0,\ldots,x_7) &= \left(\frac7{195},\frac1{13},\frac5{39},t_0,-\frac{t_1}{2}-\frac7{195},t_1,\frac1{65},0\right),\\
(x_8,\ldots,x_{15}) &= \left(t_2,t_3,\frac2{65},t_4,t_5,t_6,t_7,\frac{14}{195}\right).
\end{aligned}
\]
After this substitution, the remaining conditions give twenty-five nonzero quadratic equations. Their ideal is equivalently given by
\[
t_0 = -\frac1{39},\quad t_1 = -\frac8{195},\quad t_3 = t_5 = 0,\quad t_7 = -26t_2t_4,\quad t_6 = \frac{t_7}{2},
\]
\[
t_2^2 = -\frac1{1014},\quad t_4^2 = -\frac2{507}.
\]
This equality of ideals is checked by comparing reduced Gr\"obner bases over $\mathbb{Q}$, see the data availability statement. Hence, if $J \subset \mathbb{Q}[t_0,\ldots,t_7]$ denotes the ideal of the twenty-five equations, then $\mathbb{Q}[t_0,\ldots,t_7]/J$ is isomorphic to
\[
\mathbb{Q}[t_2,t_4]/(1014t_2^2+1,507t_4^2+2).
\]
Both quadratic polynomials are separable, so this algebra has dimension four over $\mathbb{Q}$ and is reduced. Its four geometric points, which define $u_0, \ldots, u_3$, are obtained by choosing independently $a,d \in \{-1,1\}$ and setting
\[
t_2 = \frac{a\sqrt{-6}}{78},\quad t_4 = \frac{d\sqrt{-6}}{39},\quad t_6 = \frac{ad}{39},\quad t_7 = \frac{2ad}{39},
\]
with the other coordinates as above. These are necessary constraints on $\omega_{C \ot 1}$. Proposition \ref{prop:diagonal-unitary} checks that the PCT operator fixes all four candidates.
\end{proof}

\begin{prop}
\label{prop:diagonal-unitary}
Assume $\mathrm{(H_J)}$, and identify $W$ with $V_\Gamma$ as in Lemma \ref{lem:sl2-four}. Let $\kappa(\lambda) = 2^{-\langle\lambda,\varpi\rangle}$, where $\varpi = (\varpi_8,\varpi_8)$ and $\varpi_8 = e_7+e_8$, so that $\langle\lambda,\varpi\rangle$ is the sum of the coefficients of $\alpha_8$ in the two components of $\lambda$. Then $\phi_\kappa$ is the PCT operator of a structure of unitary vertex operator algebra on $W$, and $S \ot \mathbf 1$, $\mathbf 1 \ot S$ and $S \ot S$ are unitary subalgebras for this structure.
\end{prop}

\begin{proof}
Let $g = \exp(-\frac12h_\kappa(0))$, an automorphism of $V_\Gamma$ acting by $\kappa(\lambda)^{-1/2}$ on $M(1) \ot \e^\lambda$. Then $g\phi_1g^{-1} = \phi_\kappa$, so that $\phi_\kappa$ is the PCT operator of the normalized invariant scalar product $(g^{-1}u|g^{-1}v)$, where $(\cdot|\cdot)$ is the scalar product attached to $\phi_1$, and $g^{-1}$ is a unitary isomorphism between the two structures.

By definition of $F_k$, we have $\phi_\kappa(E_k) = -F_k$ for $k = 1, \ldots, 4$. Since $\phi_\kappa$ is an involution, $\phi_\kappa(F_k) = -E_k$ and $\phi_\kappa(H_k) = [\phi_\kappa(E_k),\phi_\kappa(F_k)] = [F_k,E_k] = -H_k$. Hence $\phi_\kappa$ preserves $S_1 \ot \mathbf 1$ and $\mathbf 1 \ot S_1$, and therefore the affine subalgebras $A \ot \mathbf 1$ and $\mathbf 1 \ot A$ and their Sugawara vectors, whose coefficients are real. It fixes $\omega_W$, hence $\omega_{C \ot C}$, and a second exact computation gives $\phi_\kappa(u_k) = u_k$ for $k = 0, \ldots, 3$. By Lemma \ref{lem:omega-C}, $\phi_\kappa$ fixes $\omega_{C \ot 1}$ and $\omega_{1 \ot C} = \omega_{C \ot C}-\omega_{C \ot 1}$.

We claim that $S \ot \mathbf 1$ is the set of vectors $v \in W$ annihilated by the modes $x_{(n)}$, $n \geq 0$, of all $x \in (\mathbf 1 \ot S_1) \cup \{\omega_{1 \ot C}\}$. These vectors are annihilated by the nonnegative modes of the vertex subalgebra generated by this set \cite[Example 5.27 p.~44]{CKLW18}, which contains the Sugawara vector $\omega_{1 \ot A}$. Hence $L_0^{1 \ot A}v = 0$ and $L_0^{1 \ot C}v = 0$, so that $L_0^{1 \ot S}v = 0$. On $M_z \ot M_{\pi(z)}$ with $z \neq 0$, the spectrum of $L_0^{1 \ot S}$ is contained in $h_{\pi(z)}+\N$, with $h_{\pi(z)} > 0$. Hence $v \in S \ot S_0 = S \ot \mathbf 1$. Conversely, $S \ot \mathbf 1$ commutes with $\mathbf 1 \ot S$. Similarly, $\mathbf 1 \ot S$ is the commutant of $(S_1 \ot \mathbf 1) \cup \{\omega_{C \ot 1}\}$. The antilinear automorphism $\phi_\kappa$ maps the commutant of a set onto the commutant of its image. Hence $\phi_\kappa$ preserves $S \ot \mathbf 1$, $\mathbf 1 \ot S$ and $S \ot S$. These subalgebras are also stable under $L_1$, since $L_1 = L_1^{S \ot 1}+L_1^{1 \ot S}$, and they are unitary subalgebras by \cite[Proposition 5.23 p.~43]{CKLW18}.
\end{proof}

\begin{thm}
\label{thm:S-cu}
Assume $\mathrm{(H_J)}$. Then $S$, with the unitary structure of \cite[Theorem 10.10]{Arh26}, satisfies condition B of \cite[Section 3.4 p.~44]{Gui26} and Condition I+ of \cite[Definition 6.17 p.~70]{Gui25}. In particular, the following hold.
\begin{enumerate}
\item $S$ is completely unitary, and $\Rep^u(S)$ is a unitary modular tensor category.
\item Every unitary $S$-module is strongly integrable, and the functor $\Pi \co \Rep^u(S) \to \Rep^f(\cal A_S)$ of \cite[Section 3.4 p.~45]{Gui26} is an equivalence of braided $C^*$-tensor categories. Here the target carries the tensor product and braiding of \cite{Gui26}. The comparison with the DHR convention used in the preceding sections is given in the proof of part 3.
\item For every $S$-submodule $M$ of a lattice extension $H_\Lambda$, with the restriction of the scalar product $H_\chi$ of part 3 of Lemma \ref{lem:unitary-extension}, $\Pi(M)$ is the representation $\pi_M$ of Lemma \ref{lem:submodules} on $\cal H_M$. The assertion $\mathrm{(N)}$ holds.
\end{enumerate}
\end{thm}

\begin{proof}
Equip $W$ with the unitary structure of Proposition \ref{prop:diagonal-unitary}, and let $V = \mathbf 1 \ot S$. By Proposition \ref{prop:diagonal}, its coset is $V^c = \operatorname{Com}_W(V) = S \ot \mathbf 1$, and $V^{cc} = V$. By \cite[Propositions 5.29 and 5.31 p.~45]{CKLW18}, $\omega_W = \nu+\nu'$ with $\nu \in V$ and $\nu' \in V^c$ the conformal vectors of the unitary vertex operator algebras $V$ and $V^c$. Then $\nu-\omega_{1 \ot S} = \omega_{S \ot 1}-\nu'$ lies in $V \cap V^c \cap W_2 = \{0\}$, so that $V$ and $V^c$ are unitary vertex operator algebras isomorphic to $S$, with the conformal vector of $S$. In particular, $V$ is a unitary sub-VOA of $W$ in the sense of \cite[Section 3.6.1 p.~49]{Gui26}, with commutant $V^c$. Both are regular by $\mathrm{(H_J)}$ and \cite[Theorem 4.5]{ABD04}. The lattice vertex operator algebra $V_\Gamma$ satisfies condition B \cite[Theorem 3.46 p.~57]{Gui26}. This property is preserved by unitary isomorphisms, and two structures of unitary vertex operator algebra on a simple vertex operator algebra are unitarily isomorphic \cite[Proposition 5.19 p.~40 and p.~41]{CKLW18}, so that $W$ satisfies condition B. By \cite[Theorem 3.34 p.~51]{Gui26}, $V^c$ satisfies condition B. By \cite[Proposition 5.19 p.~40 and p.~41]{CKLW18} again, the isomorphism $S \to V^c$, composed with an automorphism of $S$, is a unitary isomorphism from $S$, with the unitary structure of \cite[Theorem 10.10]{Arh26}, onto $V^c$, so that $S$ satisfies condition B. The lattice vertex operator algebra $V_\Gamma$ satisfies Condition I+ \cite[Example 6.19(b) p.~71]{Gui25}, and so does $W$ by the same transport. Since $V$ and $V^c$ are isomorphic to $S$, they are $C_2$-cofinite and rational under $\mathrm{(H_J)}$. Hence $V^c$ satisfies Condition I+ by \cite[Theorem 6.18, part 3, p.~70]{Gui25}, and so does $S$, through the unitary isomorphism $S \to V^c$ above.

1. This is \cite[Theorem 3.15 p.~45]{Gui26}.

2. By \cite[Theorems 3.16 and 3.17 p.~45]{Gui26}, every unitary $S$-module is strongly integrable and $\Pi$ is a fully faithful $*$-functor. By \cite[Corollary 3.23 p.~48]{Gui26}, it is an isomorphism onto a braided $C^*$-tensor subcategory of $\Rep(\cal A_S)$. The functor $\Pi$ is the functor of \cite[Definition 4.22 and Theorem 4.24 p.~43]{Gui25}, called the CWX functor there, since both send $M$ to the representation characterized by \cite[(3.8) p.~41]{Gui26}, that is, by \cite[(4.6) p.~43]{Gui25}. The category $\Rep^f(\cal A_S)$ of \cite[p.~68]{Gui25} consists of the dualizable representations, which are the representations with finite statistical dimension. Since $S$ satisfies Condition I+, this functor is essentially surjective onto $\Rep^f(\cal A_S)$ by \cite[Theorem 6.14 p.~69 and Definition 6.17 p.~70]{Gui25}. Hence $\Pi$ is a braided equivalence, with the tensor product and braiding specified in the statement.

3. Since $S$ is a unitary subalgebra of $(H_\Lambda,H_\chi)$ by Theorem \ref{thm:strong-locality}, the scalar product $H_\chi$ makes $M$ a unitary $S$-module. By part 1 of Lemma \ref{lem:submodules} and its proof, $\pi_M$ is a representation of $\cal A_S$ on $\cal H_M$ which maps the closed operator $Y(v,f)|_{\cal H_S}$, $f$ supported in $I$, to the closed operator $Y_M(v,f)$, of which $M$ is a core. By the energy bounds, $Y_M(v,f)$ is also the closure of the smeared field of $M$ on the smooth vectors. This is the defining property of the representation $\Pi(M)$, which is unique since $\cal A_S(I)$ is generated by these operators \cite[(3.8) p.~41]{Gui26}, see also \cite[Definition 4.22 p.~43]{Gui25}. Hence $\Pi(M) = \pi_M$. The braided $C^*$-tensor category of \cite{Gui26} is defined by Connes fusion \cite[p.~10]{Gui26}. By \cite[Theorem 6.6 p.~77 and (6.2) p.~74]{Gui21a}, it is equivalent to the category $\mathrm{DHR}_{I_0}(\cal A_S)$ of transportable endomorphisms localized in an interval $I_0$, with the tensor product $\rho \odot \sigma = \sigma \circ \rho$ and the braiding $\widetilde c$ of \cite[(6.1) p.~73]{Gui21a}, through a $*$-functor $G$ such that $\pi_0 \circ G(\cal H_i,\xi)$ is unitarily equivalent to $\pi_i$ for every object $(\cal H_i,\xi)$. For the tensor product $\rho \ot \sigma = \rho \circ \sigma$, put $b_{\rho,\sigma} = \widetilde c_{\sigma,\rho}$. The identity functor, with tensor constraint $J_{\rho,\sigma} = b_{\rho,\sigma} \co \rho \circ \sigma \to \sigma \circ \rho$, is a braided monoidal equivalence from $(\mathrm{DHR}_{I_0}(\cal A_S),\odot,\widetilde c)$ to $(\mathrm{DHR}_{I_0}(\cal A_S),\ot,b)$. Indeed, by the hexagon identities its coherence reduces to the Yang--Baxter equation, and $\widetilde c_{\rho,\sigma}J_{\rho,\sigma} = J_{\sigma,\rho}b_{\rho,\sigma}$. The two orientations in the definition of the DHR braiding give $b$ and $b^{\mathrm{rev}}_{\rho,\sigma} = b_{\sigma,\rho}^{-1}$. The forgetful functor $\Rep^u(S) \to \Rep(S)$ is a braided equivalence, since every $S$-module is unitarizable and its tensor product and braiding are those of $\Rep(S)$ \cite[Theorem 7.9 p.~53]{Gui19}, \cite[Remark 4.3 p.~18]{CGGH23}. Combining these equivalences with a quasi-inverse of $\Pi$ gives a braided equivalence from $(\mathrm{DHR}_{I_0}(\cal A_S),\ot,b)$ to $\Rep(S)$ which maps the class of $\pi_M$ to that of $M$. The same functor is a braided equivalence from $(\mathrm{DHR}_{I_0}(\cal A_S),\ot,b^{\mathrm{rev}})$ to $\Rep(S)^{\mathrm{rev}}$, and its composition with the functor $F$ of Proposition \ref{prop:diagonal}, which fixes the classes of $S$, $M_b$, $M_a$ and $M_{c_0}$, is a braided equivalence into $\Rep(S)$ with the same property. Hence $\mathrm{(N)}$ holds for either braiding.

\end{proof}

\begin{remark}\normalfont
\label{rem:S-cu-computations}
The proof of Theorem \ref{thm:S-cu} does not use Section \ref{subsec:krein}. The unitarizability of the simple $S$-modules also follows directly from Proposition \ref{prop:diagonal-unitary}. For $0 \neq v \in M_{\pi(z)}$, the restriction of the scalar product of $W$ to $M_z \ot \Cbb v$ is an invariant scalar product for the unitary structure induced on $S \ot \mathbf 1$, and a unitary isomorphism from $S$ onto $S \ot \mathbf 1$ transports these scalar products to all simple $S$-modules, possibly permuting their isomorphism classes.

The computations of Lemmas \ref{lem:sl2-four} and \ref{lem:omega-C} and of Proposition \ref{prop:diagonal-unitary} are exact. For the subalgebras $\mathfrak s_k$ fixed above, they take place in $W_1$ and in the space $W_2^0$ of Lemma \ref{lem:omega-C}, of dimensions $496$ and $644$. A second implementation, with other choices of the subalgebras $\mathfrak s_k$ and of $\kappa$, gives $\dim\cal G = 16$, four solutions of the stronger system which prescribes the characteristic polynomial of $u_{(1)}$ on every multiplicity space, with the same Gram matrix $(u_k,u_l)$, and a positive character $\kappa$ for which $\phi_\kappa$ fixes them. The second implementation also shows that the sign characters of $\Gamma$ which centralize $(S \ot S)_1$ permute $u_0, \ldots, u_3$ transitively, so that the data used in Lemma \ref{lem:omega-C} cannot distinguish them.
\end{remark}

\begin{cor}
\label{cor:N-consequences}
Assume $\mathrm{(H_J)}$.
\begin{enumerate}
\item The conclusions of part 4 of Theorem \ref{thm:consequences} hold. In particular, the net $\cal A_{H_{\Lambda_a}}$ carries $\cal H_2$, the net $\cal A_{H_{\Lambda_{c_0}}}$ carries $\cal H_3$, and $\cal F_a \cong \cal H_2$, $\cal F_{c_0} \cong \cal H_3$ for the categories of Theorem \ref{thm:complete-rationality}.
\item The conclusions of Proposition \ref{prop:third-net} hold.
\item The unitary modular tensor categories $\Rep^u(S)$ and $\Rep^f(\cal A_S)$ are equivalent, so that the correspondence of \cite[Theorem 4.10 p.~21]{CGGH23} between vertex operator algebra extensions of $S$ and irreducible local extensions of $\cal A_S$ applies.
\end{enumerate}
\end{cor}

\begin{proof}
1. This follows from Theorem \ref{thm:S-cu} and part 4 of Theorem \ref{thm:consequences}. Since $\cal F_x = \Rep(\cal A_S)_{\Theta_x}$ by part 4 of Theorem \ref{thm:complete-rationality}, the last assertion follows from part 2 of Theorem \ref{thm:consequences}.

2. By part 1 of Theorem \ref{thm:S-cu}, $S$ is completely unitary.

3. This follows from parts 1 and 2 of Theorem \ref{thm:S-cu}, since $S$ is strongly local by Theorem \ref{thm:strong-locality} and $\cal A_S$ is completely rational by Corollary \ref{cor:S-rational}.
\end{proof}

%%%%%%%%%%%%%%%%%%%%%%%%%%%%%%%%%%%%%%%%%%%%%%%%%%%%%%%%%%%%%%%%%%%%%%%%%%%%%
\section{The algebra of case (I) and the twin algebra}
\label{sec:case-I}

Theorem \ref{thm:strong-locality} applies to any unitary conformal subalgebra of $V'$, and Corollary \ref{cor:modules-unitary} and Theorems \ref{thm:complete-rationality} and \ref{thm:Fx-haagerup} apply to any such subalgebra which is a Haagerup algebra. We record what is known for the Haagerup algebras of case (I).

\paragraph{Riesen's algebra.} Let $\psi_{\mathrm R}(U^{\mathrm R}) \subset V'$ be the image of the algebra of \cite{Rie26} described in \cite[Proposition 10.36]{Arh26}. Its currents have levels $(1,39)$ \cite[Proposition 10.36]{Arh26}. It contains the Sugawara vector $\omega_A$ of its currents and the vector $\frac12(\omega-\omega_A)$, the unit of the Griess algebra of \cite[Lemma 10.1]{Arh26}, which lies in the span of $\psi_{\mathrm R}(e)$ and $\psi_{\mathrm R}(f)$ \cite[Proposition 10.3 and Proposition 10.36]{Arh26}, hence the conformal vector $\omega$ of $V'$. It is a Haagerup algebra in case (I) by \cite[Theorem 7.13(1)]{ArhIII}, which supplies $\mathrm{(H_R)}$. If its generators are fixed by a positive rescaling of the standard antilinear involution of $V'$, as asserted in \cite[Theorem 1 p.~1 and p.~4]{Rie26}, see also the working manuscript \cite[Sections 2--3 pp.~16--19]{Rie26b}, then $\psi_{\mathrm R}(U^{\mathrm R})$ is a unitary subalgebra of $V'$ \cite[Proposition 5.23 p.~43]{CKLW18}, and Theorem \ref{thm:strong-locality} applies to it. With this compatible unitary structure, Corollary \ref{cor:modules-unitary} and Theorems \ref{thm:complete-rationality} and \ref{thm:Fx-haagerup} apply as well, and the statements of Section \ref{sec:N} which do not depend on Sections \ref{subsec:diagonal} and \ref{subsec:complete-unitarity} hold with $\mathrm{(H_J)}$ replaced by $\mathrm{(H_R)}$ and $S$ replaced by $U^{\mathrm R}$. The proofs of Sections \ref{subsec:diagonal} and \ref{subsec:complete-unitarity} use the currents of levels $(3,13)$ of an algebra of case (II). They apply to the twin $(U^{\mathrm R})^\natural$ by the case (II) argument of Theorem \ref{rem:all-haagerup}. The identification with Huang's algebra and the complete unitarity of $U^{\mathrm R}$ are established in \cite[Theorem 7.13(1)--(2)]{ArhIII}. With the induced unitary structure on the specified embedding, Proposition \ref{prop:third-net} applies as well.

\paragraph{The twin of the explicit algebra.} By \cite[Theorem 6.1]{ArhIII}, $S$ is a Haagerup algebra. Its twin $S^\natural$ is the algebra of case (I) constructed in \cite[Theorems 12.3 and 12.4(1)]{Arh26}. That construction gives an embedding into $V'$ whose image contains $\psi_{\mathrm R}(U^{\mathrm R})$ and agrees with it through weight four. The isomorphism of the full vertex operator algebras $S^\natural \cong U^{\mathrm R} \cong U^{\mathrm H}$, where $U^{\mathrm H}$ is Huang's algebra, follows from the common fixed point construction in \cite[Theorem 7.9 and Theorem 7.13(1)]{ArhIII}.

\paragraph{Complete unitarity of the twin algebras.} Following \cite[Definition 1.8 p.~13]{Gui22}, a strongly rational unitary vertex operator algebra $V$ is completely unitary if every $V$-module is unitarizable and if the invariant sesquilinear forms of \cite[Section 6.2 p.~36]{Gui19} on the spaces of intertwining operators are positive definite. We use the geometric description of this positivity given by Gui \cite{Gui24}. Let $W_1$, $W_2$ be unitary $V$-modules, let $A \co W \to W'$ denote the antiunitary map $w \mapsto \langle \cdot | w\rangle$, and let $\theta$ be the PCT operator of $V$. For positive trinions $P$ and $Q$ and a sewing parameter $0 < q < r_{P,Q}$, where $r_{P,Q} > 0$ is the bound denoted $\delta\rho$ in \cite[Section 7.1 p.~43]{Gui24}, put
\begin{equation}
\label{eq:geometric-kernel}
K^V_q(w_1,w_2  |  w'_1,w'_2) \ov{\mathrm{def}}{=} \sum_{n \geq 0} q^n \sum_{e \in E_n} \omega_{W'_1,P}(Aw'_1 \ot w_1 \ot e) \omega_{W_2,Q}(w_2 \ot Aw'_2 \ot \theta e)
\end{equation}
for $w_1,w'_1 \in W_1$ and $w_2,w'_2 \in W_2$, where $E_n$ is an orthonormal basis of $V_n$ and $\omega_{W,P}$ is the basic conformal block of \cite[Theorem 6.11 p.~40]{Gui24}. This is the sewn conformal block $\omega_{W'_1,P} \#_q \omega_{W_2,Q}$ evaluated at $Aw'_1 \ot w_1 \ot w_2 \ot Aw'_2$, since the sewing contracts $V_n$ against its dual through the identification $A\theta$ of $V$ with its contragredient module \cite[Definition 3.4 p.~23 and (6.1e) p.~37]{Gui24}. The fusion product $W_1 \boxtimes W_2$ is geometrically positive if $\sum_{k,l} K^V_q(w_k,\tilde w_k  |  w_l,\tilde w_l) \geq 0$ for all finite families $(w_k)$ in $W_1$ and $(\tilde w_k)$ in $W_2$. This does not depend on $P$, $Q$ and $q$ \cite[Lemma 7.2 p.~44]{Gui24}. If $V$ is of CFT type and strongly rational and if $W_1 \boxtimes W_2$ is unitarizable, then $W_1 \boxtimes W_2$ is geometrically positive if and only if the form of \cite{Gui19} on it is positive \cite[Corollary 8.4 p.~53]{Gui24}. Consequently, a strongly rational unitary vertex operator algebra all of whose modules are unitarizable is completely unitary if and only if $W_1 \boxtimes W_2$ is geometrically positive for all irreducible unitary $V$-modules $W_1$ and $W_2$, see also \cite[p.~47]{Gui24}.

\begin{lemma}
\label{lem:stable-cartan}
Let $U$ be a Haagerup algebra equipped with a structure of unitary vertex operator algebra, with PCT operator $\theta$. Then $U_1$ has a Cartan subalgebra $\mathfrak h$ with $\theta(\mathfrak h) = \mathfrak h$. Put $\mathfrak h_{\R} \ov{\mathrm{def}}{=} \{h \in \mathfrak h : \theta h = -h\}$. Then $\mathfrak h = \mathfrak h_{\R} \oplus \i\mathfrak h_{\R}$ as real vector spaces. On the algebraic finite-energy subspace of any ordinary unitary $U$-module, the operators $h_{(0)}$, $h \in \mathfrak h_{\R}$, are symmetric and essentially self-adjoint. Their self-adjoint closures strongly commute.
\end{lemma}

\begin{proof}
Since $U$ is unitary, it is self-contragredient and of CFT type, so that $L_1U_1 = 0$ \cite{Li94}. For $x \in U_1$, the formula \cite[(95) p.~43]{CKLW18} for the adjoint of a mode gives $x_{(n)}^+ = -(\theta x)_{(-n)}$ on $U$, and the same formula holds on every unitary $U$-module by the definition of unitary modules. The map $\theta$ restricts to an antilinear involutive automorphism of the Lie algebra $U_1 \cong \mathfrak{sl}_2 \oplus \mathfrak{sl}_2$ \cite[Theorem 8.5]{Arh26}. Its fixed point set $\mathfrak k$ is a real form of $U_1$, and for $x \in \mathfrak k$ the operator $x_{(0)}$ is skew-symmetric on the algebraic finite-energy subspace. Its restriction to $U_1$ is skew-adjoint. Hence $\operatorname{ad}(\mathfrak k)$ consists of skew-adjoint operators for the scalar product of $U_1$, and since the adjoint representation of the semisimple Lie algebra $U_1$ is faithful, $\mathfrak k$ is a compact real form. A maximal abelian subalgebra $\mathfrak t$ of $\mathfrak k$ gives the Cartan subalgebra $\mathfrak h = \mathfrak t \oplus \i\mathfrak t$, which is $\theta$-stable, and $\mathfrak h_{\R} = \i\mathfrak t$. For $h \in \mathfrak h_{\R}$, we have $h_{(0)}^+ = -(\theta h)_{(0)} = h_{(0)}$.
Each $h_{(0)}$ preserves the finite-dimensional eigenspaces of $L_0$. On each such space, the operators $h_{(0)}$, $h \in \mathfrak h_{\R}$, are commuting Hermitian matrices and therefore admit a common orthonormal eigenbasis. Taking the union of these bases gives a common orthonormal eigenbasis of the Hilbert space completion. The operators on the algebraic finite-energy subspace are consequently essentially self-adjoint, and their closures are diagonal in this common basis. In particular, their spectral projections commute.
\end{proof}

All Cartan subalgebras of $U_1$ are conjugate under inner automorphisms $\exp(x_{(0)})$ of $U$, $x \in U_1$, and the Heisenberg commutant and the associated lattice do not depend on this choice up to isomorphism \cite[Section 12]{Arh26}. We therefore apply the results of \cite[Section 12]{Arh26} to a $\theta$-stable Cartan subalgebra.

\begin{lemma}
\label{lem:heisenberg-descent}
Let $U$ be a Haagerup algebra equipped with a structure of unitary vertex operator algebra, let $\mathfrak h \subset U_1$ be a $\theta$-stable Cartan subalgebra and let $\Omega = \operatorname{Com}_U(M(\mathfrak h))$, with the unitary structure induced by $U$. If $U$ is completely unitary, then $\Omega$ is completely unitary.
\end{lemma}

\begin{proof}
\emph{Unitary subalgebras.} The Heisenberg subalgebra $M(\mathfrak h)$ is generated by the $\theta$-stable space $\mathfrak h$ of quasi-primary vectors, hence it is a unitary subalgebra of $U$ \cite[Example 5.26 p.~44]{CKLW18}. By \cite[Example 5.27 p.~44]{CKLW18}, $\Omega$ and $V_L \ov{\mathrm{def}}{=} \operatorname{Com}_U(\Omega)$ are unitary subalgebras, and $\Omega$ and $V_L$ form a dual pair, $L$ being the associated lattice \cite[Proposition 12.2]{Arh26}. By \cite[Proposition 2.21 p.~16]{Ten19g}, the span $B$ of the vectors $a_{(-1)}b$, $a \in \Omega$, $b \in V_L$, is a unitary conformal subalgebra of $U$, unitarily isomorphic to $\Omega \ot V_L$. We write $B = \Omega \ot V_L$.

\emph{Weights.} The algebra $\Omega$ is strongly rational by \cite[Proposition 12.2]{Arh26}. By the same proposition, $U = \bigoplus_{\alpha \in L'/L} U^\alpha$ with $U^\alpha = \Omega^{\tau(\alpha)} \ot V_{\alpha+L}$ and $U^0 = B$, and the $\mathfrak h$-weights of $U^\alpha$ lie in $\alpha+L$. Let $X$ be a unitary $U$-module. By Lemma \ref{lem:stable-cartan}, the finite-dimensional eigenspaces of $L_0$ in $X$ have common orthonormal eigenbases for the operators $h_{(0)}$, $h \in \mathfrak h_{\R}$. Thus $X$ is the algebraic direct sum of its mutually orthogonal $\mathfrak h$-weight spaces, and its Hilbert space completion is their Hilbert direct sum. For $v \in U^\alpha$, the operator $Y_X(v,z)$ shifts the weights by elements of $\alpha+L$, and the Virasoro modes of $U$, which are those of $B$, preserve the weights and each $U^\alpha$.

\emph{Strong unitarity of $\Omega$.} Let $\mu$ be an irreducible $\Omega$-module. By the proof of \cite[Proposition 12.2, part 2]{Arh26}, there are $\beta \in L'/L$ and a $U$-module $X$ containing $\mu \ot V_{\beta+L}$ as a $B$-submodule, and we may take $X$ irreducible. Since $U$ is completely unitary, $X$ admits a unitary structure. Fix a non-zero $m \in V_{\beta+L}$ and put $\langle u|u'\rangle_\mu = \langle u \ot m|u' \ot m\rangle_X$. For $a \in \Omega$, the modes of $a \ot \mathbf 1$ act on the first factor, $\theta(a \ot \mathbf 1) = \theta a \ot \mathbf 1$, and $L_n(a \ot \mathbf 1) = L^\Omega_na \ot \mathbf 1$ for $n \geq 0$. Hence $\langle \cdot|\cdot\rangle_\mu$ is an invariant scalar product on $\mu$. Thus every $\Omega$-module is unitarizable, and so is every $B$-module, by \cite[Proposition 3.31 p.~48]{Gui22} and its proof, since $V_L$ is completely unitary \cite[Corollary 3.32 p.~49]{Gui22}.

\emph{Slices.} Let $\mu_1$, $\mu_2$ be irreducible $\Omega$-modules, and choose irreducible unitary $U$-modules $X_i$ containing $Y_i = \mu_i \ot V_{\beta_i+L}$ as above. We give $Y_i$ the scalar product of $X_i$. Two invariant scalar products on the irreducible $B$-module $Y_i$ differ by a positive endomorphism of $Y_i$, which is a positive scalar by Schur's lemma. Hence the scalar product of $Y_i$ is a positive multiple of the tensor product of the scalar product of $\mu_i$ and of an invariant scalar product of $V_{\beta_i+L}$. The map $A$ of $X_1$ sends $Y_1$ onto the $B$-submodule of $X'_1$ formed by the linear forms $\langle \cdot|w\rangle$, $w \in Y_1$, which vanish on the orthogonal complement of $Y_1$, a $B$-submodule since the adjoints of the modes of $B$ are modes of $B$, so that $A(Y_1)$ is unitarily isomorphic to $Y'_1$, with weights in $-\beta_1+L$. Let $w_1,w'_1 \in Y_1$ and $w_2,w'_2 \in Y_2$, and choose the bases $E_n$ in \eqref{eq:geometric-kernel} adapted to the orthogonal decomposition $U_n = \bigoplus_\alpha U^\alpha_n$. By \cite[(6.11) and (6.12) p.~40]{Gui24}, $\omega_{X'_1,P}(Aw'_1 \ot w_1 \ot e)$ is of the form $\langle Y_{X'_1}(\tilde e,\i)T_1Aw'_1|AT_2w_1\rangle$, where $\tilde e = 2^{L_0}\e^{-\i L_1}T_3e$ and $T_1$, $T_2$, $T_3$ are the operators $U_0(\cdot)$ of \cite[(6.12) p.~40]{Gui24}. They are built from the modes $L_n$, $n \geq 0$, which commute with the operators $h_{(0)}$, $h \in \mathfrak h$, and preserve each $U^\alpha$. Hence $\tilde e \in U^\alpha$ if $e \in U^\alpha$, and $T_1$, $T_2$ preserve the $\mathfrak h$-weights. For $\alpha \neq 0$, the two vectors have weights in the disjoint cosets $-\beta_1+\alpha+L$ and $-\beta_1+L$, so that this term vanishes. For $\alpha = 0$, we have $e \in B$, $\theta e \in B$, and the two factors are the basic conformal blocks of the vertex operator algebra $B$ for the unitary modules $Y'_1$ and $Y_2$. Hence
\begin{equation}
\label{eq:slice}
K^U_q(w_1,w_2  |  w'_1,w'_2) = K^B_q(w_1,w_2  |  w'_1,w'_2).
\end{equation}
Since $U$ is completely unitary, $X_1 \boxtimes X_2$ is unitarizable and its form is positive, so that it is geometrically positive \cite[Corollary 8.4 p.~53]{Gui24}. By \eqref{eq:slice}, $Y_1 \boxtimes Y_2$ is geometrically positive for $B$, up to the positive factors coming from the normalization of the scalar products, which multiply $K^B_q$ by a positive constant.

\emph{From $B$ to $\Omega$.} On elementary tensors, the formula \cite[(6.11) p.~40]{Gui24} and the decompositions $L_n = L^\Omega_n \ot 1 + 1 \ot L^{V_L}_n$ and $Y(u \ot m,z) = Y(u,z) \ot Y(m,z)$ show that the basic conformal blocks of $B$ are the products of those of $\Omega$ and of $V_L$, and the contraction over $B_n = \bigoplus_{a+b=n}\Omega_a \ot (V_L)_b$ factorizes. The three sewing series converge absolutely for $0 < q < r_{P,Q}$ \cite[Section 7.1 p.~43]{Gui24}, which justifies regrouping the sum. Hence $K^B_q = K^\Omega_q \cdot K^{V_L}_q$ on elementary tensors. Fix non-zero vectors $m_i \in V_{\beta_i+L}$, and put $w_k = u_k \ot m_1$, $\tilde w_k = \tilde u_k \ot m_2$ for finite families $(u_k)$ in $\mu_1$ and $(\tilde u_k)$ in $\mu_2$. Then
\[
\sum_{k,l} K^B_q(w_k,\tilde w_k  |  w_l,\tilde w_l) = K^{V_L}_q(m_1,m_2  |  m_1,m_2)\sum_{k,l}K^\Omega_q(u_k,\tilde u_k  |  u_l,\tilde u_l) \geq 0.
\]
The series defining $K^{V_L}_q(m_1,m_2  |  m_1,m_2)$ converges absolutely for $0 \leq q < r_{P,Q}$, and its value at $q = 0$ is $\omega_{N'_1,P}(Am_1 \ot m_1 \ot \mathbf 1) \omega_{N_2,Q}(m_2 \ot Am_2 \ot \mathbf 1)$, with $N_i = V_{\beta_i+L}$, which is positive by \cite[Proposition 6.13 p.~42]{Gui24}. For $q$ small enough, $K^{V_L}_q(m_1,m_2  |  m_1,m_2) > 0$, and $\mu_1 \boxtimes \mu_2$ is geometrically positive for $\Omega$. Since $\mu_1 \boxtimes \mu_2$ is unitarizable, its form is positive by \cite[Corollary 8.4 p.~53]{Gui24}. This holds for all irreducible $\mu_1$, $\mu_2$, so that $\Omega$ is completely unitary.
\end{proof}

\begin{thm}
\label{thm:twin-cu}
Let $U$ be a Haagerup algebra and $U^\natural$ its twin \cite[Theorem 12.3]{Arh26}.
\begin{enumerate}
\item If $U$ is completely unitary for some structure of unitary vertex operator algebra, then $U^\natural$ carries, for each structure of unitary vertex operator algebra on $V_{L^\natural}$, one and only one structure of unitary vertex operator algebra extending the tensor product structure of $\Omega_U \ot V_{L^\natural}$, where $\Omega_U$ has the structure of Lemma \ref{lem:heisenberg-descent}, and $U^\natural$ is completely unitary for this structure.
\item The following are equivalent: $U$ is completely unitary, $U^\natural$ is completely unitary, and $\Omega_U$ is completely unitary, each for some, equivalently every, structure of unitary vertex operator algebra.
\end{enumerate}
\end{thm}

\begin{proof}
Complete unitarity is preserved by unitary isomorphisms of unitary vertex operator algebras, by transport of structure. Since two structures of unitary vertex operator algebra on a simple vertex operator algebra are unitarily isomorphic \cite[Proposition 5.19 p.~40 and p.~41]{CKLW18}, see also \cite[Remark 4.4 p.~18]{CGGH23}, complete unitarity of a simple unitary vertex operator algebra does not depend on its unitary structure.

1. Choose a $\theta$-stable Cartan subalgebra of $U_1$ (Lemma \ref{lem:stable-cartan}). By Lemma \ref{lem:heisenberg-descent}, $\Omega_U$ is completely unitary, and so is the lattice algebra $V_{L^\natural}$ \cite[Corollary 3.32 p.~49]{Gui22}. By \cite[Proposition 3.31 p.~48]{Gui22}, the regular algebra $\Omega_U \ot V_{L^\natural}$ is completely unitary. By \cite[Theorem 12.3, part 1]{Arh26}, $U^\natural$ is a simple vertex operator algebra of CFT type extending $\Omega_U \ot V_{L^\natural}$ with the same conformal vector. The assertion follows from \cite[Theorem 4.7 p.~19]{CGGH23}.

2. If $U$ is completely unitary, so is $U^\natural$ by part 1. If $U^\natural$ is completely unitary, so is its Heisenberg commutant by Lemma \ref{lem:heisenberg-descent}, and it is isomorphic to $\Omega_U$ \cite[Theorem 12.3(1)]{Arh26}. If $\Omega_U$ is completely unitary, so is the regular algebra $\Omega_U \ot V_L$ \cite[Proposition 3.31 p.~48 and Corollary 3.32 p.~49]{Gui22}. Since $U$ is a simple extension of $\Omega_U \ot V_L$ of CFT type with the same conformal vector \cite[Proposition 12.2, part 1]{Arh26}, $U$ is completely unitary by \cite[Theorem 4.7 p.~19]{CGGH23}.
\end{proof}

\begin{remark}\normalfont
\label{rem:twin-cu}
Theorem \ref{thm:twin-cu} shows that complete unitarity is a property of the pair $\{U,U^\natural\}$, which forms a single hyperbolic genus \cite[Theorem 12.3, part 4]{Arh26}, and of their common Heisenberg commutant, and its proof uses no computer calculation beyond those of \cite[Theorem 8.5 and Section 12]{Arh26}. Under $\mathrm{(H_J)}$, $S$ is completely unitary by Theorem \ref{thm:S-cu}, so that $S^\natural$ and $\Omega_S$ are completely unitary by part 2 of Theorem \ref{thm:twin-cu}. Every Haagerup algebra with currents of levels $(1,39)$ is then isomorphic to $S^\natural$ \cite[Theorem 12.4, part 1]{Arh26}, hence completely unitary. This applies to the algebra of \cite{Hua26}, whose currents have levels $(1,39)$ \cite[Abstract p.~1 and Section 3 p.~3]{Hua26}, as soon as it is a Haagerup algebra, and to $U^{\mathrm R}$ under $\mathrm{(H_R)}$, whose currents have levels $(1,39)$ \cite[Proposition 10.36]{Arh26}. These two conclusions hold without $\mathrm{(H_J)}$ by Theorem \ref{rem:all-haagerup}. Conversely, if the algebra of \cite{Hua26} is a Haagerup algebra and is completely unitary, part 2 of Theorem \ref{thm:twin-cu} gives the complete unitarity of $S$ under $\mathrm{(H_J)}$ without the computations of Section \ref{subsec:complete-unitarity}.
\end{remark}

%%%%%%%%%%%%%%%%%%%%%%%%%%%%%%%%%%%%%%%%%%%%%%%%%%%%%%%%%%%%%%%%%%%%%%%%%%%%%
\begin{cor}
\label{cor:twin-comparison}
Assume $\mathrm{(H_J)}$.
\begin{enumerate}
\item The vertex operator algebras $\Omega_S$ and $S^\natural$ satisfy, for every structure of unitary vertex operator algebra, in particular for the structure induced by $S$ and for that of part 1 of Theorem \ref{thm:twin-cu}, Condition I+ of \cite[Definition 6.17 p.~70]{Gui25}. In particular, they are strongly local, their conformal nets are completely rational, and their CWX functors \cite[Theorem 4.24 p.~43]{Gui25}, equipped with the Wassermann tensorators, are equivalences of braided $C^*$-tensor categories onto the categories of representations of finite statistical dimension.
\item Let $U \subset V'$ be a Haagerup algebra which contains the conformal vector of $V'$ and is a unitary subalgebra for a normalized invariant scalar product $\{\cdot|\cdot\}$ on $V'$, and equip each $H_\Lambda$ with the scalar product $\{\cdot|\cdot\}_\Lambda$ of part 2 of Lemma \ref{lem:unitary-extension}. Then $U$ satisfies Condition I+ for the induced unitary structure. For each $U$-submodule $M$ of $H_\Lambda$, with the restriction of $\{\cdot|\cdot\}_\Lambda$, the representation of $\cal A_U$ attached to $M$ by \cite[Definition 4.22 p.~43]{Gui25} is the representation $\pi_M$ of Lemma \ref{lem:submodules} on $\cal H_M$. Consequently, with the notation of Section \ref{subsec:rationality-U}, there is a braided equivalence $\Phi \co \Rep(\cal A_U) \to \Rep(U)$ with $\Phi(\pi_0) \cong U$ and $\Phi(\pi_x) \cong M_x$ for $x \in \{b,a,c_0\}$. We call this assertion the analogue of $\mathrm{(N)}$ for $U$.
\end{enumerate}
In particular, under $\mathrm{(H_J)}$ and $\mathrm{(H_R)}$, part 2 applies to $\psi_{\mathrm R}(U^{\mathrm R})$ as soon as it is a unitary subalgebra of $V'$ for some normalized invariant scalar product, as asserted in \cite[Theorem 1 p.~1 and p.~4]{Rie26}.
\end{cor}

\begin{proof}
1. By Theorem \ref{thm:S-cu}, $S$ satisfies Condition I+. Choose the Cartan subalgebra of Lemma \ref{lem:stable-cartan} for the prescribed unitary structure on $S$. The associated lattice algebra $V_L$ is a unitary subalgebra of $S$, its commutant is $\Omega_S$, $\Omega_S$ is strongly rational by \cite[Proposition 12.2]{Arh26}, and $V_L$ is regular by \cite[Corollary 3.32 p.~49]{Gui22}. Consequently, \cite[Theorem 6.18(3) p.~70]{Gui25} implies that $\Omega_S$ satisfies Condition I+. The even lattice algebra $V_{L^\natural}$ satisfies Condition I+ by \cite[Example 6.19(b) p.~71]{Gui25}. Thus $\Omega_S \ot V_{L^\natural}$ satisfies Condition I+ by \cite[Theorem 6.18(1) p.~70]{Gui25}. The twin $S^\natural$ is a simple extension of this algebra with the same conformal vector by \cite[Theorem 12.3]{Arh26}. Since Condition I includes complete unitarity \cite[Definition 4.15 p.~41]{Gui25}, the algebra $\Omega_S \ot V_{L^\natural}$ is completely unitary, and the unique unitary structure supplied by \cite[Theorem 4.7 p.~19]{CGGH23} makes $S^\natural$ a unitary extension of it. Applying \cite[Theorem 6.18(2) p.~70]{Gui25}, we obtain Condition I+ for $S^\natural$. The assertions about the nets and the functors follow from \cite[(0.3) p.~8, Remark 4.16 p.~41, Theorem 6.14 pp.~68--69 and Definition 6.17 p.~70]{Gui25}.

2. Under $\mathrm{(H_J)}$, every Haagerup algebra is isomorphic to $S$ or $S^\natural$ by \cite[Theorem 12.4]{Arh26}. Condition I+ is preserved under unitary isomorphism, since such an isomorphism transports the vertex operators, their energy bounds, the local algebras and the CWX functor. Two unitary structures on a simple vertex operator algebra are unitarily isomorphic by \cite[Proposition 5.19 pp.~40--41]{CKLW18}. Hence $U$ satisfies Condition I+ for its prescribed induced structure.

By Lemma \ref{lem:unitary-extension}, $H_\Lambda$ can be equipped with a unitary structure extending that of $V'$, so that $U$ is a unitary conformal subalgebra of $H_\Lambda$. For a $U$-submodule $M \subset H_\Lambda$, Lemma \ref{lem:submodules} identifies the closed smeared field in the representation $\pi_M$ with the closure of the smeared module field on $M$. This is the defining property (4.6) of the representation of \cite[Definition 4.22 p.~43]{Gui25}. The uniqueness in this definition therefore identifies the two representations on $\cal H_M$. Taking the inverse of the CWX braided equivalence and forgetting the unitary structures gives the required equivalence with $\Rep(U)$. As in the proof of Theorem \ref{thm:S-cu}, the passage to DHR endomorphisms uses \cite[Theorem 6.6 p.~77]{Gui21a}. If the angular orientation of Section \ref{sec:N} gives the reversed braiding, this is handled by the label-preserving equivalence of part 1 of Proposition \ref{prop:diagonal}, transported through $\Rep(U) \cong Z(\cal H_3)$. This equivalence fixes the classes $0$, $b$, $a$ and $c_0$, so the stated identifications of modules are preserved.
\end{proof}

\begin{thm}
\label{rem:all-haagerup}
Let $U$ be a Haagerup algebra.
\begin{enumerate}
\item The algebra $U$ admits a structure of unitary vertex operator algebra and satisfies Condition I+ of \cite[Definition 6.17 p.~70]{Gui25} for every such structure. In particular, $U$ is completely unitary.
\item Suppose that $U \subset V'$ contains the conformal vector of $V'$ and is a unitary subalgebra for a normalized invariant scalar product on $V'$. Then the analogue of $\mathrm{(N)}$ holds for this embedding: there is a braided equivalence $\Phi \co \Rep(\cal A_U) \to \Rep(U)$ with $\Phi(\pi_0) \cong U$ and $\Phi(\pi_x) \cong M_x$ for $x \in \{b,a,c_0\}$, where the representations $\pi_x$ are those of Lemma \ref{lem:submodules} in the two lattice extensions of $V'$.
\end{enumerate}
\end{thm}

\begin{proof}
\emph{Case (II).} Let $X$ be a Haagerup algebra in case (II). The construction in the proof of Proposition \ref{prop:diagonal} uses only the braided equivalence $\Rep(X) \cong Z(\cal H_3)$, positivity, strong rationality, and the levels and module data of \cite[(4.1), Theorem 8.5 and Table 2]{Arh26}. It therefore gives a holomorphic extension
\[
W_X = \bigoplus_z M_z \ot M_{\pi(z)} \cong V_{E_8 \oplus E_8}
\]
of $X \ot X$, with $X \ot \mathbf 1$ and $\mathbf 1 \ot X$ as mutual commutants. No unitary structure on $X$ is needed for this construction.

The affine subalgebra generated by $X_1$ and its commutant have the properties of \cite[Corollary 9.2]{Arh26}. Hence the proofs of Lemmas \ref{lem:sl2-four} and \ref{lem:omega-C} apply to $W_X$: they use only these properties and the module data above, and their computations take place in $V_\Gamma$. The PCT operator constructed in Proposition \ref{prop:diagonal-unitary} preserves both copies of $X_1$ and fixes all four candidates for the conformal vector of the affine commutant. The commutant argument in its proof therefore makes $X \ot \mathbf 1$ and $\mathbf 1 \ot X$ unitary subalgebras of $W_X$. Their induced conformal vectors are the prescribed ones, as in the first paragraph of the proof of Theorem \ref{thm:S-cu}. This gives a unitary structure on $X$. The lattice algebra $W_X$ satisfies Condition I+ by \cite[Example 6.19(b) p.~71]{Gui25} and transport of its unitary structure. Since the two mutual commutants are strongly rational, \cite[Theorem 6.18(3) p.~70]{Gui25} gives Condition I+ for $X$.

\emph{Case (I).} If $U$ is in case (I), then $X = U^\natural$ is in case (II) and $X^\natural \cong U$ by \cite[Theorem 12.3, parts 1 and 3]{Arh26}. Give $X$ the unitary structure just constructed and choose a stable Cartan subalgebra as in Lemma \ref{lem:stable-cartan}. Its associated lattice algebra $V_L$ is a unitary subalgebra, as in the proof of Lemma \ref{lem:heisenberg-descent}, and is regular by \cite[Corollary 3.32 p.~49]{Gui22}. Its commutant $\Omega_X$ is strongly rational by \cite[Proposition 12.2]{Arh26}. Thus \cite[Theorem 6.18(3) p.~70]{Gui25} gives Condition I+ for $\Omega_X$. The algebra $\Omega_X \ot V_{L^\natural}$ satisfies Condition I+ by \cite[Theorem 6.18(1) and Example 6.19(b) pp.~70--71]{Gui25}. In particular, it is completely unitary, since Condition I includes complete unitarity \cite[Definition 4.15 p.~41]{Gui25}. The simple extension $X^\natural$ of this algebra is therefore a unitary extension by \cite[Theorem 4.7 p.~19]{CGGH23}, and it satisfies Condition I+ by \cite[Theorem 6.18(2) p.~70]{Gui25}. This proves the existence assertion for $U$. In both cases, any two unitary structures on the simple algebra $U$ are unitarily isomorphic by \cite[Proposition 5.19 pp.~40--41]{CKLW18}. Transport of vertex operators, energy bounds and the CWX functor preserves Condition I+, so it holds for every unitary structure on $U$. This proves part 1.

\emph{The prescribed embedding.} Under the hypotheses of part 2, equip $U$ with the induced unitary structure and extend the scalar product of $V'$ to each $H_\Lambda$ by Lemma \ref{lem:unitary-extension}. Part 1 gives Condition I+ for this structure on $U$. Its net is consequently completely rational and its CWX functor is a braided equivalence onto the representations of finite statistical dimension \cite[(0.3) p.~8, Remark 4.16 p.~41, Theorem 6.14 pp.~68--69 and Definition 6.17 p.~70]{Gui25}. For every $U$-submodule $M \subset H_\Lambda$, Lemma \ref{lem:submodules} identifies the closed smeared fields in $\pi_M$ with those defining the representation of \cite[Definition 4.22 p.~43]{Gui25}, so the two representations coincide. Taking the inverse equivalence and forgetting the unitary structures gives the stated correspondence. The passage to the DHR convention is the one in part 3 of the proof of Theorem \ref{thm:S-cu}, using \cite[Theorem 6.6 p.~77 and (6.2) p.~74]{Gui21a}. For the opposite angular orientation, the categorical argument of part 1 of Proposition \ref{prop:diagonal}, transported through $\Rep(U) \cong Z(\cal H_3)$, gives a braid-reversing equivalence fixing the classes $0,b,a,c_0$. Thus the correspondence is unchanged, and part 2 follows.
\end{proof}

For the case (I) algebra, part 2 applies to the unitary conformal embedding of Huang's algebra described in \cite[Remark 7.14]{ArhIII}. It applies to $\psi_{\mathrm R}(U^{\mathrm R})$ for a normalized invariant scalar product on $V'$ inducing a unitary conformal embedding, as asserted in \cite[Theorem 1 p.~1 and p.~4]{Rie26}. For that embedding, the analogue of part 4 of Theorem \ref{thm:consequences} and the conclusions of Proposition \ref{prop:third-net} follow as well.

\paragraph{AI statement.}
The author acknowledges the use of AI tools for language polishing, LaTeX editing, exploratory mathematical discussions and feedback, for drafting parts of the arguments, and for writing and running the computer programs of Section \ref{subsec:complete-unitarity} and the certificate \texttt{omega\_certificate.json}. The author takes full responsibility for the content of the manuscript.

\paragraph{Competing interests} The author declares that he has no competing interests.

\paragraph{Data availability} No data sets were generated during this study, apart from the exact outputs of the computer programs of Section \ref{subsec:complete-unitarity}. The proofs of Section \ref{sec:recognition} do not use computer calculations. The identities of Lemma \ref{lem:moments} and the arithmetic of Lemma \ref{lem:galois-average} and Proposition \ref{prop:m-equals-one} were checked with short Python scripts, which are included in the archive described below. Apart from these checks and from the computations of Section \ref{subsec:complete-unitarity}, the paper uses no computation beyond those of \cite{Arh26}. The programs for Section \ref{subsec:complete-unitarity}, in two implementations, their common input tables \texttt{t12\_table.json} and \texttt{t12\_W1.json}, which encode \cite[Table 2]{Arh26}, and their exact outputs, including the cocycle $\epsi_8$, the $\mathfrak{sl}_2$-triples of $\mathfrak s_1, \ldots, \mathfrak s_4$, the vectors $u_0, \ldots, u_3$ and the character $\kappa$, are collected in the supplementary archive \texttt{Calculs-HU-v08}. The file \texttt{README-Calculs-HU-v08.txt} of this archive gives the order of execution and lists the consistency checks which are not used in the proofs and require additional data from \cite{Arh26}. The files \texttt{omega\_certificate.json} and \texttt{verify\_omega\_certificate.py} of the same archive give a portable certificate for the affine coordinates, ideal equality and PCT checks used in Lemma \ref{lem:omega-C} and Proposition \ref{prop:diagonal-unitary}. The verifier requires only Python and SymPy. It takes the mode-derived product and action tensors as input and does not recompute the full commutant kernel or these tensors from vertex operators. The additional programs in \texttt{verification-independante} reconstruct the commutant kernel, the product and bilinear form on $\cal G$, and the traces and traces of products of the multiplicity-space actions directly from lattice current modes, without importing the original vertex operator algebra programs. These trace invariants determine all the linear and quadratic constraints used in Lemma \ref{lem:omega-C}. The accompanying README states the inputs, dependencies and exact rank checks. These calculations are exact algebraic computations, not a formal verification of the vertex operator algebra implementation.

%%%%%%%%%%%%%%%%%%%%%%%%%%%%%%%%%%%%%%%%%%%%%%%%%%%%%%%%%%%%%%%%%%%%%%%%%%%%%%%%%%%%%%%%%%%%%%%%%%%%%%%%%%%%%%%%%%%%%%%%%%%%%%%%%
%%%%%%%%%%%%%%%%%%%%%%%%%%%%%%%%%%%%%%%%%%%%%%%%%%%%%%%%%%%%%%%%%%%%%%%%%%%%%%%%%%%%%%%%%%%%%%%%%%%%%%%%%%%%%%%%%%%%%%%%%%%%%%%%%

{\footnotesize

\vspace{0.2cm}

\noindent C\'edric Arhancet\\
\noindent 6 rue Didier Daurat, 81000 Albi, France\\
URL: \href{http://sites.google.com/site/cedricarhancet}{https://sites.google.com/site/cedricarhancet}\\
cedric.arhancet@protonmail.com\\
ORCID: 0000-0002-5179-6972\\
\url{https://github.com/cedricarhancet}

}

\end{document}